\documentclass[11pt]{article}
\usepackage[english]{babel}
\usepackage{amsmath,amsthm,amssymb}
\usepackage{mathrsfs}
\usepackage{amsfonts}
\usepackage[margin=0.8in]{geometry}
\usepackage{graphicx}
\usepackage{caption}
\usepackage{subcaption}
\usepackage{wasysym}
\usepackage{enumerate}
\usepackage{algorithm2e}
\usepackage{tabularx}
\usepackage{color}

\newtheorem{thm}{Theorem}[section]
\newtheorem{ass}[thm]{Assumption}
\newtheorem{cor}[thm]{Corollary}
\newtheorem{lem}[thm]{Lemma}
\newtheorem{prop}[thm]{Proposition}
\theoremstyle{definition}
\newtheorem{defn}[thm]{Definition}
\theoremstyle{remark}
\newtheorem{rem}[thm]{Remark}
\numberwithin{equation}{section}

\newcommand{\R}{\mathbb R}

\newcommand{\Node}{\mathcal N}
\newcommand{\mcF}{\mathcal F}
\newcommand{\mcJ}{\mathcal J}
\newcommand{\mcV}{\mathcal V}
\newcommand{\mcT}{\mathcal T}
\newcommand{\mcP}{\mathcal P}

\newcommand{\bigO}{\mathcal{O}}
\newcommand{\E}{\mathbb{E}}
\newcommand{\Prob}{\mathbb{P}}
\newcommand{\hatE}{\hat{\mathbb{E}}}
\newcommand{\Zed}{\mathbb{Z}}

\newcommand{\C}{\mathbb C}

\newcommand{\ie}{\textit{i.e.\ }}
\newcommand{\eg}{\textit{e.g.\ }}

\begin{document}

\title{A Control-variable Regression Monte Carlo Technique for Short-term Electricity Generation Planning}%
\author{Magnus Perninge\footnote{M.\ Perninge is with the Department of Automatic Control, Lund University, Lund,
Sweden. e-mail: magnus.perninge@control.lth.se.} \footnote{M. Perninge is a member of the LCCC Linnaeus Center and the eLLIIT Excellence Center at Lund University. This work was supported by the Swedish Research Council through the grant NT-14, 2014-3774}}%
\maketitle
\begin{abstract}
In the day-to-day operation of a power system, the system operator repeatedly solves short-term generation planning problems. When formulating these problems the operators have to weigh the risk of costly failures against increased production costs. The resulting problems are often high-dimensional and various approximations have been suggested in the literature.

In this article we formulate the short-term planning problem as an optimal switching problem with delayed reaction. Furthermore, we proposed a control variable technique that can be used in Monte Carlo regression to obtain a computationally efficient numerical algorithm.
\end{abstract}


\section{Introduction}
Present operating conditions for power markets, with less expensive renewable generation facilities located far from the major loads, create a situation where market forces often drive the operation towards heavily loaded power grids. During operation of a power system, the system operator acts to maintain secure deliverance of electricity by re-dispatching the generation to avoid overloading the network. The optimization problem that the system operator is faced with is thus that of balancing the need for secure operation against economic efficiency.

The objective of this article is to formulate the problem of specification of optimal generation re-dispatch strategies, referred to as the short-term generation planning problem, as an optimal switching problem~\cite{Dixit}, where the system operator switches between different production set-points in a group of controllable power plants. It is, thus, a generalization of the work in~\cite{ObactMMOR}.

The resulting optimal switching problem is characterized by a high-dimensional driving noise process and execution delays. As computational time is of the essence in power system operation planning we also propose an efficient numerical algorithm based on applying a control variable in Monte Carlo regression~\cite{Longstaff}.

The use of control variables in Monte Carlo regression was suggested (but not implemented) in~\cite{Bouchard2012} as a means of improving the efficiency of American options pricing in high dimensions. In the American options case there is an obvious candidate for a control variable in the European option price, but for general optimal switching problems there is no natural counterpart.

To facilitate the algorithm we develop upper and lower bounds of the value function through lower dimensional approximations of the initial problem. The lower bound is then used as control variable when solving the initial problem.


The proposed method defines a novel approach to power system operation. Previously, applications of optimal switching to problems related to electricity markets and operation has mainly been aimed at the optimal operation of a single power plant to maximize the income of the producer~\cite{CarmLud,BollanMSwitch1,ElAsri,RAid}.

The remainder of the article is structured as follows. In the next section the system operators planning problem (Problem 1) is posed. In Section~\ref{sec:DynAlg} we show that a dynamic programming principle holds for the value function of Problem 1. Then, in Section~\ref{Sec:bnds}, the lower and the upper bounds are presented. Section~\ref{Sec:Num} is devoted to numerical solution. Here, we explain how a mix of value function approximation and policy approximation can be applied to increase computational efficiency for optimal switching problems with delays. In Section~\ref{Sec:NumEX} the proposed algorithm is tested on two numerical examples.

\section{Problem setup\label{Sec:Pset}}
In this section we formulate a risk adequate version of the short-term electricity generation planning problem. The power system model that we adopt is a static AC-network model that is common in optimal power flow formulations~\cite{LowOPFsurv1}.

\subsection{Network model}
Assume an AC-power network consisting of a set of nodes $\Node=\{0,1,\ldots,n_N\}$ connected by a set $E_0\subset \Node\times \Node$ of edges. In each node, $j\in\Node$, we define a complex voltage $V_j\in \C$ such that $|V_j|$ is the voltage magnitude at the node and $\angle V$ is the phase angle relative to the reference node $0\in \Node$.

The system nodes $\Node$ contain a set of nodes $\Node_G=\{1,\ldots,n_G\}$ where controllable generators are connected and a set of nodes, $\Node_L$, where loads are connected. For the purpose of power flow analysis, each node $j\in \Node$ is characterized by two complex variables $s_j$ and $V_j$, or, alternatively, by four real variables $p_j$, $q_j$, $|V_j|$ and $\angle V_j$, where $s_j=p_j+iq_j\in \C$ is the (complex) power injected into the system at Node $j$.

For a generator node $j\in \Node_G$ the injected active power $p_{j}$ is controlled by the producer and automated systems in the generator tries to keep the nodal voltage magnitude at a reference value $V_j^{ref}$. Nodes where $p_{j}$ and $|V_j|$ are known are referred to as PV-nodes. The production of reactive power is closely connected to the output voltage and the reactive power that a generator can produce is limited. Hence, the generator may not manage to keep the voltage magnitude at the reference value under stressed conditions. In this case the node becomes what is known as a PQ-node with $p_{j}$ still controlled by the producer, but with $q_j=\bar q_j$ and $|V_j|<V_j^{ref}$.

All other nodes are also PQ-nodes since, in steady state, the injected active and reactive power in these nodes are either set by the consumer demand (load nodes) or equal zero (empty nodes).

In the reference node, $j=0$, the mismatch of power injected into the system is balanced. Although this is a generator node the active power production in this node is not actively controlled. This node is referred to as the slack node and the convention of setting $\angle V_0=0$ is generally applied.\\

The Branch injection model (BIM) is a result of applying Kirchoff's current law to each node of the system. For each edge $(j,k)\in E_0$ let $y_{jk}\in \C$ be the admittance of the corresponding power line. The power injected into the network at node $j$ is then given by
\begin{equation}\label{PFekv}
s_j=\sum\limits_{k,\,(j,k)\in E_0} y_{jk}^H V_j (V_j^H-V_k^H),\: \forall j\in \Node,
\end{equation}
where the superscript $H$ denotes the Hermitian (complex conjugate transpose). The system of equations \eqref{PFekv} is known as the power flow equations for the BIM. Together with the constraints
\begin{subequations}\label{PVQekv}
\begin{align}
|V_j|&\leq V_j^{ref},\quad \forall j\in \Node_G,\\
q_j&\leq \bar q_j,\quad \forall j\in \Node_G,\\
(|V_j| - V_j^{ref})(q_j - \bar q_j)&=0,\quad \forall j\in \Node_G,
\end{align}
\end{subequations}
the power flow equations give a criterion for feasibility of operating points.
\begin{rem}
Additional operation constraints such as line-thermal limits and limits on the voltage magnitudes in load-nodes are often added to the power flow equations.
\end{rem}

To simplify notation we will use the common assumption that $\varphi_j=\angle s_j$ is constant in all load nodes. We thus end up with a vector of system parameters $(p_G,p_D)\in \R^{n_G}_+\times \R^{n_L}_+$, to which the nodal-voltages will adapt.

\begin{rem}
We take the demand, $p_D$, as a parameter rather than the load, $p_L$, since the demand can be assumed to take any value in $\R_+^{n_L}$ while the actual nodal-loading is restricted by power flow feasibility.
\end{rem}

\subsection{Feasible sets and contingencies}
By the physical laws of electricity transfer, the set of parameters for which the extended set of power flow equations \eqref{PFekv}-\eqref{PVQekv} have a solution is limited~\cite{VCutsem}. We define
\begin{equation*}
G_0:=\{(p_G,p_D)\in \R^{n_P}_+:\exists V\in \C^{n+1}\text{ that solves \eqref{PFekv} and \eqref{PVQekv}}\},
\end{equation*}
so that $G_0$ is the set of all parameter vectors for which there is a solution to the power flow equations.

When an event occurs that leads to failure of one or more of the power lines, the network topology will change. A typical such event is a short-circuit after which a power line is disconnected. The failure of a power line is referred to as a contingency. Assume that a set of $n_c$ contingencies have been deemed plausible, after which the remaining power lines make up the edges $E_i\subset E_0$ for $i=1,\ldots,n_c$ in the network graph.

To each post-contingency network graph corresponds a new feasible set and we define
\begin{equation*}
G_i:=\{(p_G,p_D)\in \R^{n_P}:\exists V\in \C^{n+1}\text{ that solves \eqref{PFekv} (with $E_0$ replaced by $E_i$) and \eqref{PVQekv}}\}.
\end{equation*}
We denote by $\partial G_i$ the boundary of the set $G_i$. Let $D_G\subset \R^{n_G}_+$ be a convex and compact subset of the production levels that $p_G$ can take and define, for each $p_G\in D_G$, the set $G_i(p_G):=\{p_D\in\R^{n_L}_+:\:(p_G,p_D)\in G_i\}$.

\begin{ass}\label{AssG}
We make the following assumptions/modifications of the closed sets $G_i\subset \R^{n_G}\times \R^{n_L}_+$, for $i=0,1,\ldots,n_c$,
\begin{enumerate}[A)]
  \item $G_i$ is convex and compact.
  \item For all $p_G\in D_G$ and $p_D\in \R^{n_L}_+$ there exists $z\geq 0$ such that $(p_G,p_D-z)\in G_i$.
  \item For all $p_G\in D_G$ we have $\left(\{p_G\}\times G_i(p_G)\right)\cap \emph{int}\,G_i\neq\emptyset$.
  \item We extend $G_i$ using the load over-satisfaction concept~\cite{Baldick}: if $(p_G,p_D)\in G_i$, then $(p_G,p_D')\in G_i$ for all $p_D'\in \R^{n_L}_+$ with $p_D'\leq p_D$.\label{assDOMIN}
\end{enumerate}
\end{ass}

\begin{rem}
Assumption~\ref{AssG}.A was long believe to hold for all systems until a counterexample was presented in~\cite{Hiskens}. Since then a number of counterexamples have been presented, and finding conditions under which Assumption~\ref{AssG}.A (or slightly weaker versions) holds is today an active field of research (see \eg~\cite{ZhangTse,LowOPFsurv1,LowOPFsurv2}).

Assumptions~\ref{AssG}.B is natural as the operation should never be dispatched in a way that no load disruption can save the system. Assumptions~\ref{AssG}.C is purely technical and should hold in all realistic cases.

The rationality behind Assumption~\ref{AssG}.\ref{assDOMIN} is that, in case of an emergency, it should be a possible to increase the load (without a substantial cost) in order to save the system.
\end{rem}

\subsection{Voltage stability}
The power consumed at a load node, $j\in\Node_L$, is given by $s_j=-y_{L,j}^H|V_j|^2$, where $y_{L,j}\in \C$ is the load admittance. To account for voltage drops, mechanisms embedded either in the connection to lower level subsystems or in the loads themselves alter the load admittance, to maintain balance between the consumed power and the demanded power\footnote{An example would be the charger of a laptop that can handle the different voltage levels in distribution grids of different regions by adapting the internal impedance.}. For pedagogical purpose we can model this behavior by the simplified characteristic
\begin{equation*}
\frac {d}{dt} \Re[y_{L,j}]=\frac{1}{K_{\text{LR},j}}(p_{Dj}-p_{Lj}),
\end{equation*}
where $p_{Lj}=-\Re[s_j]$ and $K_{\rm{LR}}$ a positive constant. If, at any time, the active power demand is above what the system can supply we must have $p_{Dj}>p_{Lj}$. This gives rise to an increasing load conductance $\Re[y_{L,j}]$. By the nature of power flow, increasing a sufficiently high load conductance leads to an even smaller load~\cite{VCutsem}. This makes power flow feasibility critical in the sense that unless actions are taken a voltage collapse will occur in the manner of an accelerating real part of $y_{L,j}$ with the implication that $|V_j|\to 0$.
If, on the other hand, the vector of active power productions and complex power demands is inside the feasible domain, a voltage stable operating point exists~\cite{VCutsem}.

\subsection{Load disruption}
In light of the previous discussion we realize that operating the system in a way that leads to infeasible demand may have catastrophical consequences. However, there are countermeasures that can be taken if such an event should occur. The most effective countermeasure is to disconnect parts of the consumers to salvage the remainder of the system. As modern society is very dependent on having a continuous supply of electricity, there is a large societal cost associated to disconnection of consumers.

We assume that if contingency $i$ occurs, load is immediately disconnected to return operation to the feasible set. We assume that the cost of disconnecting loads is linear in the amount of disconnected load, but may differ between different nodes. The immediate cost of disconnecting loads to satisfy feasibility constraints following contingency $i$ is, thus, the minimum to
\begin{subequations}
\label{LSprob}
\begin{align}
\min\limits_{z_i\in \R^{n_L}_+}  \quad\quad & c^\top z_i\\
\text{subj. to}\quad\:\: & (p_G,p_D-z_i) \in G_i´. \label{bivEENS}
\end{align}
\end{subequations}
Here, $c\in \R^{n_L}$ is the vector of nodal costs for energy not served. We thus have $c>0$.



Also in the non-contingent situation the operating point has to stay within the region of power flow feasibility.
\begin{defn}
Let $z_i^*:D_G\times \R^{n_L}_+\to \R^{n_L}_+$ be such that $z_i^*(p_G,p_D)$ is a solution to \eqref{LSprob}, for $i=0,\ldots,n_c$.
\end{defn}

\begin{prop}\label{prop:zLIP}
There are constants $C_{1}^{LS}>0$ and $C_{2}^{LS}>0$ such that, for $i=0,\ldots,n_c$, we have
\begin{equation}
|c^\top z^*_i(p_G,p_D)|\leq C_{1}^{LS} \|p_D\|,\quad \forall p_G\in D_G,\, p_D\in\R^{n_L}_+,
\end{equation}
and
\begin{align}\label{zprop2}
|c^\top z^*_i(p_G,p_D)-c^\top z^*_i(p_G',p_D')|\leq C_{2}^{LS} (\|p_G-p_G'\|+\|p_D-p_D'\|),
\end{align}
$\forall (p_G,p_D),\,(p_G',p_D')\in D_G\times \R^{n_L}_+$.
\end{prop}

\bigskip
\bigskip

\noindent\emph{Proof.} First we note that by Assumptions~\ref{AssG}.B a solution $z^*_i(p_G,p_D)$ always exists, but is not necessarily unique. Specifically, by Assumptions~\ref{AssG}.\ref{assDOMIN} we have that $(p_G,0)\in G_i$, for all $p_G\in D_G$. Hence,
\begin{align*}
|c^\top z^*_i(p_G,p_D)|\leq c^\top p_D\leq \|c\| \|p_D\|.
\end{align*}
For the second inequality we pick $p_G,\,p_G'\in D_G$ such that $c^\top z^*_i(p_G,p_D)>0$ and $c^\top z^*_i(p_G,p_D)\geq c^\top z^*_i(p_G',p_D')$ .

By Assumption~\ref{AssG}.A and the optimality of $z_i^*(p_G,p_D)$ there is a unit vector $\eta:=[\eta_G^\top\: \eta_L^\top]^\top \in S^{n_p-1}$ such that
\begin{equation}\label{GetaDEF}
\eta_G^\top(\hat p_G-p_G)+\eta^\top_L(\hat p_D-(p_D-z_i^*(p_G,p_D))) \leq 0,\quad\forall(\hat p_G,\hat p_D)\in G_i
\end{equation}
and $z^*(p_G,p_D)$ is a solution to the optimization problem
\begin{align*}
  \min_{z\in\R^{n_L}_+} \:\:&\quad \quad c^\top z, \\
  \text{subj. to}  &\quad \quad\eta^\top_L(z_i^*(p_G,p_D)-z)\leq 0.
\end{align*}
Furthermore, $\eta_L\geq 0$ and by Assumption~\ref{AssG}.C, $\eta_{L,j}>0$ for some $j\in\{1,\ldots,n_L\}$. Let $z'(p_G',p_D')$ solve
\begin{align*}
  \min_{z\in\R^{n_L}_+} \:\:&\quad \quad c^\top z, \\
  \text{subj. to}  &\quad \quad\eta^\top_L(z_i^*(p_G,p_D)-z)\leq \eta_G^\top(p_G-p_G')+\eta^\top_L(p_D-p_D').
\end{align*}
Then by \eqref{GetaDEF} we have $c^\top z'(p_G',p_D') \leq c^\top z_i^*(p_G',p_D')$ and, hence,
\begin{equation*}
|c^\top z^*_i(p_G,p_D)- c^\top z^*_i(p_G',p_D')|\leq |c^\top z^*_i(p_G,p_D)-c^\top z'(p_G',p_D')|.
\end{equation*}
It is now a simple exercise to show that
\begin{equation*}
|c^\top z^*_i(p_G,p_D)-c^\top z'(p_G',p_D')|\leq\left(\eta_G^\top(p_G-p_G')+\eta^\top_L(p_D-p_D')\right)\min_{j\in\{1,\ldots,n_G\}}\frac {c_j}{\eta_{L,j}},
\end{equation*}
which finishes the proof of the proposition.\qed\\

\subsection{Stochastic model}
We assume that the system operator observes the demand $p_D$ and models its stochastic behavior in the time period $[0,T]$ with the stochastic process $X_t$ on the filtered probability space $(\Omega,\mcF,\{\mcF_t\}_{0\leq t \leq T},\Prob)$. Here, $\mcF_t=\sigma(X_s:0\leq s\leq t)$ is the filtration generated by $X_t$. We assume that $X_t=m(t)+Z_t$, where $m:[0,T]\to\R^{n_L}_+$ is a smooth function (with $\|m(t)-m(s)\|\leq C^m|t-s|$ for some $C^m>0$) representing a forecast, and the forecast error $Z_t$ is the strong solution to a linear Stochastic Differential Equation (SDE) driven by an $n_L$-dimensional Brownian motion $W_t$, as follows:
\begin{align*}
dZ_t & = -\gamma Z_tdt + \sigma dW_t;\quad 0<t\leq T,\\
Z_0 &= x_0-m(0).
\end{align*}
Here, $\gamma \geq 0$ is a scalar and $\sigma\in\R^{n_L\times n_L}$ is a non-singular matrix.

\begin{rem}
With the above model the demand process will take negative values with a positive probability. As negative demand is not physically reasonable we take $X_t^+:=\max(X_t,0)$ to be the actual demand, where the maximum is taken componentwise, but we still assume that the operator observes $X_t$.
\end{rem}

\subsection{Control model}
We will assume that the operator will control production by switching between a number of pre-specified production set-points. This gives a generalization of the regulating market setup that is common in deregulated power markets~\cite{ObactMMOR}.

\subsubsection{Operator controlled production}
We assume that, for each node $k \in \Node_G$, the system operator can control $p_{G,k}$ by ordering the owner of the corresponding power plant to switch the production set-point. The set-points are restricted to a number of pre-specified production levels\footnote{If the operator cannot control the production in a subset of the nodes then we simply let $\Theta_k:=\{\theta_{k,1}\}$ for these nodes.}  $\Theta_k:=\{\theta_{k,1},\ldots, \theta_{k,n_k^\theta}\}$, with $\theta_{k,1}<\cdots<\theta_{k,n_k^\theta}$. At all times $t\in [0,T]$ the production $p_G(t)$ will thus take values in $D_G:=[\theta_{1,1},\theta_{1,n_1^\theta}]\times\cdots\times [\theta_{n_G,1},\theta_{n_G,n_{n_G}^\theta}]$.

\begin{defn}
We define the set of production set-points $\Theta:=\Theta_1 \times \cdots \times \Theta_{n_G}$.
\end{defn}

If a switch from set-point $l\in\Theta_k$ to set-point $m\in\Theta_k$ at node $k \in \Node_G$ is ordered\footnote{We will constantly abuse notation and write $l$ for $\theta_{k,l}$, when the context is obvious.}, the production will not immediately jump to the new production level. Rather, the production will follow a monotone continuous function $\theta^T_{k,(l,m)}:[0,\delta_{k,(l,m)}]\to [\theta_{k,l},\theta_{k,m}]$, with $\theta^T_{k,(l,m)}(0)=\theta_{k,l}$ and $\theta^T_{k,(l,m)}(\delta_{k,(l,m)})=\theta_{k,m}$ at $\delta_{k,(l,m)}>0$ time units after the switching order. The rate of production change is bounded by physical constraints that render limits called maximal ramp rates of the power plants.


The restriction to a finite set of production set-points is not unrealistic as power plants with multiple turbines usually have multiple production levels of locally optimal efficiency, where producers prefer to operate their units. 

There may also be restrictions on the order in which producers are willing (and able) to switch between different levels of production. To model this we assume that, for each $k\in \Node_G$, there is a directed graph $\Lambda_k \subset \Theta_k\times \Theta_k$, such that $(l,m)\in \Lambda_k$ iff the owner of the power plant at node $k$ allows a switch from level $\theta_{k,l}$ to level $\theta_{k,m}$.

The physically constrained ramp-rates now implies that there exists a $C^{\textit{RR}}>0$ such that, for each $k\in \Node_G$ and each $(l,m)\in \Gamma_k$, we have
\begin{equation}\label{ekv:maxRR}
|\theta_{k,(l,m)}(s)-\theta_{k,(l,m)}(s')|\leq C^{\textit{RR}}|s-s'|,\quad \forall s,\,s'\in [0,\delta_{k,(l,m)}].
\end{equation}

As altering the production level, possibly by starting up a new turbine, generally leads to additional fuel costs and/or wear and tear on equipment there is a fixed switching cost $K_{k,(l,m)}\geq 0$ associated to the switch $(l,m)\in \Lambda_k$, for $k=1,\ldots,n_G$.

The controllable generator at node $k\in \Node_G$ is either running in a steady mode, producing a constant active power output $\theta_{k,l}\in\Theta_k$, or is in the midst of a transition from production level $\theta_{k,l}\in\Theta_k$ to level $\theta_{k,m}\in\Theta_k$, with $(l,m)\in \Lambda_k$. To describe the operation mode of the system we use the state of each controllable generator.

\begin{defn}
We define the set of operation modes to be $\Gamma:=(\Theta_1\cup\Lambda_1)\times\cdots\times (\Theta_{n_G}\cup\Lambda_{n_G})$.
\end{defn}

The operator will thus control the production by switching between different production set-points in the set $\Theta$. Due to delays, mainly caused by maximal ramp-rates of production facilities, the operation mode will take values in the set $\Gamma$. For each operation mode $a\in\Gamma$ the set of switches available to the operator is constrained by the allowed switches as defined in the sets $\Lambda_k$, but also by the fact that a generator that is in transition between two production levels cannot be switched to a new level. We thus define the set of admissible switches.

\begin{defn}\label{def:Lambda}
We define the set of all possible switches as $\Lambda=\cup_{k\in \Node_G} (\{k\} \times \Lambda_k)$. For each operation mode $a\in\Gamma$, the set of admissible switches is
\begin{equation}
\Lambda(a):=\{(k,(l,m))\in \Lambda|\:a_k=l,\,(l,m)\in\Lambda_k\}.
\end{equation}
\end{defn}

In addition to the operator controlled switches between modes, the completion of transitions will render alterations in the operation mode. To keep track of the progression of transitions we introduce the \emph{transition progression vector}.
\begin{defn}
For each operation mode $a\in\Gamma$, the transition progression vector $r\in H(a)\subset \R^{n_G}_+$ specifies the time spent in transition for each of the transiting generators. With the convention that the transition progression is always zero in stationary production we thus get
\begin{equation}
H(a):=\{(r_1,\ldots,r_{n_G})\in \R^{n_G}_+|\: r_k=0,\text{ if } a_k\in \Theta_k,\:r_k\in[0,\delta_{k,(l,m)})\text{ if }a_k=(l,m)\}.
\end{equation}
\end{defn}

\subsubsection{Operation strategies}
To control the power flows through the network, the system operator re-dispatches generation by requesting that the producers switch between their announced production levels. The control for our problem will be a sequence of \emph{switching requests}\footnote{The term request is used here to emphasize that productional changes are not instantaneous.} $\{\beta_1, \beta_2,\ldots\} $ and the corresponding \emph{intervention times} $\{\tau_1, \tau_2,\ldots\}$.

Here, $\beta_j\in\Lambda$ will represent the $j^{th}$ request that the operator makes and $\tau_j\in[0,T]$ the time that the request is made.

\begin{defn}
For each $\beta=(k,(l,m))\in\Lambda$ we let $\beta^{\Node}=k$, $\beta^{f}=l$ and $\beta^t=m$, so that $\beta^\Node$ represents the generator the request concerns, $\beta^f$ the level that the transition is from and $\beta^t$ the level that the transition is to.
\end{defn}

Assume that at time $t=0$ production starts in operation mode $\alpha_0$, with transition progression $\xi_0\in H(\alpha_0)$.

Applying the (stochastic) control $v=(\tau_1,\ldots,\tau_N;\beta_1,\ldots,\beta_N)$ will lead to a process of operation modes $\alpha(\cdot,v):[0,T]\times \Omega \to\Gamma$. For each $k\in \Node_G$, let $v^k=(\tau^k_1,\ldots,\tau^k_{N_k};\beta_{k,1},\ldots,\beta_{k,N_k})$ be the subsequence of the control $v$ containing all $\beta_j$ for which $\beta_j^{\Node}=k$ and the corresponding intervention times. Then
\begin{equation*}
\alpha_k(t,v)=\sum\limits_{j=0}^{N_k}\left((\beta^f_{k,j},\beta^t_{k,j})1_{[\tau^k_j,\tau^k_j+\delta_{\beta_{k,j}})}(t) +\beta_{k,j}^t 1_{[\tau^k_j+\delta_{\beta_{k,j}},\tau^k_{j+1})}(t)\right),
\end{equation*}
with $\tau^k_{0}:=-\xi_{0,k}$, $\beta_{k,0}:=(k,\alpha_{0,k})$ and $\tau^k_{N_k+1}:=\infty$.

The control $v$ is an \emph{admissible control} if $0\leq\tau_1\leq\cdots\leq\tau_N$ is a sequence of $\mcF_t$--stopping times and the sequence $(\beta_1,\ldots,\beta_N) $ is such that $\beta_j$ is $\mcF_{\tau_j}$--measurable and $\beta_j\in \Lambda(\alpha(\tau_j^-))$ for $j=1,\ldots,N$, and $\tau^k_j<\tau^k_{j+1}$, for $j=1,\ldots,N_k$, $\forall k\in \Node_G$.

\begin{defn}
The set of all admissible controls is denoted $\mcV$.
\end{defn}

As the length of the operating period is finite, $T<\infty$, and the minimal delay $\min\limits_{\beta\in \Lambda} \delta_{\beta}=:\delta_{\min}>0$, the number of switches is bounded,
\begin{equation}\label{NswBOUND}
N\leq \frac{T}{\delta_{\min}}n_G.
\end{equation}
This means that the set of admissible controls only contains controls that make a finite number of switching requests.

\subsubsection{Production process}
Using the definition of an admissible control we can express the production in the controllable generation nodes as a $\mcF_{t}$--measurable process $\zeta:[0,T]\times \Omega \to \R^{n_G}$. First, we express the transition progression vector, $\xi(t,v)\in H(\alpha(t,v))$, as
\begin{equation*}
\xi_k(t,v)=\sum\limits_{j=0}^{N_k}(t-\tau^k_j)1_{[\tau^k_j,\tau^k_j+\delta_{\beta_{k,j}})}(t),
\end{equation*}
$\forall k\in \Node_G$. Now, we define the function $\zeta:\cup_{a\in\Gamma}\{a\}\times H(a)\to U_G$ as
\begin{equation*}
\zeta_k(a,r):=\left\{\begin{array}{l l} \theta_{k,a_k} \quad & \text{if }a_k\in \Theta_k \\
\theta^T_{k,a_k}(r_k)\quad & \text{if }a_k\notin \Theta_k\end{array}\right.
\end{equation*}
for each $k\in \Node_G$. This allows us to define the \emph{production process} corresponding to the control $v\in\mcV$ as $(\zeta(\alpha(t,v),\xi(t,v));\: 0\leq t\leq T)$.

\subsection{Operational cost}
There are two different types of operational costs that have to be considered when re-dispatching the generation. On one side we have the cost of production that may differ significantly between different production facilities. On the other side we have the expected cost from having to disconnect parts of the consumers if the system becomes too stressed, generally as the result of a contingency.

\subsubsection{Production cost}
The expected production cost when applying the control $v\in\mcV$ is
\begin{equation}
J_{\text{PC}}(v)=\E\left[\int_0^T C_G(\zeta(\alpha(t,v),\xi(t,v)))dt + \sum\limits_{\tau_i\leq T}K_{\beta_i}\right],
\end{equation}
where $\E$ is expectation with respect to the probability measure $\Prob$ and the costs are the continuous function $C_G:D_G\to \R$ representing the additional running cost induced by utilizing the controllable power plants and the second part is the sum of all the switching costs.

%
%

\begin{rem}
Note that $C_G$ is not the running cost of the power plants, but rather the result of a difference in running costs between the controllable production facilities and the slack node, where the production is automatically dispatched to maintain balance between production and demand. This difference may account for changes in losses as the production is re-dispatched, as long as the losses do not depend on the demand.
\end{rem}

\subsubsection{Cost of unserved demand}
The other part of the expected operational cost comes from the societal cost of unserved demand. To make a realistic model of the actual load disruption process we would have to augment the filtration to include a stochastic process that describes the state of the grid as it jumps between different contingency settings. Furthermore, a number of emergency actions may become available to the operator when a contingency occurs and the system enters a state termed \emph{emergency state}. As the number of contingencies is often very large, this approach is far from being computationally tractable. Instead we model the cost of unserved demand as a weighted sum of the minimal economic value of load shedding necessary to return the parameter vector to the feasible region. The weights, $q_i\geq 0$, corresponding to each of the different contingencies, should be chosen to appropriately represent the failure rates of components and the expected length of an outage.

The expected cost of energy not served during the time period $[0,T]$, when applying the control $v\in\mcV$, is thus modeled as
\begin{equation*}
J_{\text{EENS}}(v)=\E\left[\int_0^T  \sum_{i=0}^{n_c}q_i c^\top z_i^*(\zeta(\alpha(t,v),\xi(t,v)),X_t^+)dt\right].
\end{equation*}
Here, as detailed above,
\begin{itemize}
  \item $n_c$ is the number of contingencies that we consider,
  \item $q_i$ is the weight of system state $i$, for $i=0,\ldots,n_c$,
  \item $c\in \R^{n_L}_+$ is the linear cost of load disruption in the different nodes, $c>0$, and
  \item $z_i^*:\R^{n_G}\times\R_+^{n_L}\to \R_+^{n_L}$ is the optimal load disruption when the system is in configuration $i$, for $i=0,\ldots,n_c$.
\end{itemize}

\begin{rem}
For generation expansion planning, adequacy measures such as loss of load probability (LOLP) and expected energy not served (EENS), have long been used as constraints, but surprisingly it is not, with a few notable exceptions, common practice in power system operation planning.
\end{rem}

\subsection{Control problem}
The control problem that the operator is faced with is that of weighing productional costs against the expected cost of unserved demand.\\

\textbf{Problem 1.} Find $v^*$ and $J(v^*)$ that solves
\begin{equation}
\inf\limits_{v\in\mcV}J(v),
\end{equation}
where $J(v):=J_{\text{PC}}(v)+J_{\text{EENS}}(v)$.\qed\\

\begin{rem}
Note that Problem 1 is not a Markovian problem in $(X_t,\alpha(t,v))$ as it depends on past control actions. However, a simple augmentation of the state space to $(X_t,\alpha(t,v),\xi(t,v))$ renders the Markovian property which is essential for solution by dynamic programming.
\end{rem}
To make notation shorter we introduce
\begin{equation*}
f(x,a,r):=C_G(\zeta(a,r)) + \sum_{i=0}^{n_c}q_i c^\top z_i^*(\zeta(a,r),x^+),
\end{equation*}
which allows us to write
\begin{equation}
J(v)=\E\left[\int_0^T  f(X_t,\alpha(t,v),\xi(t,v))dt+\sum\limits_{\tau_i\leq T}K_{\beta_i}\right].
\end{equation}


\begin{prop}\label{fprop}
There exist constants $C_1^f>0$ and $C_2^f>0$ such that for all $a,a'\in\Gamma$, we have
\begin{align}\label{fprop1}
|f(x,a,r)|&\leq C_1^f(1+\|x\|),\quad \forall (x,r)\in \R^{n_L}\times  H(a)\\
\label{fprop2}
|f(x,a,r)-f(x',a',r')| &\leq C_2^{f}(\|x-x'\|+\|\zeta(a,r)-\zeta(a',r')\|),
\end{align}
for all $x,x'\in \R^{n_L}$ and all $r\in H(a)$ and $r'\in H(a')$.
\end{prop}

\noindent\emph{Proof.} We have that
\begin{equation*}
|f(x,a,r)|\leq c^\top (x\vee 0)\sum_{i=0}^{n_c}q_i+ \sup_{p_G\in D_G} |C_G(p_G)|.
\end{equation*}
Since $C_G$ is continuous and $D_G$ is compact it follows that the supremum is attained and bounded, $|f(x,a,r)|\leq\max\{\max_{p_G\in D_G} |C_G(p_G)|,\sum_{i=0}^{n_c}q_i\|c\|\}(1+\|x\|)=:C^f_1(1+\|x\|)$.

The second property follows from \eqref{zprop2} and the fact that $C_G$ is Lipschitz.\qed\\


\section{Dynamic programming algorithm\label{sec:DynAlg}}

\begin{defn}
For $t\in [0,T]$ and $x\in \R^{n_L}$ we let, for all $s\in [t,T]$, $X_s^{t,x}:=X_s|X_t=x$. Furthermore, we let $\mcV_t$ be the set of all admissible controls $v\in\mcV$ with $\tau_1\geq t$ and define $\alpha^{t,a,r}(s,v)$ and $\xi^{t,a,r}(s,v)$, for all $v\in \mcV_t$ and $s\in[t,T]$, as
\begin{equation}\label{alphaDEF2}
\alpha^{t,a,r}_k(s,v)=\sum\limits_{j=0}^{N_k}\left((\beta^f_{k,j},\beta^t_{k,j})1_{[\tau^k_j,\tau^k_j+\delta_{\beta^k_j})}(s) +\beta_{k,j}^t 1_{[\tau^k_j+\delta_{\beta^k_j},\tau^k_{j+1})}(s)\right),
\end{equation}
and
\begin{equation}\label{xiDEF2}
\xi^{t,a,r}_k(s,v)=\sum\limits_{j=0}^{N_k}(s-\tau^k_j)1_{[\tau^k_j,\tau^k_j+\delta_{\beta^k_j})}(s),
\end{equation}
where $\tau^k_{0}:=t-r_k$, $\beta_{k,0}:=(k,a_k)$ and $\tau^k_{N_k+1}:=\infty$. Finally, we let $\mcV_{t,a,r}$ be the set of controls $v\in\mcV_t$ that are admissible for $\alpha(t^-,v)=a$ and $\xi(t^-,v)=r$.
\end{defn}

\begin{defn}
For all $(t,x,a,r)\in [0,T]\times \R^{n_L}\times \mathop{\cup}\limits_{a\in \Gamma}\{a\times H(a)\}$, we define the cost-to-go starting the problem afresh at time $t\in [0,T]$, with $X_t=x$, $\alpha(t^-)=a$, $\xi(t^-)=r$ and applying the control $v\in \mcV_t$ as
\begin{equation}
J(t,x,a,r;v):=\E\left[\int_t^T  f(X_s^{t,x},\alpha^{t,a,r}(s,v),\xi^{t,a,r}(s,v))ds + \sum\limits_{\tau_i\leq T}K_{\beta_i}\right].
\end{equation}
Furthermore, we define the value function $V:[0,T]\times\R^{n_L}\times \cup_{a\in\Gamma}(\{a\}\times H(a))\to\R$ as
\begin{equation}\label{VFdef}
V(t,x,a,r):=\inf\limits_{v\in \mcV_{t,a,r}}J(t,x,a,r;v).
\end{equation}
\end{defn}

For each admissible control $v\in\mcV$ and each $\nu\in \mcT_t$ with $\nu\leq T$, where $\mcT_t$ is the set of $\mcF_t$--stopping times $\tau$ with $\tau\geq t$, we have
\begin{align*}
J(t,x,a,r;v)=\E\left[\int_t^\nu  f(X_s^{t,x},\alpha^{t,a,r}(s,v),\xi^{t,a,r}(s,v))ds + \sum\limits_{t\leq \tau_i \leq \nu}K_{\beta_i}+J(\nu,X_\nu^{t,x},\alpha^{t,a,r}(\nu,v),\xi^{t,a,r}(\nu,v);v)\right],
\end{align*}
by the strong Markov property (see \eg \cite{XYZbook}).

The aim of this section is to prove that the principle of optimality holds for the value function of Problem 1:
\begin{thm} \label{Popt1}(Principle of optimality)
For each $t\in[0,T]$, $x\in\R^{n_L}$, $a\in\Gamma$ and $r\in H(a)$, we have
\begin{align*}
V(t,x,a,r)=\inf\limits_{v\in \mcV_t}\E\left[\int_t^\nu  f(X^{t,x}_s,\alpha^{t,a,r}(s,v),\xi^{t,a,r}(s,v))ds + \sum\limits_{\tau_i \leq \nu}K_{\beta_i}+ V(\nu,X^{t,x}_\nu,\alpha^{t,a,r}(\nu,v),\xi^{t,a,r}(\nu,v)) \right],
\end{align*}
for each $\nu\in \mcT_t$, where $\mcT_t$ is the set of $\mcF_t$--stopping times $\tau$ with $\tau\geq t$.
\end{thm}

\bigskip

The principle of optimality is essential in dynamic programming and is the foundation on which most numerical solution algorithms are based. The proof will be divided into a number of lemmas and propositions.


\begin{lem}\label{zetaLEM}
There exists a constant  $C^\zeta>0$ such that for every $t^A,t^B\in [0,T]$, $a\in \Gamma$ and $r^A,r^B\in H(a)$ and every $v:=(\tau_1,\ldots,\tau_{N};\beta_1,\ldots,\beta_{N})\in \mcV_{t^A,a,r^A}$, there is a $\tilde v:=(\tilde\tau_1,\ldots,\tilde\tau_{\tilde N};\tilde\beta_1,\ldots,\tilde\beta_{\tilde N})\in \mcV_{t^B,a,r^B}$ such that
\begin{equation}\label{ekv:zetaBND}
\max_{s\in [t^A\vee t^B,T]}\|\zeta(\alpha^{t^A,a,r^A}(s,v),\xi^{t^A,a,r^A}(s,v))-\zeta(\alpha^{t^B,a,r^B}(s,\tilde v),\xi^{t^B,a,r^B}(s,\tilde v))\|\leq C^\zeta(\|r^A-r^B\|+|t^A-t^B|),
\end{equation}
and
\begin{equation}\label{ekv:lowerSWcost}
\sum_{\tilde\tau_j \leq T}K_{\tilde\beta_j}\leq \sum_{\tau_j \leq T}K_{\beta_j},
\end{equation}
for all $\omega\in\Omega$.
\end{lem}

\noindent\emph{Proof.}
Let $d:=\inf\{s\geq 0: r^B-r^A+se_a\geq 0\}+(t^B-t^A)^+\leq \|r^A-r^B\|+|t^A-t^B|$, where $e_a\in\{0,1\}^{n_G}$, with $(e_a)_k:=1_{\Lambda_k}(a_k)$ for $k=1,\ldots,n_G$. We define the control $\tilde v$ as $\tilde \tau_j:=\tau_j+d$, and $\tilde \beta_j:=\beta_j$, for $j=1,\ldots,\tilde N$, where $\tilde N:= \max \{j\in\{0,\ldots,N\}:\: \tau_j+d<T\}$.

Since $v$ is an admissible control and $d\geq 0$ we have that $\tilde\tau_j$ is a $\mcF_t$--stopping time and $\tilde \beta_j$ is $\mcF_{\tilde\tau_j}$--measurable for $j=1,\ldots,\tilde N$. Furthermore, $\tilde \tau_j^k<\tilde \tau_{j+1}^k$, for $j=1,\ldots,\tilde N^k-1$, for all $k\in\Node_G$.

What remains is to show that $\alpha^{t^B,a,r^B}(s,\tilde v)$ is well defined for $s\in [t^B,T]$, and that $\tilde \beta_j\in\Lambda(\alpha^{t^B,a,r^B}(\tilde \tau_j^-,\tilde v))$ for $j=1,\ldots,\tilde N$.

Clearly, $\alpha^{t^B,a,r^B}(\tilde \tau_1^-\wedge T,\tilde v)$ exists. For all $a'\in \Gamma$ we let $\mathcal{K}^T(a'):=\{k\in \Node_G:\: a'_k \in \Lambda_k\}$. We have that $(\alpha^{t^A,a,r^A}(\tau^-_1\wedge T,v))_k=(\alpha^{t^B,a,r^B}(\tilde \tau^-_1\wedge T,\tilde v))_k$, for all $k\in \mathcal{K}^T(\alpha^{t^A,a,r^A}(\tau_1^-\wedge T,v))$, $\forall \omega\in\Omega$. Now, by Definition~\ref{def:Lambda} we find that $\Lambda(\alpha^{t^A,a,r^A}(\tau^-_1\wedge T,v))\subset\Lambda(\alpha^{t^B,a,r^B}(\tau^-_1\wedge T,\tilde v))$. Hence, $\tilde\beta_1\in \Lambda(\alpha^{t^B,a,r^B}(\tilde \tau^-_1,\tilde v))$, whenever $\tilde N\geq 1$.

By a similar argument $\alpha^{t^B,a,r^B}(\tilde \tau^-_2\wedge T,\tilde v)$ is well defined and $\tilde\beta_2\in \Lambda(\alpha^{t^B,a,r^B}(\tilde \tau^-_2,\tilde v))$, whenever $\tilde N\geq 2$. By induction it follows that $\tilde v \in \mcV_{t^B,a,r^B}$.\\

To prove the inequality we note that for all $s\in [t^A+d,T]$ we have
\begin{align*}
&\|\zeta(\alpha^{t^A,a,r^A}(s,v),\xi^{t^A,a,r^A}(s,v))-\zeta(\alpha^{t^B,a,r^B}(s,\tilde v),\xi^{t^B,a,r^B}(s,\tilde v))\|\leq
\\
&\qquad\quad\|\zeta(\alpha^{t^A,a,r^A}(s-d,v),\xi^{t^A,a,r^A}(s-d,v))-\zeta(\alpha^{t^B,a,r^B}(s,\tilde v),\xi^{t^B,a,r^B}(s,\tilde v))\|
\\
&\qquad+\|\zeta(\alpha^{t^A,a,r^A}(s-d,v),\xi^{t^A,a,r^A}(s-d,v))-\zeta(\alpha^{t^A,a,r^A}(s,v),\xi^{t^A,a,r^a}(s,v))\|.
\end{align*}
For $k=1,\ldots,n_G$ we have
\begin{align*}
\|\zeta_k(\alpha^{t^A,a,r^A}(s-d,v),\xi^{t^A,a,r^A}(s-d,v))-\zeta_k(\alpha^{t^B,a,r^B}(s,\tilde v),\xi^{t^B,a,r^B}(s,\tilde v))\|\leq C^{RR}|r_k^B -r_k^A +t^A+d-t^B|
\end{align*}
and
\begin{align*}
\|\zeta_k(\alpha^{t^A,a,r^A}(s-d,v),\xi^{t^A,a,r^A}(s-d,v))-\zeta_k(\alpha^{t^A,a,r^A}(s,v),\xi^{t^A,a,r^A}(s,v))\|\leq C^{RR} d.
\end{align*}
Hence,
\begin{align*}
\|\zeta(\alpha^{t^A,a,r^A}(s,v),\xi^{t^A,a,r^A}(s,v))-\zeta(\alpha^{t^B,a,r^B}(s,\tilde v),\xi^{t^B,a,r^B}(s,\tilde v))\|&\leq n_G C^{RR}(\|r^B -r^A \| +|t^A-t^B|+2d)
\\
&\leq 3 n_G C^{RR}(\|r^B -r^A \| +|t^A-t^B|)
\end{align*}
for all $s\in [t^A+d,T]$. Furthermore,
\begin{align*}
&\max_{s\in [t^A\vee t^B,t^A+d]}\|\zeta(\alpha^{t^A,a,r^A}(s,v),\xi^{t^A,a,r^A}(s,v))-\zeta(\alpha^{t^B,a,r^B}(s,\tilde v),\xi^{t^B,a,r^B}(s,\tilde v))\|
\\
&\leq \|\zeta(\alpha^{t^A,a,r^A}(t^A\vee t^B,v),\xi^{t^A,a,r^A}(t^A\vee t^B,v))-\zeta(\alpha^{t^B,a,r^B}(t^A\vee t^B,\tilde v),\xi^{t^B,a,r^B}(t^A\vee t^B,\tilde v))\|+2n_G C^{RR}d
\\
&\leq n_G C^{RR}(\|r^A-r^B\|+|t^A-t^B|+2d)
\\
&\leq 3 n_G C^{RR}(\|r^A-r^B\|+|t^A-t^B|)
\end{align*}
and \eqref{ekv:zetaBND} follows.\\

Finally we have,
\begin{equation*}
\sum_{\tilde\tau_j \leq T}K_{\tilde\beta_j}= \sum_{\tau_j \leq T-d}K_{\beta_j}\leq \sum_{\tau_j \leq T}K_{\beta_j},
\end{equation*}
for all $\omega\in\Omega$.\qed\\

\bigskip

In the proof of the next proposition we will use two well known results for solutions of SDEs. First, there is a constant $C_1^X>0$, such that
\begin{equation}\label{Xmax}
\E\left[\sup\limits_{s\in [t,T]}\|X_s^{t,x}\|\right]\leq C_1^X\left(1+\|x\|\right),
\end{equation}
for all $(t,x)\in\R^{n_L}$. Second, there is a constant $C_2^X>0$, such that,
\begin{equation}\label{Xspread}
\E\left[\sup\limits_{s\in [t^A\vee t^B,T]}\|X_s^{t^A,x^A}-X_s^{t^B,x^B}\|\right]\leq C_2^X\left\{\|x^A-x^B\|+(1+\|x^A\|\vee \|x^B\|)|t^A-t^B|^{1/2}\right\},
\end{equation}
for all $(t^A,x^A),(t^B,x^B)\in\R^{n_L}$.

\bigskip

\begin{prop}\label{Jprop}
There exists a constant $C^J>0$ such that, for each $t^A,t^B\in [0,T]$, $a\in\Gamma$, $r^A,r^B\in H(a)$ and for each $v:=(\tau_1,\ldots,\tau_{N};\beta_1,\ldots,\beta_{N})\in \mcV_{t^A,a,r^A}$, there is a control $\tilde v:=(\tilde\tau_1,\ldots,\tilde\tau_{\tilde N};\tilde\beta_1,\ldots,\tilde\beta_{\tilde N})\in \mcV_{t^B,a,r^B}$ such that,
\begin{align}\nonumber
J(t^B,x^B,a,r^B,\tilde v)\leq J(t^A,x^A,a,r^A,v)+&C^J(\|x^A-x^B\|+\|r^A-r^B\|
\\&\qquad+(1+\|x^A\|\vee \|x^B\|)|t^A-t^B|^{1/2}),\label{ekx:Jprop}
\end{align}
for all $x^A,x^B\in\R^{n_L}$.
\end{prop}

\noindent\emph{Proof.} Let $\tilde v$ be as in the proof of Lemma~\ref{zetaLEM}. Applying the results of Lemma~\ref{zetaLEM} and Proposition~\ref{fprop} combined with \eqref{Xmax}--\eqref{Xspread} we find
\begin{align*}
&J(t^B,x^B,a,r^B,\tilde v)- J(t^A,x^A,a,r^A,v)\\
&=\E\Bigg[\int_{t^B}^T f(X_s^{t^B,x^B},\alpha^{t^B,a,r^B}(s,\tilde v),\xi^{t^B,a,r^B}(s,\tilde v))ds + \sum_{\tilde \tau_j\leq T} K_{\tilde \beta_j}
\\
&\quad -\left(\int_{t^A}^T f(X_s^{t^A,x^A},\alpha^{t^A,a,r^A}(s,v),\xi^{t^A,a,r^A}(s,v))ds + \sum_{\tau_j\leq T} K_{\beta_j}\right)\Bigg]
\\
&\leq\E\Bigg[\int_{t^B}^{t^A\vee t^B} f(X_s^{t^B,x^B},\alpha^{t^B,a,r^B}(s,\tilde v),\xi^{t^B,a,r^B}(s,\tilde v))ds - \int_{t^A}^{t^A\vee t^B} f(X_s^{t^A,x^A},\alpha^{t^A,a,r^A}(s,v),\xi^{t^A,a,r^A}(s,v))ds
\\
&\quad\quad\quad \int_{t^A\vee t^B}^T \left\{f(X_s^{t^B,x^B},\alpha^{t^B,a,r^B}(s,\tilde v),\xi^{t^B,a,r^B}(s,\tilde v))
- f(X_s^{t^A,x^A},\alpha^{t^A,a,r^A}(s,v),\xi^{t^A,a,r^A}(s,v))\right\}ds\Bigg]
\\
&\leq C^f_1(1+C^X_1(1+\|x^A\|\vee\|x^B\|))|t^B-t^A|+T C_2^f\Big(\|r^A-r^B\|+|t^A-t^B|
\\
&\quad\quad\quad+C_2^X\left\{\|x^A-x^B\|+(1+\|x^A\|\vee \|x^B\|)|t^A-t^B|^{1/2}\right\}\Big)
\\
&\leq C^J\left(\|x^A-x^B\|+\|r^A-r^B\|+(1+\|x^A\|\vee \|x^B\|)|t^A-t^B|^{1/2}\right)
\end{align*}
with $C^J=\sqrt{T} C^f_1(1+C^X_1) + T C_2^f( 1\vee \sqrt{T}\vee C_2^X)$.\qed\\

\bigskip

\begin{prop}\label{VFprop}
There exists constants $C_1^V>0$ and $C_2^V>0$ for which the value function satisfies the following
\begin{align}\label{ekv:VFprop1}
|V(t,x,a,r)|&\leq C_1^V(1+\|x\|),\quad \forall (t,x,a,r)\in [0,T]\times \R^{n_L}\times \mathop{\cup}\limits_{a\in \Gamma}\{a\times H(a)\} \\
|V(t,x,a,r)-V(t',x',a,r')| &\leq C_2^V(\|x-x'\|+\|r-r'\|+(1+\|x\|\vee \|x'\|)|t-t'|^{1/2}),\nonumber\\
&\quad\quad\quad\forall (t,x,r),\,(t',x',r')\in [0,T]\times \R^{n_L}\times H(a),\:\forall a\in \Gamma.\label{ekv:VFprop2}
\end{align}
\end{prop}

\noindent\emph{Proof.} We have that
\begin{align*}
|V(t,x,a,r)|&\leq (T-t)\, \E\left[ \sup\limits_{s\in [t,T]} C^f_1(1+\|X^{t,x}_s\|)\right]\\
&\leq TC^f_1(1+C^X_1(1+\|x\|)),
\end{align*}
and \eqref{ekv:VFprop1} follows.The second inequality is an immediate consequence of Proposition~\ref{Jprop} and the proof is competed.\qed\\

\bigskip

\noindent\emph{Proof of principle of optimality (Theorem~\ref{Popt1}).} To prove the principle of optimality we define
\begin{align*}
\bar V(t,x,a,r):=\inf\limits_{v\in \mcV_{t,a,r}}\E\Bigg[\int_t^\nu  f(X^{t,x}_s,\alpha^{t,a,r}(s,v),\xi^{t,a,r}(s,v))ds + \sum\limits_{\tau_i \leq \nu}K_{\beta_i}\\
+ V(\nu,X^{t,x}_\nu,\alpha^{t,a,r}(\nu,v),\xi^{t,a,r}(\nu,v)) \Bigg].
\end{align*}
By \eqref{VFdef} and the definition of infimum there exists an $\varepsilon >0$ and a control $v^\varepsilon:=(\tau^\varepsilon_1,\ldots,\tau^\varepsilon_{N^\varepsilon};\beta^\varepsilon_1,\ldots,\beta^\varepsilon_{N^\varepsilon})$ such that
\begin{equation*}
V(t,x,a,r)+\varepsilon > J(t,x,a,r;v^\varepsilon).
\end{equation*}
By the strong Markov property of $(X_t,\alpha(t,v),\xi(t,v))$ corresponding to any control $v\in\mcV$ we have
\begin{align*}
J(t,x,a,r;v^\varepsilon)&=\E\Bigg[\int_t^\nu  f(X^{t,x}_s,\alpha^{t,a,r}(s,v^\varepsilon),\xi^{t,a,r}(s,v^\varepsilon))ds + \sum\limits_{t\leq \tau^\varepsilon_i \leq \nu}K_{\beta^\varepsilon_i}+
\\
&\quad\quad\quad\quad J(\nu,X^{t,x}_\nu,\alpha^{t,a,r}(\nu,v^\varepsilon),\xi^{t,a,r}(\nu,v^\varepsilon);v^\varepsilon) \Bigg]
\\
&\geq \E^{t,x,a,r}\Bigg[\int_t^\nu  f(X^{t,x}_s,\alpha^{t,a,r}(s,v^\varepsilon),\xi^{t,a,r}(s,v^\varepsilon))ds + \sum\limits_{t\leq\tau^\varepsilon_i \leq \nu}K_{\beta^\varepsilon_i}+
\\
&\quad\quad\quad\quad V(\nu,X^{t,x}_\nu,\alpha^{t,a,r}(\nu,v^\varepsilon),\xi^{t,a,r}(\nu,v^\varepsilon)) \Bigg]
\\
&\geq \bar V(t,x,a,r)
\end{align*}
Since $\varepsilon>0$ was arbitrary we have $\bar V(t,x,a,r)\leq V(t,x,a,r)$.\\

Conversely, for any $\varepsilon>0$ and any $v\in \mcV_{t'}$ there is, by Proposition~\ref{Jprop} and Proposition~\ref{VFprop}, a $\delta>0$ such that, whenever $\|x-x'\|+\|r-r'\|+|t-t'|<\delta$ we have,
\begin{equation}
|V(t,x,a,r)-V(t',x',a,r')|\leq \varepsilon,\quad \forall a\in\Gamma,
\end{equation}
and
\begin{equation}\label{Jdiff}
J(t,x,a,r;\tilde v)-J(t',x',a,r';v)\leq \varepsilon,\quad \forall a\in\Gamma,
\end{equation}
where $\tilde v$ is the control used in the proof of Proposition~\ref{Jprop}. Let $(D_j^a)_{j\geq 1}$ be a Borel partition of $[0,T]\times \R^{n_L}\times H(a)$ with diameter diam$(D_j^a)<\delta$. For each $a\in \Gamma$ and each subset in the corresponding partition choose $(t_j,x_j,r_j)\in D_j^a$. Then, there is a control $v_{a,j}^\varepsilon\in\mcV_{t_j}$ such that
\begin{equation*}
J(t_j,x_j,a,r_j;v_{a,j}^\varepsilon)\leq V(t_j,x_j,a,r_j)+\varepsilon
\end{equation*}
and for any $(t,x,r)\in D^a_j$ a corresponding control $\tilde v_{a,j}^\varepsilon$ for which
\begin{equation*}
J(t,x,a,r;\tilde v_{a,j}^\varepsilon)\leq J(t_j,x_j,a,r_j;v_{a,j}^\varepsilon)+\varepsilon.
\end{equation*}
Hence, we have
\begin{equation*}
J(t,x,a,r;\tilde v_{a,j}^\varepsilon)\leq J(t_j,x_j,a,r_j;v_{a,j}^\varepsilon)+\varepsilon\leq V(t_j,x_j,a,r_j)+2\varepsilon \leq V(t,x,a,r)+3\varepsilon.
\end{equation*}
Now let $\hat v=(\tau_1,\ldots,\tau_{N(\nu)},\tilde \tau^\varepsilon_{a',j,1}\ldots \tilde \tau^\varepsilon_{a',j,\tilde N(\nu)};\beta_1,\ldots,\beta_{N(\nu)},\tilde \beta^\varepsilon_{a',j,1}\ldots \tilde \beta^\varepsilon_{a',j,\tilde N(\nu)})\in \mcV_t$, where $t\leq \tau_1\leq\cdots\leq\tau_{N(\nu)}\leq \nu \leq \tilde\tau^\varepsilon_{a',j,1}\leq\cdots\leq\tilde \tau^\varepsilon_{a',j,\tilde N(\nu)}\leq T$ and $a'=\alpha^{t,a,r}(\nu,\hat v)$ and $j$ is such that $(\nu,X_{\nu}^{t,x},\xi^{t,a,r}(\nu,\hat v))\in D_j^{a'}$. Then,
\begin{align*}
&V(t,x,a,r)\leq J(t,x,a,r;\hat v)\\
&=\E\left[\int_t^T  f(X_s^{t,x},\alpha^{t,a,r}(s,\hat v),\xi^{t,a,r}(s,\hat v))ds + \sum\limits_{\hat \tau_i\leq T}K_{\hat \beta_i}\right]
\\
&=\E\Bigg[\int_t^\nu  f(X_s^{t,x},\alpha^{t,a,r}(s,v),\xi^{t,a,r}(s,v))ds + \sum\limits_{t\leq \tau_i\leq \nu}K_{\beta_i}
\\
&\quad +\E\left[\int_\nu^T  f(X^{\nu,X_\nu^{t,x}}_s,\alpha^{\nu,\alpha^{t,a,r}(\nu,v),\xi^{t,a,r}(\nu, v)}(s,\tilde v^\varepsilon_{a',j}),\xi^{\nu,\alpha^{t,a,r}(\nu,v),\xi^{t,a,r}(\nu, v)}(s,\tilde v^\varepsilon_{a',j}))ds + \sum\limits_{\nu\leq \tilde \tau^\varepsilon_{a',j,i}\leq T}K_{\tilde \beta^\varepsilon_{a',j,i}}\right] \Bigg]
\\
&=\E\left[\int_t^\nu  f(X^{t,x}_s,\alpha^{t,a,r}(s,v),\xi^{t,a,r}(s,v))ds + \sum\limits_{t\leq \tau_i\leq \nu}K_{ \beta_i}+J(\nu,X_\nu,\alpha^{t,a,r}(\nu,v),\xi^{t,a,r}(\nu,v),\tilde v^\varepsilon_{a',j})\right]\\
&\leq \E^{t,x,a,r}\left[\int_t^\nu  f(X_s^{t,x},\alpha^{t,a,r}(s,v),\xi^{t,a,r}(s,v))ds + \sum\limits_{t\leq \tau_i\leq \nu}K_{ \beta_i}+V(\nu,X_\nu,\alpha^{t,a,r}(\nu,v),\xi^{t,a,r}(\nu,v))\right]+3\varepsilon.
\end{align*}
By taking the infimum over all controls $v=(\tau_1,\ldots,\tau_{N(\nu)};\beta_1,\ldots,\beta_{N(\nu)})\in \mcV_{t,a,r}$, with $t\leq \tau_1\leq\cdots\leq\tau_{N(\nu)}\leq \nu$ we get $V(t,x,a,r)\leq \bar V(t,x,a,r) +3\varepsilon$. Since $\varepsilon>0$ was arbitrary, the result follows.\qed\\

\bigskip

As mentioned above, Theorem~\ref{Popt1} has an important implication for building algorithms. This will be detailed in the following corollary. But first we need a few new definitions.

\begin{defn}
For each $a\in\Gamma$ and each $\beta\in \Lambda(a)$ we define the sum $a+\beta$, such that $(a+\beta)_k:=a_k$, for $k\neq \beta^\Node$ and $(a+\beta)_{\beta^\Node}:=(\beta^f,\beta^t)$.
\end{defn}

\begin{defn}
For each $t\in [0,T]$, each $a\in \Gamma$ and $r\in H(a)$ we let $A_{s}^{t,a,r}:=\alpha^{t,a,r}(s,\emptyset)$, and $R_{s}^{t,a,r}:=\xi^{t,a,r}(s,\emptyset)$, for all $s\in [t,T]$.
\end{defn}

\begin{defn}
Given $a\in\Gamma$ and $r\in H(a)$, we let $\psi(a,r):=\max_{j\in\Node_G}(\delta_{j,a_j}-r_j)$ and $\psi(a):=\psi(a,0)$.
\end{defn}

We are now ready to state the corollary which implies that Problem 1 can be written as an iterated optimal stopping problem.
\begin{cor}
For each $t\in[0,T]$, $x\in\R^{n_L}$, $a\in\Gamma$ and $r\in H(a)$, we have
\begin{align}\nonumber
V(t,x,a,r)=\inf\limits_{\tau \in \mcT_{t,\psi(a,r)}}\E^{t,x}\Bigg[\int_t^{\tau \wedge T}   f(X_s,a,R_{s}^{t,a,r})ds + \min\limits_{\beta\in\overline \Lambda(a)}\{K_\beta+V(\tau,X_\tau,A_{\tau}^{t,a,r}+\beta, R_{\tau}^{t,a,r})\} \Bigg],
\label{dpALG}\end{align}
where $\mcT_{t,\psi(a,r)}$ is the subset of $\tau\in\mcT_t$ with $\tau\leq t+\psi(a,r)$ and $\overline \Lambda(a):=\Lambda(a)\cup\emptyset$.
\end{cor}

\noindent\emph{Proof.} This follows from the fact that $\mcT_{t,\psi(a,r)}$ is a subset of $\mcT_{t}$.\qed\\

\section{Bounds on the value function\label{Sec:bnds}}
The numerical algorithm that we propose in this article will be facilitated by upper and lower bounds on the value function. To provide adequate such bounds we solve two simplified versions of Problem~1.

The key concept of the simplified problems is that the dimensionality is reduced to the extent that conventional numerical methods, such as Markov-Chain approximation algorithms\cite{numSCbok}, become tractable. Compared to the regression Monte Carlo approach, that will be used to solve Problem 1, the results produced by Markov-Chain approximations can be considered exact.

\subsection{A lower bound}
The lower bound will be obtained by making an outer approximation of the feasible regions following each of the specified contingencies (including the base case) and representing the transition progression vector by its sum. Combined with a relaxation on the set of admissible controls this will lead to a two-dimensional switching problem that can be solved using conventional numerical methods.\\

\subsubsection{Running cost approximation}
We first pick a normalized vector $\eta_L\in \R^{n_L}_+$ and $N_{\eta_G}\geq 1$ vectors $\eta^j_G\in \R^{n_G}$, with $j\in \mcJ:=\{1,\ldots,N_{\eta_G}\}$. From these we define $\eta^j:=(\eta_G^{j},\eta_L)\in\R^{n_G}\times\R^{n_L}_+$. Since $G_i$ is non-empty and compact, the maximal values, $d_{i,j}^\eta$, to the non-linear optimization problems
\begin{subequations}
\begin{align}
\max\limits_{p_G\in D_G,\:p_D\in \R^{n_L}_+} \:\:\quad (\eta_G^j)^\top p_G + \eta_L^\top p_D \\
\text{subj. to} \quad\quad\quad\quad\:  (p_G,p_D)\in G_i
\end{align}
\end{subequations}
are well defined, for $i=0,1,\ldots,n_c$, for all $j\in\mcJ$. Now, each of the sets
\begin{equation}\label{ekv:GetaDEF}
G_{i,j}^\eta:=\{(p_G,p_D)\in D_G\times\R^{n_L}:\, (\eta_G^j)^\top p_G + \eta_L^\top p_D \leq d_{i,j}^\eta\}
\end{equation}
are outer approximations of $G_i$ ($G_i\subset G_{i,j}^\eta$).

For all $t\in[0,T]$, we define the process $X^\eta_t=\eta_L^\top X_t$. We note that $X^\eta_t=\eta_L^\top m(t)+Z^\eta_t$ is a Markov process, with
\begin{align*}
dZ^\eta_t&=-\alpha dZ^\eta_t dt+\sigma^\eta d W^\eta_t,
\\
Z^\eta_0&=\eta_L^\top(x_0-m(0)),
\end{align*}
where $\sigma^\eta:=\|\eta^\top_L\sigma\|$ and $W^\eta_t:=\frac{1}{\sigma^\eta}\eta^\top_L\sigma W_t$, for all $t\in[0,T]$. Let $\mcF^\eta_t$ be the filtration generated by $(W^\eta_t,0\leq t\leq T)$. 

Similarly, for all $j\in\mcJ$, $a\in \Gamma$ and $r\in H(a)$, we define $\zeta^\eta_j(a,r):=(\eta^j_G)^\top\zeta(a,r)$. By \eqref{ekv:GetaDEF} we thus have $(\zeta(a,r),X_t)\in G_{i,j}^\eta$ iff $X^\eta_t+\zeta^\eta_j(a,r)\leq d_{i,j}^\eta$. This leads us to the definition $z_{i,j}^\eta(y,x):=\mathop{\arg\min}\limits_{z\in \R^{n_L}_+}\{ c^\top z : y+x-\eta^\top_L z\leq d_{i,j}^\eta\}$, which gives
\begin{align}\label{zETAdef}
c^\top z_{i,j}^\eta(x,y) =\left[(y+x-d_{i,j}^\eta)\vee 0\right]\min\limits_{k\in\Node_G}\frac{c_k}{\eta_k}.
\end{align}

For each $y\in\R$ we define the running cost, corresponding to the approximated region of feasibility, as
\begin{equation*}
f^\eta_1(y,a,r):=C_G(\zeta(a,r)) + \max_{j\in \mcJ}\sum_{i=0}^{n_c}q_i c^\top z_{i,j}^\eta(\zeta_j^\eta(a,r),y).
\end{equation*}

Having introduced the above approximations we find that if we replace $f$ in Problem~1 with $f^\eta_1$ we can reduce the dimensionality of the system operators planning problem by only observing $X^\eta_t$ instead of the full demand process. To obtain a further reduction we introduce the new variable $\rho$ taking values in $[0,\varsigma_{\max})$, with $\varsigma_{\max}:=\max_{a\in\Gamma}\varsigma(a)$, where $\varsigma(a):=\sum_{k=1}^{n_G}\delta_{(k,a_k)}$, that will represent the sum of the transition progression vector and define a new running cost that is dominated by $f^\eta_1$ and, thus, also by the actual running cost. We will then propose to use the solution of an optimal switching problem in the two-dimensional $(X^\eta_t,\rho(t))$--space as a lower bound to the value function in Problem 1.\\

\begin{defn}
For each $y\in\R$, $a\in\Gamma$ and $\rho\in [0,\varsigma(a))$, we define the running cost,
\begin{equation*}
f^\eta_2(y,a,\rho):=\mathop{\inf}\limits_{r\in H(a)}\{f^\eta_1(y,a,r):\,S(r)=\rho\}.
\end{equation*}
where $S(r)=\sum_{k=1}^{n_G} r_k$.
\end{defn}
Then, $f^\eta_2$ is convex and monotonically increasing in $y$. Furthermore, it satisfies the conditions of the following proposition.
\begin{prop}
For all $a\in\Gamma$,
\begin{align}\label{ekv:fETApropDOM}
f^\eta_2(\eta _L^\top x,a,S(r))&\leq f(x,a,r),\quad \forall (x,r)\in \R^{n_L}\times  H(a).
\end{align}
Furthermore, there are constants $C_1^{f,\eta}>0$ and $C_2^{f,\eta}>0$ such that,
\begin{align}\label{ekv:fETAprop1}
|f^\eta_2(y,a,w)| &\leq C_1^{f,\eta}(1+|y|),\quad \forall (y,w)\in \R\times  [0,\varsigma(a))
\end{align}
and
\begin{align}\label{ekv:fETAprop2}
|f^\eta_2(y^A,a,w)-f^\eta_2(y^B,a,w)| &\leq C_2^{f,\eta}|y^A-y^B|, \quad \forall y^A,\,y^B\in \R\text{ and }w\in[0,\varsigma(a)).
\end{align}
\end{prop}

\noindent\emph{Proof.} First we note that, for $i=0,\ldots,n_c$, for all $j\in\mcJ$, we have
\begin{equation*}
c^\top z^\eta_{i,j}((\eta_G^j)^\top p_G,\eta_L^\top x)\leq c^\top z^\eta_{i,j}((\eta_G^j)^\top p_G,\eta_L^\top (x\vee 0))\leq c^\top z^*_{i}( p_G,x\vee 0),\quad \forall\, (p_G,x) \in D_G\times\R^{n_L},
\end{equation*}
where the first inequality follows from the fact that $\eta_L\geq 0$ and the second follows from $G_{i}\subset G_{i,j}$. Now, for each $a\in\Gamma$ and each $\rho\in [0,\varsigma(a))$, we take the minimum over all transition progression vectors with $\sum r_k =\rho$ and \eqref{ekv:fETApropDOM} follows.

Inequality \eqref{ekv:fETAprop1} follows from noting that
\begin{align*}
|f^\eta_1(y,a,r)|\leq \max_{p_G\in D_G}|C_G(p_G)|+\sum_{i=0}^{n_c}q_i|y|
\end{align*}
and \eqref{ekv:fETAprop2} follows immediately from \eqref{zETAdef} by noting that
\begin{align*}
|f^\eta_2(y^A,a,w)-f^\eta_2(y^B,a,w)|&\leq \sum_{i=0}^{n_c}q_i\max_{j\in\mcJ}\max_{p_G\in D_G}|c^\top z_{i,j}^\eta((\eta_G^j)^\top p_G,y^A)-c^\top z_{i,j}^\eta((\eta_G^j)^\top p_G,y^B)|
\\
&\leq \sum_{i=0}^{n_c}q_i |y^A-y^B|.
\end{align*}
which finishes the proof.\qed\\

\subsubsection{Admissible controls}
Using the scalar $\rho$ to replace the $n_G$--dimensional vector $r$ causes problems when trying to define the set of admissible switching requests, as $A^{t,a,r}_s$ and $R^{t,a,r}_s$ generally cannot be uniquely defined as a functions of $S(r)$, for all $a\in\Gamma$. Instead the system operator in the lower bound problem will be given the opportunity to (under suitable conditions) decide when a transition is completed, so that she may continuously choose the vector $r$, based on $S(r)= \rho$. This gives a relaxation on the set of admissible controls that allows us to further reduce the dimension of the state space. First we define the number of transiting states.
\begin{defn}
For each $a\in\Gamma$, the number of transiting states is defined as $N_T(a):=\sum_{k=1}^{n_G} 1_{\Lambda_k}(a_k)$. Furthermore, we let $N_T^0:=N_T(\alpha_0)$ and define $N_T^{0,k}:=1_{\Lambda_k}(a_k)$, $\forall k\in\Node_G$.
\end{defn}
An admissible control $v:=\{\tau_1,\ldots,\tau_N;\beta_1,\ldots,\beta_N;\vartheta_{-N_T^0+1},\ldots,\vartheta_{0},\vartheta_1,\ldots,\vartheta_N\}$ for the lower bound problem will be a sequence of $\mcF^\eta_t$--stopping times $0\leq\tau_1\leq\cdots\leq \tau_N\leq T$ and corresponding $\mcF^\eta_{\tau_j}$--measurable switching requests $\beta_j\in\Lambda(\alpha^\eta(\tau_j^-))$. What is new in the lower bound problem is that we introduce the execution times $\vartheta_{-N_T^0+1},\ldots,\vartheta_{0}$ and $\vartheta_1,\ldots,\vartheta_N$, such that $\vartheta_j$ dictate when the transition initiated at time $\tau_j$ is completed (where we assume that a sorting is made of the states that are in transition at the start of the operating period and let $j\leq 0$ corresponds to a execution times of these transitions). As for Problem~1 we define the operating mode process through the sequence\\ $v^k:=(\tau_1^k,\ldots,\tau_{N_k}^k;\beta_{k,1},\ldots,\beta_{k,N_k};\vartheta_{-N_T^{0,k}+1},\ldots,\vartheta_{N_k}^k)$, this time as
\begin{align}
\alpha^\eta_k(t,v)=\alpha_{0,k} 1_{[0,\vartheta^k_0)}(t)+\alpha_{0,k}^t 1_{[\vartheta^k_0,\tau_1^k)}(t)+\sum_{j=1}^{N_k}\left((\beta_{k,j}^f, \beta_{k,j}^t) 1_{[\tau_j^k,\vartheta^k_j)}(t)+\beta_{k,j}^t 1_{[\vartheta^k_j,\tau^k_{j+1})}(t)\right),\label{alphaETA_DEF}
\end{align}
with $\vartheta^k_0=0$ whenever $N_T^{0,k}=0$. From \eqref{alphaETA_DEF} we conclude that $\vartheta_j$ plays the role that $\tau_j+\delta_{\beta_j}$ has in Problem~1.

To keep track of the progression of transitions we introduce the total transition progression process
\begin{align}
\rho(t,v)=\sum_{k=1}^{n_G}\xi_{0,k}+\int_{0}^t N_T(\alpha^\eta(s,v))ds-\sum_{\vartheta_{j} \leq t}\delta_{\beta_j}\label{rho_DEF}
\end{align}

Note that both $\alpha^\eta$ and $\rho$ all $\mcF_t^\eta$--adapted.
\begin{defn}
The set of all admissible controls for the lower bound problem, denoted $\mcV^{\eta,R}$, are all the controls $v:=\{\tau_1,\ldots,\tau_N;\beta_1,\ldots,\beta_N;\vartheta_{-N_T^0+1},\ldots,\vartheta_{0},\vartheta_1,\ldots,\vartheta_N\}$ as specified above with
\begin{equation}\label{execBDS}
0\leq \rho(t,v)\leq \varsigma (\alpha^\eta(t,v)),\quad \forall t\in[0,T].
\end{equation}
\end{defn}

Another way of writing \eqref{rho_DEF} is\footnote{See proof of Proposition~\ref{lemLB}.}
\begin{equation}\label{rho_DEFalt}
\rho(t,v)=\sum_{\tau_j\leq T}1_{[\tau_j,\vartheta_j)}(t)(t-\tau_j)
\end{equation}
Assume that $n_G=1$, then by plugging \eqref{rho_DEFalt} into \eqref{execBDS} we find that $\vartheta_j=\tau_j+\delta_{\beta_j}$ and the set of admissible controls reduces to the subset of $\mcF^\eta_t$--adapted admissible controls for Problem~1. When the system has more generators the situation is rather different and the admissible controls can be seen as a relaxation of the set of $\mcF^\eta_t$--adapted controls in $\mcV$.

\subsubsection{The lower bound problem}
\begin{defn}
For $t\in [0,T]$ and $y\in \R$ we let $X_s^{\eta,t,y}=X^\eta_s|X^\eta_t=y$, for all $s\in [t,T]$. Furthermore, we let $\mcV_t^{\eta,R}$ be the set of all admissible controls $v\in\mcV^{\eta,R}$ with $\tau_1\geq t$ and $\vartheta_j\geq t$, for $j=-N_T(\alpha^\eta(t^-,v))+1,\ldots,N$. We define $\alpha^{\eta,t,a}(s,v)$ and $\rho^{t,a,w}(s,v)$, for all $v\in \mcV_t^{\eta,R}$, $a\in\Gamma$ and $w\in [0,\varsigma(a))$, as
\begin{equation}\label{alphaETAdef2}
\alpha^{\eta,t,a}_k(s,v):=a_{k} 1_{[0,\vartheta^k_0)}(t)+a_{k}^t 1_{[\vartheta^k_0,\tau_1^k)}(t)\sum\limits_{j=0}^{N_k}\left((\beta^f_{k,j},\beta^t_{k,j})1_{[\tau^k_j,\vartheta^k_j)}(t) +\beta_{k,j}^t 1_{[\vartheta^k_j,\tau^k_{j+1})}(t)\right),
\end{equation}
with $\vartheta^k_0=0$ whenever $a_k\notin\Lambda_k$; and
\begin{equation}\label{rhoDEF2}
\rho^{t,a,w}(s,v)=w+\int_{t}^s N_T(\alpha^{\eta,t,a}(u,v))du-\sum_{t\leq \vartheta_{j} \leq s}\delta_{\beta_j}
\end{equation}
Finally, we let $\mcV_{t,a,r}^{\eta,R}$ be the set of controls $v\in\mcV_t^{\eta,R}$ that are admissible for $\alpha^\eta(t^-,v)=a$ and $\rho(t^-,v)=w$.
\end{defn}

For each $t\in[0,T]$, $y\in\R$ and $a\in \Gamma$, the expected operating cost-to-go associated with the control $v\in \mcV^{\eta,R}_t$ is
\begin{equation*}
J^\eta(t,y,a,w;v):=\E\left[\int_t^T f^\eta_2(X_s^{\eta,t,y},\alpha^{\eta,t,a}(s),\rho^{t,a,w}(s))ds+ \sum\limits_{t<\tau_i\leq T}K_{\beta_i}\right].
\end{equation*}
To find the lower bound corresponding to $\eta_L\in\R^{n_L}_+$ we need to solve the following problem.\\

\noindent\textbf{Problem LB.} For all $t\in [0,T]$, all $a\in\Gamma$ and all $x\in \R$, find
\begin{equation*}
V^{\eta}(t,y,a,w)=\inf\limits_{v\in\mcV^{\eta,R}_{t,a,w}} J^\eta(t,y,a,w;v).
\end{equation*}
for all $w\in [0,\varsigma (a))$. \qed\\
\bigskip

The remainder of this Section will be devoted to showing that the value function of Problem LB can be used to define a lower bound for the value function of Problem 1.

\begin{defn}
Let $\mcV^\eta\subset\mcV$ be the subset of all $\mcF^\eta_t$--adapted controls.
\end{defn}

The following lemma specifies that $\mcV^{\eta,R}$ is a relaxation of $\mcV^\eta$ in the sense that there is an injective map from $\mcV^\eta$ to $\mcV^{\eta,R}$ and proves that the cost-to-go is lower for Problem LB.
\begin{lem}\label{lemLB}
Given $t\in [0,T]$, $a\in\Gamma$ and $r\in H(a)$. For each $v=(\tau_1,\ldots,\tau_N;\beta_1,\ldots,\beta_N)\in\mcV^\eta_{t,a,r}$ there is a sequence of $\mcF^\eta_t$--stopping times $\vartheta_{-N_T(a)+1},\ldots,\vartheta_N$, such that $v^R:=(\tau_1,\ldots,\tau_N;\beta_1,\ldots,\beta_N;\\ \vartheta_{-N_T(a)+1},\ldots,\vartheta_N)\in \mcV^{\eta,R}_{t,a,S(r)}$ and
\begin{equation*}
J^\eta(t,\eta_L^\top x,a,S(r);v^R)\leq J(t,x,a,r;v),
\end{equation*}
for all $(t,x,a,r)\in [0,T]\times \R^{n_L}_+ \times \cup_{a\in \Gamma}\{a\times H(a)\}$.
\end{lem}

\noindent\emph{Proof.} Let $\kappa_1,\ldots,\kappa_{N_T(a)}$ be an enumeration of the generators that are in transition in operation mode $a$ and let $\beta^0_{j}:=(\kappa_j,\alpha_{0,\kappa_j})$ be the corresponding switching requests, for $j=1,\ldots,N_T(a)$.

Pick $\vartheta_{-N_T(a)+j}=t+\delta_{\beta^0_{j}}-r_{\kappa_j}$, for $j=1,\ldots,N_T(a)$, and $\vartheta_j=\tau_j+\delta_{\beta_j}$, for $j=1,\ldots,N$. Then $\vartheta_{-N_T(a)+1},\ldots,\vartheta_{0}$ are $\mcF^\eta_0$--measurable and $\vartheta_1,\ldots,\vartheta_N$ are $\mcF^\eta_t$--stopping times. By definition we have $\alpha^{t,a,r}(s,v)=\alpha^{\eta,t,a}(s,v^R)$, for all $s\in[t,T]$.

We will show that \eqref{execBDS} holds by proving that \eqref{rho_DEFalt} is in fact equivalent to \eqref{rho_DEF}. We let
\begin{equation*}
\tilde\rho(s,v^R)=\sum_{j=-N_T(a)+1}^N 1_{[\tau_j,\vartheta_j)}(t)(s-\tau_j)
\end{equation*}
where $\tau_{-N_T(a)+j+1}:=t-r_{\kappa_j}$, for $j=1,\ldots,N_T^0$. We have $\tilde\rho(t^-,v^R)=\sum_{j=1}^{N_T(a)}r_{\kappa_j}=S(r)$. Now, for all $s\in [t,T]$ for which $s\notin \{\vartheta_{-N_T(a)+1},\ldots,\vartheta_{N}\}$ we have $\frac{d}{dt} \tilde\rho(s,v^R)=N_T(\alpha^\eta(s,v^R))$ and for $s\in \{\vartheta_{-N_T(a)+1},\ldots,\vartheta_{N}\}$ we have $\tilde\rho(s,v^R)-\tilde\rho(s^-,v^R)=\sum_{\vartheta_j=s}\vartheta_j$. Hence, $\tilde\rho(s,v^R)=\rho^{t,a,S(r)}(s,v^R)$, for all $s\in [t,T]$. By \eqref{xiDEF2} it now follows that $\rho^{t,a,S(r)}(s,v^R)=S(\xi^{t,a,r}(s,v))$. Hence, \eqref{execBDS} holds for $v^R$ and we conclude that $v^R\in\mcV_{t,a,S(r)}^{\eta,R}$.

For the second part we note that, for all $(t,x,a,r)\in [0,T]\times \R^{n_L}_+ \times \cup_{a\in \Gamma}\{a\times H(a)\}$, and all $v\in\mcV^\eta$,
\begin{align*}
J(t,x,a,r;v)&=\E\left[\int_t^T  f(X_s^{t,x},\alpha^{t,a,r}(s,v),\xi^{t,a,r}(s,v))ds + \sum\limits_{\tau_i\leq T}K_{\beta_i}\right]\\
&\geq \E\left[\int_t^T  f_2^\eta (\eta_L^\top X_s^{t,x},\alpha^{t,a,r}(s,v),S(\xi^{t,a,r}(s,v)))ds + \sum\limits_{\tau_i\leq T}K_{\beta_i}\right]\\
&= \E \left[\int_t^T  f_2^\eta ( X^{\eta,t,\eta_L^\top x}_s,\alpha^{\eta,t,a}(s,v^R),\rho^{t,a,S(r)}(t,v^R))ds + \sum\limits_{\tau_i\leq T}K_{\beta_i}\right]\\
&= J^\eta(t,\eta_L^\top x,a,S(r);v^R),
\end{align*}
which finishes the proof.\qed\\

\bigskip

This leads us to the following theorem and the subsequent corollary which are the main results of this section.

\begin{thm}
For each $t\in[0,T]$, $x\in\R^{n_L}$ and $a\in \Gamma$ we have
\begin{equation}
V^{\eta}(t,\eta_L^\top x,a,S(r))\leq V(t,x,a,r),
\end{equation}
$\forall r\in H(a)$.
\end{thm}

\noindent\emph{Proof.} Follows from Lemma~\ref{lemLB} by noting that
\begin{align*}
\inf\limits_{v\in\mcV^{\eta,R}_{t,a,S(r)}} J^\eta(t,y,a,w;v)\leq \inf\limits_{v\in\mcV^{\eta}_{t,a,r}} J^\eta(t,y,a,w;v^R)\leq \inf\limits_{v\in \mcV_{t,a,r}} J(t,x,a,r;v),
\end{align*}
with $y=\eta_L^\top x$ and $w=S(r)$.\qed\\

\begin{cor}
Let $U\subset \{\eta_L\in S^{n_L-1}:\: \eta_L>0\}$ be a compact set. For each $t\in[0,T]$, $x\in\R^{n_L}$ and $a\in \Gamma$ we have
\begin{equation}
V^{\rm{LB}}(t,x,a,r):=\min_{\eta\in U} V^{\eta}(t,\eta_L^\top x,a,S(r))\leq V(t,x,a,r),
\end{equation}
$\forall r\in H(a)$.
\end{cor}

To improve the tightness of the lower bound we can thus pick a number of vectors $\eta_L$ and use the maximum of the corresponding lower bounds as a lower bound for the value function of Problem 1.

\subsubsection{Dynamic programming algorithm}
We will now show that the principle of optimality holds for the value function of Problem LB. This is a slightly more complicated matter that in the setting of Problem 1 since we may be forced to take actions to satisfy \eqref{execBDS}.

We will need the additional well known fact that tells us that there is a constant $C^{X^\eta}_3>0$ such that, for any $t\in [0,T]$ and any $d\in (0,T-t]$, we have
\begin{align*}
\E\left[\sup_{s\in [t,\,T-d]}|X^{\eta,t,x}_{s}-X^{\eta,t,x}_{s+d}|\right]&\leq C^{X^\eta}_3(1+\|x\|)\sqrt{d}.
\end{align*}
We start by showing that an equivalent of Proposition~\ref{Jprop} holds for the cost-to-go in Problem LB.
\begin{prop}\label{JETApropAB}
Given $a\in \Gamma$ there exists a constant $C^{J,\eta}>0$ such that, for any $t^A,\,t^B \in [0,T]$ and $w^A,w^B\in [0,\varsigma(a))$, and for each $v:=(\tau_1,\ldots,\tau_{N};\beta_1,\ldots,\beta_{N};\vartheta_{-N_T(a)+1},\ldots,\vartheta_{N})\in \mcV^{\eta,R}_{t^A,a,w^A}$ there is a control $\tilde v:=(\tilde\tau_1,\ldots,\tilde\tau_{\tilde N};\tilde\beta_1,\ldots,\tilde\beta_{\tilde N}; \tilde\vartheta_{-N_T(a)+1},\ldots,\tilde\vartheta_{\tilde N})\in \mcV^{\eta,R}_{t^B,a,w^B}$ such that,
\begin{align}\nonumber
J^\eta(t^B,y^B,a,w^B;\tilde v)\leq J&^\eta(t^A,y^A,a,w^A;v)
\\& + C^{J,\eta}\{|y^A-y^B|+(1+|y^A|+|y^B|)(|w^A-w^B|^{1/2}+|t^A-t^B|^{1/2})\}, \label{ekv:JETAprop}
\end{align}
for all $y^A,y^B\in\R$.
\end{prop}

\noindent\emph{Proof.} The proof will be divided into three cases depending on the relation between $w^A$ and $w^B$, and on the sign of $d^R:=\frac{w^A-w^B}{N_T(a)}+t^B-t^A\leq |w^A-w^B|+|t^A-t^B|$, where we use the convention that $\frac{w^A-w^B}{N_T(a)}=0$ whenever $N_T(a)=0$.\\

Assume first that $w^A\geq w^B$ and that $d^R\geq 0$. Then $\tilde v$ can be defined to trail behind $v$ as in the proof in Proposition~\ref{Jprop}. Define the control $\tilde v$ as $\tilde \tau_j:=\tau_j+d^R$ and $\tilde \beta_j:=\beta_j$, for $j=1,\ldots,\tilde N$, where $\tilde N:=\max\{j\in\{0,\ldots,N\}:\: \tau_j+d^R\leq T\}$ and $\tilde \vartheta_j:=(\vartheta_j+d^R)$, for $j=-N_T(a)+1,\ldots,\tilde N$. Since $w^A\geq w^B$ we have $\tilde \tau_1\geq t^A+d^R=\frac{w^A-w^B}{N_T(a)}+t^B\geq t^B$, and similar for $\tilde \vartheta_j$, for $j=-N_T(a)+1,\ldots,\tilde N$. Furthermore, $\rho^{t^B,a,w^B}(\tilde v,t)\in [w^B,w^A]$ for all $(t,\omega)\in [t^B,\,t^A+d^R]\times \Omega$ and $\rho^{t^B,a,w^B}(\tilde v,t+d^R)=\rho^{t^A,a,w^A}(v,t)$ for all $(t,\omega)\in [t^A+d^R,\,T]\times \Omega$, hence $\tilde v\in \mcV^{\eta,R}_{t^B,a,w^A}$.\\

We will use the following inequality
\begin{align*}
\E\left[\sup_{s\in [t^A,\,T-d^R]}|X^{\eta,t^A,y^A}_{s}-X^{\eta,t^B,y^B}_{s+d^R}|\right]&\leq \E\Bigg[\sup_{s\in [t^A,\,T-d^R]}|X^{\eta,t^A,y^A}_{s}-X^{\eta,t^A,y^A}_{s+d^R}|
\\
&\quad+\sup_{s\in [t^A,\,T-d^R]}|X^{\eta,t^A,y^A}_{s+d^R}-X^{\eta,t^B,y^B}_{s+d^R}|\Bigg]
\\
&\leq C^{X^\eta}_3(1+|y^A|)\sqrt{d^R}+C^{X^\eta}_2(|y^A-y^B|
\\
&\quad+(1+|y^A|\vee |y^B|)|t^A-t^B|^{1/2})\nonumber
\\
&\leq C^{X^\eta}_4\left(|y^A-y^B|+(1+|y^A|\vee |y^B|)(\sqrt{d^R}+|t^A-t^B|^{1/2})\right),
\end{align*}
with $C^{X^\eta}_4:=C^{X^\eta}_2\vee C^{X^\eta}_3$.\\

For the switching costs we have,
\begin{equation}\label{ekv:swCOSTequiv}
\sum_{t^B\leq \tilde\tau_j \leq T}K_{\tilde\beta_j}=\sum_{t^A\leq \tau_j \leq T-d^R}K_{\beta_j} \leq \sum_{t^A\leq \tau_j \leq T}K_{\beta_j}.
\end{equation}
Hence,
\begin{align*}
&J^\eta(t^B,y^B,a,w^B;\tilde v)- J^\eta(t^A,y^A,a,w^A;v)
\\
&\leq E\left[\int_{t^B}^T f^\eta_2(X_s^{\eta,t^B,y^B},\alpha^{\eta,t^B,a}(s,\tilde v),\xi^{t^B,a,w^B}(s,\tilde v))ds - \int_{t^A}^T f^\eta_2(X_s^{\eta,t^A,y^A},\alpha^{\eta,t^A,a}(s,v),\rho^{t^A,a,w^A}(s,v))ds\right]
\\
&=E\Bigg[\int_{t^B}^{t^A+d^R} f^\eta_2(X_s^{\eta,t^B,y^B},\alpha^{\eta,t^B,a}(s,\tilde v),\rho^{t^B,a,w^B}(s,\tilde v))ds\\
&\quad+\int_{t^A}^{T-d^R} \{f^\eta_2(X_{s+d^R}^{\eta,t^B,y^B},\alpha^{\eta,t^A,a}(s,v),\rho^{t^A,a,w^A}(s,v)-f^\eta_2(X_s^{\eta,t^A,y^A},\alpha^{\eta,t^A,a}(s,v),\rho^{t^A,a,w^A}(s,v)) )\}ds
\\
&\quad+\int_{T-d^R}^{T} f^\eta_2(X_s^{\eta,t^A,y^A},\alpha^{\eta,t^A,a}(s,v),\rho^{t^A,a,w^A}(s,v))ds\Bigg]
\\
&\leq C^{f,\eta}_1(1+C^{X^\eta}_1(1+|y^B|))(t^A+d^R-t^B) + C^{f,\eta}_2C^{X^\eta}_4\left(|y^A-y^B|+(1+|y^A|\vee |y^B|)(\sqrt{d^R}+|t^A-t^B|^{1/2})\right)T
\\
&\quad + C^{f,\eta}_1(1+C^{X^\eta}_1(1+|y^A|))d^R
\\
&\leq C^{J,\eta}_1 \left\{|y^A-y^B|+(1+|y^A|+|y^B|)(|w^A-w^B|^{1/2}+|t^A-t^B|^{1/2})\right\},
\end{align*}
with $C^{J,\eta}_1:=2 C^{f,\eta}_2C^{X^\eta}_4 T +3C_1^{f,\eta}(1+ C_1^{X^\eta})(\sqrt{\varsigma_{\max}}\vee\sqrt{T})$. This proves the proposition for $w^A\geq w^B$ and $d^R\geq 0$.\\

\bigskip
\bigskip

Now, assume instead that $w^A\geq w^B$ and $d^R<0$. Let $t^C:=t^B+\frac{w^A-w^B}{N_T(a)}<t^A$, then with $\tilde \tau_1\geq t^C$ and $\tilde\vartheta_j\geq t^C$, for $j=-N_T(a)+1,\ldots,\tilde N$, we get $\rho^{t^B,a,w^B}((t^C)^-,\tilde v)=w^A$, $\forall \omega\in\Omega$.

Since $v$ is $\mcF_t^\eta$--adapted we can write $v_t:=(\tau_1,\ldots,\tau_{N_t};\beta_1,\ldots,\beta_{N_t};\vartheta_{-N_T(a)+1},\ldots,\vartheta_{N_t})=\mu(t,W^\eta_{(t^A\vee\cdot)\wedge t}-W^\eta_{t^A})$, where $N_t:=\max\{j\in \{1,\ldots,N\}:\:\tau_j\leq t\}$ and $W^\eta_{(t^A\vee\cdot)\wedge t}$ is the trajectory of $W^\eta_s$ on the interval $[t^A,t]$. For all $t\in [t^B,T-(t^A-t^C)]$ we let
\begin{equation*}
\tilde v_t:=\left\{\begin{array}{l l}
\emptyset, & \text{for }t\in [t^B,t^C), \\
\mu(t+t^A-t^C,W^\eta_{(t^C\vee\cdot)\wedge t}-W^\eta_{t^C}),\quad &  \text{for }t\in [t^C,T-(t^A-t^C)].
\end{array}\right.
\end{equation*}
For $t\in (T-(t^A-t^C),T]$ we let $\tilde v_t:=(\tilde\tau_1,\ldots,\tilde\tau_{\tilde N_{T-(t^A-t^C)}};\tilde\beta_1,\ldots,\tilde\beta_{\tilde N_{T-(t^A-t^C)}}; \tilde\vartheta_{-N_T(a)+1},\ldots,\tilde\vartheta_{\tilde N_{T-(t^A-t^C)}})$, where $\tilde N_t:=\max\{j\in \{1,\ldots,\tilde N\}:\:\tilde \tau_j\leq t\}$ and the $\tilde \vartheta_j$, with $\tilde \vartheta_j > T-(t^A-t^C)$ are chosen so that \eqref{execBDS} holds for $\rho^{t^B,a,w^B}(t,\tilde v)$ (\eg we can let $\tilde \vartheta_j=\inf \{t>T-(t^A-t^C):\: \rho^{t^B,a,w^B}(t^-,\tilde v)=\varsigma(\alpha^{\eta,t^B,a}(t^-,\tilde v))$).

Now, $\rho^{t^B,a,w^B}(t,\tilde v)\in [w^B,\,w^A]$ for all $t\in [t^B,\,t^C)$ and \eqref{execBDS} holds for $\rho^{t^B,a,w^B}(t,\tilde v)$ for all $t\in [t^C,\,T-(t^A-t^C)]$ since $v$ is admissible, and for $t\in (T-(t^A-t^C),\,T]$ by the definition of the $\tilde \vartheta_j$ with $\tilde \vartheta_j>T-(t^A-t^C)$. As the measurability constraints also hold we have $\tilde v\in \mcV^{\eta,R}_{t^B}$.\\

For the switching costs we have,
\begin{equation}\label{ekv:swCOSTequiv2}
\E\left[\sum_{t^B\leq \tilde\tau_j \leq T}K_{\tilde\beta_j}\right]=\E\left[\sum_{t^C\leq \tilde\tau_j \leq T-(t^A-t^C)}K_{\tilde\beta_j}\right]=\E\left[\sum_{t^A\leq \tau_j \leq T}K_{\beta_j}\right],
\end{equation}
Hence, implementing the control $\tilde v$ we find
\begin{align*}
&J^\eta(t^B,y^B,a,w^B;\tilde v)-J^\eta(t^A,y^A,a,w^A;v)
\\
&=E\Bigg[\int_{t^B}^T f^\eta_2(X_s^{\eta,t^B,y^B},\alpha^{\eta,t^B,a}(s,\tilde v),\rho^{t^B,a,w^B}(s,\tilde v))ds-\int_{t^A}^T f^\eta_2(X_s^{\eta,t^A,y^A},\alpha^{\eta,t^A,a}(s,v),\rho^{t^A,a,w^A}(s,v))ds\Bigg]
\\
&\leq E\Bigg[\int_{t^C}^{T-(t^A-t^C)} f^\eta_2(X_s^{\eta,t^B,y^B},\alpha^{\eta,t^B,a}(s,\tilde v),\rho^{t^B,a,w^B}(s,\tilde v))ds
\\
&\quad -\int_{t^A}^T f^\eta_2(X_s^{\eta,t^A,y^A},\alpha^{\eta,t^A,a}(s,v),\rho^{t^A,a,w^A}(s,v))ds\Bigg]+C_1^{f,\eta}(1+C^{X^\eta}_1(1+|y^B|))|t^A-t^B|.
\end{align*}
Now, for $t\in[t^A,\,T]$ we let $\hat X^{\eta,t^B,y^B}_t=\eta_L^\top m(t-(t^A-t^C))+\hat Z_t^{\eta,t^B,y^B}$ where $(\hat Z_t^{\eta,t^B,y^B};\,t^A\leq t \leq T)$ is the strong solution to,
\begin{align*}
d\hat Z_{t}^{\eta,t^B,y^B}&=-\gamma \hat Z_t^{\eta,t^B,y^B} dt+\sigma^\eta dW^\eta_t,\\
\hat Z_{t^A}^{\eta,t^B,y^B}&=X^{\eta,t^B,y^B}_{t^C}-\eta_L^\top m(t^C).
\end{align*}
Then
\begin{align*}
\E\left[\sup\limits_{s\in [t^A,T]}|\hat X^{\eta,t^B,y^B}_s-X_s^{\eta,t^A,y^A}|\right] &= \E\left[\sup\limits_{s\in [t^A,T]}|\eta_L(m(s-(t^A-t^C))-m(s))+\hat Z^{\eta,t^B,y^B}_s-Z^{\eta,t^A,y^A}_s|\right]
\\
&\leq C^m|t^A-t^C|+C^{Z^\eta}_2\left(|y^A-y^B|+(1+|y^A|\vee |y^B|)|t^C-t^B|^{1/2}\right)
\\
&\leq C_5^{X^\eta}\left\{|y^A-y^B|+(1+|y^A|\vee |y^B|)(|w^A-w^B|^{1/2}+|t^A-t^B|^{1/2})\right\},
\end{align*}
where $C_5^{X^\eta}:=C^{m}\sqrt{T}+C^{Z^\eta}_2$. Hence, we have,
\begin{align*}
&E\Bigg[\int_{t^C}^{T-(t^A-t^C)} f^\eta_2(X_s^{\eta,t^B,y^B},\alpha^{\eta,t^B,a}(s,\tilde v),\rho^{t^B,a,w^B}(s,\tilde v))ds
\\
&\quad\quad-\int_{t^A}^T f^\eta_2(X_s^{\eta,t^A,y^A},\alpha^{\eta,t^A,a}(s,v),\rho^{t^A,a,w^A}(s,v))ds\Bigg]
\\
&=E\Bigg[\int_{t^A}^T \left\{f^\eta_2(\hat X_s^{\eta,t^B,y^B},\alpha^{\eta,t^A,a}(s,v),\rho^{t^A,a,w^A}(s,v)) - f^\eta_2(X_s^{\eta,t^A,y^A},\alpha^{\eta,t^A,a}(s,v),\rho^{t^A,a,w^A}(s,v))\right\}ds\Bigg]
\\
&\leq C_2^{f,\eta}C_5^{X^\eta}\left\{|y^A-y^B|+(1+|y^A|\vee |y^B|)(|w^A-w^B|^{1/2}+|t^A-t^B|^{1/2})\right\}T
\end{align*}
and the proposition follows with $C^{J,\eta}=C^{J,\eta}_2:=C_1^{f,\eta}(1+C^{X^\eta}_1) \sqrt{T} + C_2^{f,\eta}C_5^{X^\eta}T$ for $w^A\geq w^B$ and $d^R<0$.\\

\bigskip

Assume now that $w^A<w^B$. In this setting we cannot obtain a candidate for $\tilde v$ by performing a simple time-shift of $v$ as any such could lead to a violation of \eqref{execBDS}. Instead a time shift combined with a temporary fast-forward will be used to define $\tilde v$.

Let $w^C:=(2w^B-w^A) \wedge \frac{\varsigma (a)+w^B}{2}$ and define $t^C:=t^B+(w^C-w^B)/n_G$. Then, $t^C$ is the first time after $t^B$ that $\rho^{t^B,a,w^B}(\cdot,\tilde v)$ can reach $w^C$, for any $\tilde v\in\mcV_{t^B,a,w^B}^{\eta,R}$. We define the constant $c:=\frac{w^C-w^A}{t^C-t^B}=\frac{w^C-w^A}{w^C-w^B}n_G>1$.

If $t^C\geq T$, then $\emptyset \in \mcV_{t^B,a,w^B}^{\eta,R}$ and we have $J^\eta(t^B,y^B,a,w^B;\emptyset)\leq C^{f,\eta}_1 (1+C^{X^\eta}_1)\sqrt{\varsigma_{\max}}(1+|y^B|)|w^B-w^A|^{1/2}$, hence we are done. Therefore, we will from now on assume that $t^C< T$.

If, $t^A+c(t^C-t^B)\geq T$ we let $\tilde v:=(\tilde \vartheta_{-N_T(a)+1},\ldots,\tilde \vartheta_0)$, where $\tilde \vartheta_j:=\inf \{t>t^B:\: \rho^{t^B,a,w^B}(t^-,\tilde v)=\varsigma(\alpha^{\eta,t^B,a}(t^-,\tilde v))\}$. Then, $T-t^B<|t^A-t^B|+2(w^B-w^A)$ and we have $J^\eta(t^B,y^B,a,w^B;\tilde v)\leq C^{f,\eta}_1(1+ C^{X^\eta}_1)(2\sqrt{\varsigma_{\max}}\vee \sqrt{T})(1+|y^B|)(|w^B-w^A|^{1/2}+|t^A-t^B|^{1/2})$. To avoid trivialities we thus also assume that $t^A+c(t^C-t^B)< T$.

For the same reason we assume that $t^B+c(t^C-t^B) < T$.\\

As before we use the fact that the control $v$ can be represented by some function $\mu$ as $v_t=\mu(t,W^\eta_{(\cdot\wedge t)\vee t^A}-W^\eta_{t^A})$, for all $t\in[t^A,T]$. We define the control $\hat v^{\text{FF}}=(\hat \tau_1^{\text{FF}},\ldots,\hat\tau_{\hat N^{\text{FF}}}^{\text{FF}};\hat \beta_1^{\text{FF}},\ldots,\hat\beta_{\hat N^{\text{FF}}}^{\text{FF}};\hat \vartheta_{-N_T(a)+1}^{\text{FF}},\ldots,\hat\vartheta_{\hat N^{\text{FF}}}^{\text{FF}})$ and the family of controls $\hat v^s=(\hat \tau_1^s,\ldots,\hat\tau_{\hat N}^s;\hat \beta_1^s,\ldots,\hat\beta_{\hat N}^s;\hat \vartheta^s_{-N_T(a)+1},\ldots,\hat\vartheta^s_{\hat N})$ for all $s\in [t^B,T]$, as
\begin{equation}\label{hatvcDEF}
\hat v_t^{\text{FF}} :=  \mu\big(t^A+c(t-t^B),c\left(W^\eta_{\left(\frac{(t^B\vee\cdot)-t^B}{c}+t^B\right)\wedge t}-W^\eta_{t^B}\right)\big),\quad \text{if }t^B\leq t\leq t^C
\end{equation}
and
\begin{equation}\label{hatvDEF}
\hat v_t^{s} :=  \mu\big(t^A+(t-t^B),\hat W^{\eta,s}_{(t^B\vee \cdot\wedge t)}\big),\quad \text{if }t^B\leq t\leq T-(t^A-t^B)^+,
\end{equation}
where
\begin{align*}
\hat W^{\eta,s}_t:=\left\{\begin{array}{l l}
c(W^\eta_{(t-t^B)/c+t^B}-W^\eta_{t^B}), &\text{if } t\leq t^B+c(s - t^B), \\
c(W^\eta_{(s-t^B)/c+t^B}-W^\eta_{t^B})+W^\eta_{t-s+(s-t^B)/c+t^B}-W^\eta_{(s-t^B)/c+t^B},\quad &\text{if } t^B+c(s - t^B) < t \leq T. \\
\end{array}\right.
\end{align*}
Then, by an appropriate choice of the $\hat\vartheta_j$ with, according to \eqref{hatvDEF}, $\hat\vartheta_j>T-(t^A-t^B)^+$, we have $\hat v^s\in \mcV_{t^B,a,w^A}$, for all $s\in[t^B,T]$.

Assume now that $s=t^C$. From \eqref{hatvcDEF} and \eqref{hatvDEF} we find that $\hat \tau_j^{\text{FF}}-t^B = c(\hat \tau_j^{t^C}-t^B)$ and $\hat\beta_j^{\text{FF}}=\hat\beta_j^{t^C}$, for $j=1,\ldots,\hat N^{\text{FF}}$ and $\hat\vartheta^{\text{FF}}_j-t^B=c(\hat\vartheta^{t^C}_j-t^B)$ for all $j\in\{j'\in\{-N_T(a)+1,\ldots,\hat N^{\text{FF}}\}:\,\hat\vartheta^{\text{FF}}_{j'}\leq t^{\text{FF}}\}$. This means that $\hat v^{\text{FF}}_t$ can be seen as a fast-forward of the control $\hat v^{t^C}$ with rate $c>1$. Furthermore, $\hat v_t^{t^C}$ is $\mcF^\eta_{(t-t^B)/c+t^B}$--measurable, for all $t\in [t^B,\,t^B+c(t^C-t^B)]$ and $\hat \tau_j^{\text{FF}}$ is an $\mcF_t$--stopping time. Since $\alpha^{\eta,t^B,a} (t,\hat v^{\text{FF}})=\alpha^{\eta,t^B,a} (t^B+c(t-t^B),\hat v^{t^C})$, $\forall (t,\omega)\in [t^B,t^C]\times\Omega$, we have $\hat\beta_j^{\text{FF}}=\hat\beta_j^{t^C}\in\Lambda(\alpha^{\eta,t^B,a}(\hat\tau_j,\hat v^{t^C}))=\Lambda(\alpha^{\eta,t^B,a}(\hat\tau_j^{\text{FF}},\hat v^{\text{FF}}))$, $\forall\omega\in\Omega$, for $j=1,\ldots,\hat N^{\text{FF}}$.

Let
\begin{equation*}
\nu:=\inf\{t\geq t^B:\rho^{t^B,w^B}(t,\hat v^{\text{FF}})\leq \rho^{t^B,w^A}(c(t-t^B)+t^B,\hat v^{t^C})\}.
\end{equation*}
Now, since $N_T(\alpha^{\eta,t^B,a} (t,\hat v^{\text{FF}}))\geq 1$, for all $t\in[t^B,t^C]$, we have
\begin{align*}
\rho^{t^B,a,w^B}(t^C,\hat v^{\text{FF}})-\rho^{t^B,a,w^A}(t^B+c(t^C-t^B),\hat v^{t^C})&\leq w^B-w^A-(c-1)(t^C-t^B)\\
&=\left(\frac{1}{n_G}-1\right)(w^C-w^B)\leq 0.
\end{align*}
Hence, $\nu\leq t^C$, $\forall \omega\in\Omega$ and $\nu$ is an $\mcF^\eta_t$--stopping time. Furthermore, by continuity of $\rho^{t^B,a,w^B}(t,\hat v^{\text{FF}})-\rho^{t^B,a,w^A}(t^B+c(t-t^B),\hat v^{t^C})$ on $[t^B,\,t^C]$ we have $\rho^{t^B,a,w^B}(\nu,\hat v^{\text{FF}})=\rho^{t^B,a,w^A}(t^B+c(\nu-t^B),\hat v^{t^C})$, $\forall \omega\in\Omega$.

From this we conclude that
\begin{equation*}
\rho^{t^B,a,w^B}(t,\hat v^{\text{FF}})\geq \rho^{t^B,a,w^A}(t^B+c(t-t^B),\hat v^{t^C})>0,
\end{equation*}
for all $(t,\omega)\in [t^B,\nu]\times \Omega$. Furthermore, we have\footnote{Which implies that $N_T(\alpha^{\eta,t^B,a} (t,\hat v^{\text{FF}}))\geq 1$, for all $(t,\omega)\in[t^B,t^C]\times\Omega$.}
\begin{equation*}
\rho^{t^B,a,w^B}(t,\hat v^{\text{FF}})\leq \varsigma(\alpha^{\eta,t^B,a}(t,\hat v^{\text{FF}}))-(\varsigma(a)-w^C) < \varsigma(\alpha^{\eta,t^B,a}(t,\hat v^{\text{FF}})),
\end{equation*}
for all $(t,\omega)\in [t^B,t^C]\times \Omega$.


We now define the control $\tilde v$ as
\begin{equation}\label{ekv:tildvDEF3}
\tilde v_t := \left\{\begin{array}{l l}
\hat v_t^{\text{FF}} & \text{if }t^B\leq t\leq \nu \\
\hat v^\nu_{t+(c-1)(\nu-t^B)}\quad & \text{if }\nu\leq t\leq T-(c-1)(\nu-t^B)
\end{array}\right.
\end{equation}
and assume that $\tilde N:=\tilde N_{T-(c-1)(\nu-t^B)}$ (\ie no more switching requests are made in $(T-(c-1)(\nu-t^B),\,T])$, and pick $\tilde \vartheta_j=\inf \{s>T-(t^A-\nu+c(\nu-t^B))^+:\,\rho^{t^B,a,w^B}(s^-,\tilde v)=\varsigma(\alpha^{\eta,t^B,a}(s^-,\tilde v))\}$ for all $j\in\{j'\in\{-N_T(a)+1,\ldots,\tilde N\}:\,\tilde\vartheta_{j'}>T-(t^A-\nu+c(\nu-t^B))^+\text{ according to \eqref{ekv:tildvDEF3}}\}$. Then $\tilde v\in \mcV^{\eta,R}_{t^B,a,w^B}$.\\


Now, let $(\hat X_t^{\eta,\nu};\:t^B\leq t\leq T-(t^A-t^B)^+)$ be defined as $\hat X_t^{\eta,\nu}=m(t+t^A-t^B)+\hat  Z_t^{\eta,\nu}$, where $(\hat  Z_t^{\eta,\nu};\:t^B\leq t\leq T-|t^A-t^B|)$ is the strong solution to
\begin{align*}
d\hat  Z_t^{\eta,\nu}&=-\gamma \hat  Z_t^{\eta,\nu}dt+\sigma d\hat W_t^{\eta,\nu},\quad t\in [t^B,T-|t^A-t^B|]\\
\hat Z_{t^B}^{\eta,\nu}&=\eta_L^\top (x^A-m(t^A)).
\end{align*}
We have
\begin{align*}
&\E\left[ \sup\limits_{s\in [t^B,\, T-|t^A-t^B|]} | X^{\eta,t^B,y^B}_{s}-\hat X_s^{\eta,\nu}|\right] \leq C^m|t^A-t^B|+\E\left[ \sup\limits_{s\in [t^B,\, T-|t^A-t^B|]} | Z^{\eta,t^B,y^B}_{s}-\hat Z_s^{\eta,\nu}|\right]
\\
&\hspace{100pt}\leq C^m |t^A-t^B|+C^{Z^\eta}_3\left(|y^A-y^B|+(1+|y^A|\vee |y^B|)(c(t^C-t^B))^{1/2}\right)
\\
&\hspace{100pt}\leq  C^{X^\eta}_6\left(|y^A-y^B|+(1+|y^A|\vee |y^B|)(|w^A-w^B|^{1/2}+|t^A-t^B|^{1/2})\right),
\end{align*}
with $C^{X^\eta}_6:=C^m\sqrt{T}\vee \sqrt{2}C_3^{Z^\eta}$, and
\begin{equation}\label{ekv:swCOSTequiv2}
\E\left[\sum_{t^B\leq \hat\tau_j^\nu \leq T}K_{\hat\beta_j^\nu}\right]=\E\left[\sum_{t^A\leq \tau_j \leq T-(t^A-t^B)^+}K_{\beta_j}\right]\leq \E\left[\sum_{t^A\leq \tau_j \leq T}K_{\beta_j}\right].
\end{equation}
Hence,
\begin{align*}
&J^\eta(t^B,y^B,w^A;\hat v^\nu)-J^\eta(t^A,y^A,w^A;v)\\
&=\E\Bigg[\int_{t^B}^T f^\eta_2(X_s^{\eta,t^B,y^B},\alpha^{\eta,t^B,a}(s,\hat v^\nu),\rho^{t^B,a,w^A}(s,\hat v^\nu))ds-\int_{t^A}^T f^\eta_2(X_s^{\eta,t^A,y^A},\alpha^{\eta,t^A,a}(s,v),\rho^{t^A,a,w^A}(s,v))ds\Bigg]
\\
&=\E\Bigg[\int_{t^B}^{T-|t^A-t^B|}\left\{f^\eta_2(X_s^{\eta,t^B,y^B},\alpha^{\eta,t^B,a}(s,\hat v^\nu),\rho^{t^B,a,w^A}(s,\hat v^\nu))-f^\eta_2(\hat X_s^{\eta,\nu},\alpha^{\eta,t^B,a}(s,\hat v^\nu),\rho^{t^B,a,w^A}(s,\hat v^\nu))\right\}ds
\\
&\qquad+\int_{T-(t^A-t^B)^+}^T f^\eta_2(X_s^{\eta,t^B,y^B},\alpha^{\eta,t^B,a}(s,\hat v^\nu),\rho^{t^B,a,w^A}(s,\hat v^\nu))ds
\\
&\qquad- \int_{T-(t^B-t^A)^+}^T f^\eta_2(X_s^{\eta,t^A,y^A},\alpha^{\eta,t^A,a}(s,v),\rho^{t^A,a,w^A}(s,v))ds\Bigg]
\\
&\leq C^{f,\eta}_2 C^{X^\eta}_6\left(|y^A-y^B|+(1+|y^A|\vee |y^B|)(|w^A-w^B|^{1/2}+|t^A-t^B|^{1/2})\right)T\\
&\quad+C^{f,\eta}_1(1+C_1^{X^\eta})(1+|y^A|\vee |y^B|)|t^A-t^B|
\\
&\leq C^{J,\eta}_{3,1} \left(|y^A-y^B|+(1+|y^A|\vee |y^B|)(|w^A-w^B|^{1/2}+|t^A-t^B|^{1/2})\right),
\end{align*}
with $C^{J,\eta}_{3,1}:=C^{f,\eta}_2 C^{X^\eta}_6 T + C^{f,\eta}_1(1+C_1^{X^\eta})\sqrt{T}$.\\

Furthermore, we have
\begin{align*}
&J^\eta(t^B,y^B,w^B;\tilde v)-J^\eta(t^B,y^B,w^A;\hat v^\nu)\\
&=\E\Bigg[\int_{t^B}^T \left\{f^\eta_2(X_s^{\eta,t^B,y^B},\alpha^{\eta,t^B,a}(s,\tilde v),\rho^{t^B,a,w^B}(s,\tilde v))ds- f^\eta_2(X_s^{\eta,t^B,y^B},\alpha^{\eta,t^B,a}(s,\hat v^\nu),\rho^{t^B,a,w^A}(s,\hat v^\nu))\right\}ds\Bigg]
\\
&=\E\Bigg[\int_{t^B}^{\nu}f^\eta_2(X_s^{\eta,t^B,y^B},\alpha^{\eta,t^B,a}(s,\tilde v),\rho^{t^B,a,w^B}(s,\tilde v))ds
\\
&\qquad - \int_{t^B}^{t^B+c(\nu-t^B)}f^\eta_2(X_s^{\eta,t^B,y^B},\alpha^{\eta,t^B,a}(s,\hat v^\nu),\rho^{t^B,a,w^A}(s,\hat v^\nu))ds
\\
&\qquad+\int_{\nu}^{T-(c-1)(\nu-t^B)}\Big\{f^\eta_2(X_s^{\eta,t^B,y^B},\alpha^{\eta,t^B,a}(s,\tilde v),\rho^{t^B,a,w^B}(s,\tilde v))
\\
&\qquad\quad-f^\eta_2(X_{s+(c-1)(\nu-t^B)}^{\eta,t^B,y^B},\alpha^{\eta,t^B,a}(s,\tilde v),\rho^{t^B,a,w^B}(s,\tilde v))\Big\}ds
\\
&\qquad + \int_{T-(c-1)(\nu-t^B)}^T f^\eta_2(X_s^{\eta,t^B,y^B},\alpha^{\eta,t^B,a}(s,\tilde v),\rho^{t^B,a,w^B}(s,\tilde v))ds\Bigg]
\\
&\leq 3C^{f,\eta}_1(1+C_1^{X^\eta})(1+|y^B|)c(t^C-t^B)+ C^{X^\eta}_3(1+|y^B|)\{(c-1)(t^C-t^B)\}^{1/2}T,
\\
&\leq C^{J,\eta}_{3,2}(1+|y^B|)|w^A-w^B|^{1/2},
\end{align*}
where $C^{J,\eta}_{3,2}:=3C^{f,\eta}_1(1+C_1^{X^\eta})\sqrt{2\varsigma_{\max}}+\sqrt{2}C^{X^\eta}_3 T$.

Hence,
\begin{align*}
J^\eta(t^B,y^B,w^B;\tilde v)-J^\eta(t^A,y^A,w^A;v)\leq (C^{J,\eta}_{3,1}+C^{J,\eta}_{3,2})\Big(|y^A-y^B|
\\+(1+|y^A|\vee |y^B|)(|w^A-w^B|^{1/2}+|t^A-t^B|^{1/2})\Big).
\end{align*}
\bigskip

We conclude that \eqref{ekv:JETAprop} holds with $C^{J,\eta}:=C^{J,\eta}_1\vee C^{J,\eta}_2 \vee (C^{J,\eta}_{3,1}+C^{J,\eta}_{3,2})$.\qed\\

\bigskip

This leads us to the equivalent of Proposition~\ref{VFprop} for the value function of Problem LB.
\begin{prop}\label{VFpropETA}
There are constants $C^{V,\eta}_1>0$ and $C^{V,\eta}_2>0$ such that, for all $a\in\Gamma$, the following holds
\begin{align}\label{ekv:VFpropETA1}
|V^\eta(t,y,a,w)|&\leq C^{V,\eta}_1(1+|y|),\quad \forall (t,y,w)\in [0,T]\times \R\times [0,\varsigma(a)) \\
|V^\eta(t^A,y^A,a,w^A)-V^\eta(t^B,y^B,a,w^B)| &\leq C^{V,\eta}_2\big\{|y^A-y^B|+(1+|y^A|+|y^B|)(|w^A-w^B|^{1/2}\nonumber\\
&\quad\quad\quad\quad+|t^A-t^B|^{1/2})\big\}.\label{VFpropETA2}
\end{align}
for all $(t^A,y^A,w^A),\,(t^B,y^B,w^B)\in [0,T]\times \R\times [0,\varsigma(a))$.
\end{prop}

\noindent\emph{Proof.} Applying \eqref{ekv:fETAprop1} we find that \eqref{ekv:VFpropETA1} holds with $C^{V,\eta}_1:=C^{f,\eta}_1(1+C_1^{X^\eta}) T$. The second inequality is an immediate consequence of Proposition~\ref{JETApropAB}.\qed\\

\bigskip
\bigskip

\begin{thm}(Principle of optimality)
For each $t\in[0,T]$, $y\in\R$, $a\in\Gamma$ and $w\in [0,\varsigma(a))$, we have
\begin{align*}
V^\eta(t,y,a,w)&=\inf\limits_{v\in \mcV^{\eta,R}_{t,a,w}}\E\Bigg[\int_t^\nu  f^\eta_2(X^{\eta,t,y}_s,\alpha^{\eta,t,a}(s,v),\rho^{t,a,w}(s,v))ds + \sum\limits_{\tau_i \leq \nu}K_{\beta_i}
\\
&\hspace{120pt} + V^\eta(\nu,X_\nu^{\eta,t,y},\alpha^{\eta,t,a}(\nu,v),\rho^{t,a,w}(\nu,v)) \Bigg],
\end{align*}
for each $\nu\in \mcT^\eta_t$, where $\mcT_t^\eta$ is the set of $\mcF_t^\eta$--stopping times $\tau$ with $\tau\geq t$.
\end{thm}

\noindent\emph{Proof.} Follows exactly the same line as the proof of Theorem~\ref{Popt1} and is therefore omitted.\qed\\

\begin{cor}
For each $t\in[0,T]$, $y\in\R$, $a\in\Gamma$ and $w\in [0,\varsigma(a))$, we have
\begin{align*}
V^\eta(t,y,a,w)=&\inf\limits_{\tau\in \mcT^\eta_{t,a,w}}\E\Bigg[\int_t^\tau  f^\eta_2(X^{\eta,t,y}_s,a,w+N_T(a)(s-t))ds +1_{\left\{t+\frac{\varsigma(a)-w}{N_T(a)}\right\}}(\tau)V^\eta(\tau,X_\tau^{\eta,t,y},a^{t},0)\\
&+ 1_{\left[t,t+\frac{\varsigma(a)-w}{N_T(a)}\wedge T\right)}(\tau)\Big\{\Big(\min_{k\in \mathcal{K}^T_{a,w,\tau-t}}V^\eta(\tau,X_\tau^{\eta,t,y},a^{t,k},w+N_T(a)(\tau-t)-\delta_{k,a_k})\Big)
\\
&\quad\bigwedge\min_{\beta\in \Lambda(a)} \Big(K_{\beta_i} + V^\eta(\tau,X_\nu^{\eta,t,y},a+\beta,w+N_T(a)(\tau-t))\Big)\Big\} \Bigg],
\end{align*}
where $\mcT_{t,a,w}^\eta$ is the set of $\mcF_t^\eta$--stopping times $\tau$ with $t\leq \tau \leq t+\frac{\varsigma(a)-w}{N_T(a)}$, we use $\mathcal{K}^T_{a,w,s}$ to denote the set of all $k\in \mathcal{K}^T(a)$ such that $w+N_T(a)s\geq \delta_{k,a_k}$ and $a^{t,k}\in\Gamma$ is given by $(a^{t,k})_j=a_j$ for $j\neq k$ and $(a^{t,k})_k=a^t_k$.
\end{cor}

\subsection{An upper bound}
To obtain an upper bound on the value function we approximate the set of feasible operating points from the inside by a continuum of generalized ellipsoids. The generalized ellipsoids will be based on a number of signed distance functions $d_i:[0,T]\times D_G\to \R$, representing a weighted distance from the forecasted demand at time $t\in[0,T]$, $m(t)$, to the set of non-feasible operating points. Combined with the approach of solving impulse control problems with delays outlined in~\cite{OksenImpulse} this will give us an upper bound that can be computed using conventional numerical methods.

\subsubsection{Running cost approximation}
We start with the upper bound approximation of the running cost.
\begin{defn}
For each $t\in [0,T]$ and each $p_G\in D_G$, we let
\begin{equation*}
\tilde d_i(t,p_G):=\min\limits_{p_D\in\R^{N_s}}\{\|\sigma^{-1}(p_D-m(t))\|:\,(p_G,p_D)\in\partial G_i\}.
\end{equation*}
We now define $d_i(t,p_G)=\tilde d_i(t,p_G)$ if $m(t)\in G_i(p_G)$ and $d_i(t,p_G)=-\tilde d_i(t,p_G)$ when $m(t)\notin G_i(p_G)$.
\end{defn}

\begin{rem} For each $t$ in a finite subset of $[0,T]$ we may wish to decide $d_i(t,\cdot)$ on a finite set $\hat D_G:=\{p^1_G\ldots,p^M_G\}\subset D_G$ of sample points. Due to the convexity of $G_i$ the approximation $\hat d_i(t,\cdot):D_G\to\R$ given, for each $p_G\in D_G$, by
\begin{equation}\label{ekv:dBext}
\hat d_i(t,p_G):=\max\{d\in\R: (p_G,d)\in \text{Conv}(\{(p^l_G,d_i(t,p^l_G)):\,l=1,\ldots,M \})\},
\end{equation}
satisfies $\hat d_i(t,p_G)\leq d_i(t,p_G)$, $\forall p_G\in D_G$.
\end{rem}

\begin{defn}
For each $p_G\in D_G$ and each $t\in[0,T]$, we define
\begin{equation}
B^{\sigma}(m(t),d_i(t,p_G)):=\{p_D\in\R^{n_L} :\,\|\sigma^{-1}(p_D-m(t))\|\leq d_i(t,p_G)\},
\end{equation}
for $i=0,\ldots,n_c$.
\end{defn}

From the definition it is clear that $B^{\sigma}(m(t),d_i(t,p_G))\subset G_i(p_G)$.

\begin{lem}\label{lem:CB}
For $i=0,\ldots,n_c$ we let $C_{i}^B:[0,T]\times D_G \times \R_+\to \R_+$ be defined as
\begin{equation}\label{ekv:CBdef}
C_{i}^B(t,p_G,y):=\max(0,(y-d_i(t,p_G)))\max\limits_{\|\xi\|=1}c^\top(\sigma \xi)^+.
\end{equation}
Then,
\begin{equation}\label{ekv:CBdom}
c^\top z_i^*(p_G,p_D) \leq C_{i}^B(t,p_G,\|\sigma^{-1}(p_D-m(t))\|),\quad i=0,\ldots,n_c,
\end{equation}
for all $(p_G,p_D)\in D_G\times \R^{n_L}_+$.
\end{lem}

\noindent\emph{Proof.} Assume first that $d_i(t,p_G)\geq 0$. If $\|\sigma^{-1}(p_D-m(t))\|\leq d_i(t,p_G)$, then $p_D\in B^{\sigma}(m(t),d_i(t,p_G))\subset G_i(p_G)$. Hence, \eqref{ekv:CBdom} trivially holds with both sides equal to zero in this case.

Assume instead that $\|\sigma^{-1}(p_D-m(t))\|> d_i(t,p_G)$. Let $z\in \R^{n_L}$ be such that $\|\sigma^{-1}(p_D-z-m(t))\| \leq d_i(t,p_G)$, then $p_D-z\in B^{\sigma}(m(t),d_i(t,p_G))\subset G_i(p_G)$. Hence, by Assumption~\ref{AssG}.\ref{assDOMIN}, $p_D-z^+\in G_i(p_G)$. We have, with $y=\|\sigma^{-1}(p_D-m(t))\|$,
\begin{align*}
c^\top z_i^*(p_G,p_D)&\leq\max\limits_{x\in\R^{n_L}_+}\left\{\min_{z\in \R^{n_L}_+}\{c^\top z:\: (p_G,x-z)\in G_i\}:\|\sigma^{-1}(x-m(t))\|=y\right\}
\\
&\leq \max\limits_{x\in\R^{n_L}}\left\{\min_{z\in \R^{n_L}}\{c^\top z^+:\:\: \|\sigma^{-1}(x-z-m(t))\|\leq d_i(t,p_G)\}:\|\sigma^{-1}(x-m(t))\|=y\right\}
\\
&\leq \max\limits_{\tilde x\in\R^{n_L}}\left\{\min_{\hat x\in \R^{n_L}}\{c^\top (\sigma(\tilde x-\hat x))^+:\: \|\hat x\|\leq d_i(t,p_G)\}:\:\|\tilde x\|=y\right\}
\\
&=(y-d_i(t,p_G))\max_{\|\xi\|=1}c^\top (\sigma\xi)^+
\\
&=C_{i}^B(t,p_G,\|\sigma^{-1}(p_D-m(t))\|),
\end{align*}
for $i=0,\ldots,n_c$.\\

If $d_i(t,p_G) < 0$, then there is a $z_0\in\R^{n_L}$, with $\|\sigma^{-1}z_0\|=-d_i(t,p_G)$ such that $m(t)-z_0\in B^{\sigma}(m(t),d_i(t,p_G))\subset G_i(p_G)$. Now, for each $x\in\R^{n_L}$ we have that $x-(x-m(t))-z_0 \in B^{\sigma}(m(t),d_i(t,p_G))\subset G_i(p_G)$. Hence, $x-(x-m(t)+z_0)^+\in G_i(p_G)$ and we have, with $y=\|\sigma^{-1}(p_D-m(t))\|$,
\begin{align*}
c^\top z_i^*(p_G,p_D)&\leq\max\limits_{x\in\R^{n_L}_+}\left\{\min_{z\in \R^{n_L}_+}\{c^\top z:\: (p_G,x-z)\in G_i\}:\|\sigma^{-1}(x-m(t))\|=y\right\}
\\
&\leq \max\limits_{x\in\R^{n_L}}\left\{\max_{z_0\in \R^{n_L}}\{c^\top (x-m(t)+z_0)^+:\: \|\sigma^{-1}z_0\|=-d_i(t,p_G)\}:\:\|\sigma^{-1}(x-m(t))\|=y\right\}
\\
&= \max\limits_{\tilde x\in\R^{n_L}}\left\{\max_{\tilde z_0\in \R^{n_L}}\{c^\top (\sigma(\tilde x+\tilde z_0))^+:\: \|\tilde z_0\|=-d_i(t,p_G)\}:\:\|\tilde x\|=y\right\}
\\
&=(y-d_i(t,p_G))\max_{\|\xi\|=1}c^\top (\sigma\xi)^+
\\
&=C_{i}^B(t,p_G,\|\sigma^{-1}(p_D-m(t))\|),
\end{align*}
for $i=0,\ldots,n_c$.\qed\\

\begin{defn}
For each $p_G\in D_G$ and each $y\in\R_+$ we define
\begin{equation}\label{ekv:fB1DEF}
f^B_1(t,p_G,y):=C_G(p_G) + \sum_{i=0}^{n_c}q_i C_{i}^B(t,p_G,y).
\end{equation}
\end{defn}
Note that $f^B_1$ is convex in $p_G$ (whenever $C_G$ is convex) and convex and monotonically increasing in $y$. We will now use $f^B_1$ to define the running cost for the upper bound problem.
\begin{defn}
Given $a\in\Gamma$ and $w\in [0,\psi(a))$, we let
\begin{equation}\label{ekv:fB2DEF}
f^B_2(t,y,a,w):=\max_{p_G\in D_G(a,w)}f^B_1(t,p_G,y),\quad \forall (t,y)\in [0,T]\times \R_+,
\end{equation}
where $D_G(a,w):=\overline{\{\zeta(a,r)\in D_G:\: \exists r\in H(a),\:\psi(r,a)=w\}}$.
\end{defn}
The next proposition gives us some important properties of $f^B_2$.
\begin{prop}
For all $a\in\Gamma$, and each $w\in [0,\psi(a))$ we have
\begin{align}\label{ekv:fBpropDOM}
f^B_2(t,\|\sigma^{-1}(x-m(t))\|,a,w)&\geq f(x,a',r),
\end{align}
for all $(t,x,a',r)\in [0,T]\times \R^{n_L}\times \mathop{\cup}\limits_{a\in \Gamma} \{\{a\} \times H(a)\}$ such that $\zeta(a',r)\in D_G(a,w)$. Furthermore, there are constants $C_1^{f,B}>0$ and $C_2^{f,B}>0$ such that, for all $a\in\Gamma$ and all $w\in [0,\psi(a))$,
\begin{align}\label{ekv:fBprop1}
|f^B_2(t,y,a,w)| &\leq C_1^{f,B}(1+|y|),\quad \forall (t,y)\in [0,T]\times\R_+
\end{align}
and
\begin{align}\label{ekv:fBprop2}
|f^B_2(t^A,y^A,a,w)-f^B_2(t^B,y^B,a,w)| &\leq C_2^{f,B}(|t^A-t^B|+|y^A-y^B|),
\end{align}
for all $(t^A,y^A),(t^B,y^B)\in[0,T]\times \R_+$.
\end{prop}

\noindent\emph{Proof.} Combining Lemma~\ref{lem:CB} with \eqref{ekv:fB1DEF} and \eqref{ekv:fB2DEF} we find that, for all $(a',r)\in \cup_{a\in \Gamma}\{\{a\}\times H(a)\}$ such that $\zeta(a',r)\in D_G(a,w)$,
\begin{align*}
f^B_2(t,\|\sigma^{-1}(x-m(t))\|,a,w)&\geq C_G(\zeta(a',r))+\sum_{i=0}^{n_c}q_i C^B_i(t,\zeta(a',r),\|\sigma^{-1}(x-m(t)))
\\
&\geq f(x,a',r).
\end{align*}
For the second inequality, we note that
\begin{align*}
|f^B_2(t,y,a,w)|\leq \max_{t\in[0,T]}\max_{p_G\in D_G}\left\{C_G(p_G)+ \max_{\|\xi\|=1}c^\top (\sigma\xi)^+\sum_{i=0}^{n_c}q_i(|d_i(t,p_G)|+|y|)\right\}
\end{align*}
from which \ref{ekv:fBprop1} follows with $C_1^{f,B}:=\max\limits_{t\in[0,T]}\max\limits_{p_G\in D_G}\left\{C_G(p_G) + \max\limits_{\|\xi\|=1}c^\top (\sigma\xi)^+\sum_{i=0}^{n_c}q_i(|d_i(t,p_G)| \vee  1)\right\} $.\\

To prove the last inequality we let
\begin{equation*}
p_G^A:=\mathop{\arg\min}\limits_{p_G\in D_G(a,w)} f^B_1(t^A,p_G,y^A)
\end{equation*}
and note that there is a $C^d>0$ such that, for all $t^A,t^B\in[0,T]$,
\begin{equation*}
|d_i(t^A,p_G^A)-d_i(t^B,p_G^A)|\leq C^d|t^A-t^B|,\quad \text{for }i=0,\ldots,n_c.
\end{equation*}
Hence,
\begin{align*}
f^B_2(t^A,y^A,a,w)-f^B_2(t^B,y^B,a,w) &\leq f^B_1(t^A,p_G^A,y^A)-f^B_1(t^B,p_G^A,y^B)
\\
&=\max_{\|\xi\|=1}c^\top (\sigma\xi)^+\sum_{i=0}^{n_c}q_i \left\{(y^A-d_i(t^A,p_G^A))^+-(y^B-d_i(t^B,p_G^A))^+\right\}
\\
&\leq \max_{\|\xi\|=1}c^\top (\sigma\xi)^+\sum_{i=0}^{n_c}q_i \left((y^A-y^B-(d_i(t^A,p_G^A)-d_i(t^B,p_G^A))\right)^+
\\
&\leq \max_{\|\xi\|=1}c^\top (\sigma\xi)^+\sum_{i=0}^{n_c}q_i \left(|y^A-y^B|+C^d|t^A-t^B|\right).
\end{align*}
By interchanging $(t^A,y^A)$ and $(t^B,y^B)$ we find that \eqref{ekv:fBprop2} holds with $C_2^{f,B}:=\max\limits_{\|\xi\|=1}c^\top (\sigma\xi)^+\sum_{i=0}^{n_c}q_i(1\vee C^d)$.\qed\\

\subsubsection{Admissible controls}
The running-cost, $f^B_2$ only depends on the demand through the process $(X^B_t;0\leq t\leq T)$ defined as $X^B_t:=\|\sigma^{-1}Z_t\|$. It is well known that $X^B$ is a Markov process. To see this we let $\tilde Z=\sigma^{-1}Z$ and get
\begin{align*}
d\tilde Z_t=-\gamma \tilde Z_t dt +dW_t.
\end{align*}
If we consider the square process
\begin{align}\nonumber
d(\tilde Z_t^\top \tilde Z_t)&=(d\tilde Z^\top_t)\tilde Z_t+\tilde Z_t^\top d\tilde Z_t +d\tilde Z_t^\top d\tilde Z_t\\
\nonumber&=-2\gamma \tilde Z_t^\top \tilde Z_tdt+2\tilde Z^\top dW_t+n_L dt\\
&=\{n_L-2\gamma \tilde Z_t^\top \tilde Z_t\}dt+2\tilde Z^\top dW_t\label{kvadSDE}
\end{align}

We let $W^1$ be a one-dimensional Wiener process and let $R_t$ be the Cox-Ingersoll-Ross process~\cite{CIRart} that solves
\begin{align}
dR_t &=(n_L-2\gamma R_t)ds+2\sqrt{R_t} dW_t^1,\quad t\in[0,T].
\end{align}
Then $R_t | R_s=\|\sigma^{-1}Z_s\|^2$, $0\leq s\leq t\leq T$ is a weak solution to \eqref{kvadSDE} on $[s,T]$.\\

\begin{defn}
Let $\mcV^B$ be the set of controls $v=(\tau_1,\ldots,\tau_N;\beta_1,\ldots,\beta_N)\in\mcV$ that are adapted to the filtration $\mcF_t^B$ generated by $X^B_t$, with $\tau_1\geq \psi(\alpha_0,\xi_0)$ and, for $j=1,\ldots,N-1$, either $\tau_j=\tau_{j+1}$ or $\tau_j+\max\limits_{i:\,\tau_i=\tau_j} \delta_{\beta_i}\leq \tau_{j+1}$.

Furthermore, for $t\in[0,T]$, we let $\mcV_{t,a,w}^B$ be the controls $v\in\mcV^B$ with $\tau_1\geq t+w$.
\end{defn}
We thus limit the controls in a way that no switching requests can be made whilst in a transition mode, but several switching requests can be made simultaneously. This will allow us to use the method proposed to solve impulse control problems with delayed reaction in~\cite{OksenImpulse}.

\begin{defn}
For each $t\in[0,T]$, $a\in \Gamma$, and $w\in [0,\psi(a))$ we let $\alpha^{B,t,a,w}(s,v)$ be defined as
\begin{align}\nonumber
\alpha^{B,t,a,w}_k(s,v):=a_k1_{[t,t+w)}(s)+a^t_k1_{[t+w,\tau_1]}(s)+\sum_{j=1}^{N^k} \Big\{(\beta^{k,f}_j, \beta^{k,t}_j) 1_{[\tau_j^k,\tau_j^k+\max\limits_{j':\:\tau_{j'}=\tau^k_j} \delta_{\beta_{j'}})}(s)
\\
+\beta^{k,t}_j1_{[\tau_j^k+\max\limits_{j':\:\tau_{j'}=\tau^k_j} \delta_{\beta_{j'}},\tau^k_{j+1})}(s)\Big\}\label{alphaBdef}
\end{align}
and
\begin{equation}\label{alphaBdef}
\varrho^{t,a,w}(s,v):=(w-(s-t))^+ + \sum_{j=1}^N 1_{[\tau_j,\tau_{j+1})}(s)\left(\max\limits_{j':\:\tau_{j'}=\tau_j} \delta_{\beta_{j'}}-(s-\tau_j)\right)^+
\end{equation}
$\forall (s,v)\in [t,T]\times \mcV^B_{t,a,w}$.

\subsubsection{The upper bound problem}
We define the expected operating cost-to-go associated with the control $v\in \mcV^{B}_{t,a,w}$ as
\begin{equation*}
J^B(t,y,a,w;v):=\E\left[\int_t^T f^B_2(s,X^{B,t,y}_s,\alpha^{B,t,a,w}(s,v),\varrho^{t,a,w}(s,v)))ds+ \sum\limits_{t<\tau_i\leq T}K_{\beta_i}\right].
\end{equation*}
\end{defn}
To find the upper bound we solve the following problem.\\

\noindent\textbf{Problem UB.} For all $t\in [0,T]$, all $a\in\Gamma$ and all $y\in \R$, find
\begin{equation*}
V^{\rm{UB}}(t,y,a,w)=\inf\limits_{v\in\mcV^{B}_{t,a,w}} J^B(t,y,a,w;v).
\end{equation*}
for all $w\in [0,\psi(a))$. \qed\\

First we note that if $a\notin \Theta$ we can write
\begin{align*}
J^B(t,y,a,w;v)=\E\Bigg[\int_t^{(t+w)\wedge T} f^B_2(s,X^{B,t,y}_s,a,w-(s-t))ds+ J^B(t+w,X^{B,t,y}_{t+w},a^t,0;v)\Bigg]
\end{align*}
and otherwise we have
\begin{align*}
J^B(t,y,a,w;v)&=\E\Bigg[\int_t^{\tau_1\wedge T} f^B_2(s,X^{B,t,y}_s,a,0))ds
\\
& \quad + 1_{[t,T)}(\tau_1)\left\{\sum_{j:\:\tau_j=\tau_1}K_{\beta_j} + J^B(\tau_1,X^{B,t,y}_{\tau_1}, a^t+\sum_{j:\:\tau_j=\tau_1} \beta_j,\max_{j:\:\tau_j=\tau_1} \delta_{\beta_j};v)\right\}\Bigg].
\end{align*}


Problem UB can does be rewritten as a non-delay optimal switching problem, with a finite number of possible switches, and we have
\begin{align}\nonumber
V^{\rm{UB}}(t,y,a,w)&=\inf\limits_{\tau\in \mcT^B_{t+w}}\E\Bigg[\int_t^{\tau\wedge T} f^B_2(s,X^{B,t,y}_s,a1_{[t,t+w)}(s)+a^t1_{[t+w,T)}(s),(w-(s-t))^+)ds
\\
&\hspace{50pt}+1_{[t,T)}(\tau) \min_{\beta\in \Lambda^{\mcP}(a^t)}\left\{K_\beta+V^{\rm{UB}}(\tau,X^{B,t,y}_{\tau},a^t+\beta,\psi(a+\beta);v)\right\}\Bigg],\label{ekv:DPPB}
\end{align}
where, for $s\in [0,T]$, $\mcT^B_s$ is the set of $\mcF^B_t$--stopping times $\tau$, with $\tau\geq s$. We use the notation $\Lambda^{\mcP}(a^t)$ to represent the set of all vectors $\beta=(\beta_1,\ldots,\beta_m)$, with $m\in \{1,\ldots,n_G\}$, where $\beta_j\in \Lambda(a^t)$ for $j=1,\ldots,m$ and $\beta_j^{\Node}\neq \beta_k^{\Node}$ for $k\neq j$, and let $K_\beta:=\sum_{j=1}^m K_{\beta_j}$ and $a+\beta:=a+\sum_{\beta_j\in\beta}\beta_j$.


It only remains to show that the value function to Problem UB is, in fact, an upper bound to the value function of Problem 1.
\begin{thm}
For each $t\in [0,T]$, $x\in\R^{n_L}$ and $a\in\Gamma$ we have
\begin{equation}\label{ekv:isUB}
V^{\rm{UB}}(t,\|\sigma^{-1}(x-m(t))\|,a,\psi(a,r))\geq V(t,x,a,r),
\end{equation}
$\forall r\in H(a)$.
\end{thm}

\noindent\emph{Proof.} For $v\in\mcV^B_{t,a,\psi(a,r)}$ we have $v\in \mcV_{t,a,r}$ and
\begin{equation*}
\zeta(\alpha^{t,a,r}(s,v),\xi^{t,a,r}(s,v))\in D_G(\alpha^{B,t,a,\psi(a,r)}(s,v),\varrho^{t,a,\psi(a,r)}(s,v)),\quad\forall (s,\omega)\in [t,T]\times \Omega.
\end{equation*}
Now, using \eqref{ekv:fBpropDOM} we find that
\begin{align*}
J(t,x,a,r;v)&=\E\left[\int_t^T  f(X_s^{t,x},\alpha^{t,a,r}(s,v),\xi^{t,a,r}(s,v))ds + \sum\limits_{\tau_i\leq T}K_{\beta_i}\right]
\\
&\leq \E\Bigg[\int_{t}^T  f^B(s,X_s^{B,t,\|\sigma^{-1}(x-m(t))\|},\alpha^{B,t,a,\psi(a,r)}(s,v),\varrho^{t,a,\psi(a,r)}(s,v))ds + \sum\limits_{t<\tau_i\leq T}K_{\beta_i}\Bigg]
\\
&=J^B(t,\|\sigma^{-1}(x-m(t))\|,a,\psi(a,r);v).
\end{align*}
Taking the infimum over the set $\mcV_{t,a,r}$ on the left hand side of the inequality and over $\mcV^B_{t,a,\psi(a,r)}\subset \mcV_{t,a,r}$ on the right hand side \eqref{ekv:isUB} follows.\qed\\


\section{Numerical solution algorithm\label{Sec:Num}}
The general approach to solving high-dimensional optimal switching problems is based on using Monte Carlo simulation combined with regression to approximate the conditional expectation operator in the Bellman equation of the time-discretized problem~\cite{CarmLud}.

The algorithm proposed in this paper is based on the work in \cite{RAid}, where a numerical scheme for solving high dimensional optimal switching problems was outlined, and a rate of convergence was derived. By utilizing the upper and lower bounds of the value function we will try to improve the efficiency of the method from \cite{RAid} to be able to manage problems of even higher dimension.

\subsection{Approximations}
The approximations made are
\begin{itemize}
  \item \emph{Time discretization} The intervention times $\tau_i$ will be restricted to a finite set of equally spaced points and the integral of the running-cost will be approximated by an Euler scheme.
  \item \emph{Switching restriction} We restrict switching by only allowing one switching request per time-point. 
  \item \emph{Running cost approximation} Techniques from machine learning will be used to approximate the running cost.
  \item \emph{Conditional expectation approximation} The conditional expectation operator in the Bellman equation of the time-discretized problem is approximated by Monte Carlo-regression on a set of local polynomials.
\end{itemize}

\subsubsection{Time discretization}
We start by going from continuous to discrete time by introducing the grid $\Pi=\{t_0,t_1,\ldots,t_{N_{\Pi}}\}$, with $t_m= m\Delta t$ for $m=0,\ldots, N_{\Pi}$, where $\Delta t=T/N_{\Pi}$. Here, it is convenient to chose $N_{\Pi}$ so that $\Delta t$ divides $\delta_\beta$, for all $\beta\in\Lambda$. To get a discrete time problem we reduce the set of stopping times in the admissible control so that switching requests are only allowed at grid points, \ie for all discretized intervention times $\bar\tau_j$ we have $\bar\tau_j\in \Pi$.
\begin{defn}
We define the set of admissible controls $\overline{\mcV}$ as the set of all $\bar v=(\bar\tau_1,\ldots,\bar\tau_N; \bar\alpha_1,\ldots,\bar\alpha_{N})\in\mcV$, with $\bar\tau_j\in \Pi$, for $j=1,\ldots,N$. Furthermore, for all $t_m\in\Pi$, we let $\overline{\mcV}_{t_m,a,r}$ be the set of all $v\in \mcV_{t_m,a,r}$ with $\tau_j\in\Pi$, for $j=1,\ldots,N$.
\end{defn}

\begin{defn}
Let $H_\Pi(a)=H(a)\cap \Pi$ be the discrete time version of the set of all possible transition progression vectors for the operating mode $a\in\Gamma$.
\end{defn}

\begin{defn}
For each $t_m\in \Pi$ we let $\mcT_{t_m}$ be the set of all $\mcF_t$--stopping times $\tau$, with $\tau\in\Pi\cup \infty$ and $\tau\geq t_m$.
\end{defn}

\subsubsection{Switching restriction}
To make the problem more computationally tractable we only allow single switching requests at each point in time.

To simplify notation when writing out the Bellman equation for the discrete-time problem we augment the sets of admissible switching requests with the empty set, as follows.
\begin{defn}
Let $\bar \Lambda:=\Lambda \cup \emptyset$ and for each $a\in\Gamma$ let $\bar \Lambda(a):=\Lambda(a) \cup \emptyset$
\end{defn}

The Bellman equation for the discrete-time version of Problem 1 then reads
\begin{align}
V_{\Pi}(T,x,a,r)&=0,\\
V_{\Pi}(t_m,x,a,r)&=f(x,a,r)\Delta t+\min \limits_{\beta\in \overline\Lambda(a)}\left\{K_{\beta}+\E\left[V_{\Pi}(t_{m+1},X^{t_m,x}_{t_{m+1}},A_{t_{m+1}}^{t_m,a+\beta,r},R_{t_{m+1}}^{t_m,a+\beta,r})\right]\right\}, \label{ekv:diskBellman}
\end{align}
with the convention that $K_\emptyset :=0$.

\subsubsection{Running cost approximation}
To compute the running cost we have to solve $n_c$ optimization problems. Since the number of contingencies $n_c$ that we want to include in the analysis may grow considerable with the size of the system\footnote{A common praxis is to include all contingencies with $k$ component failures, where $k$ may range from one up to three or four and sometimes even higher.}, this becomes intractable when trying to find the optimal load disruption $z^*$ while solving Problem~1. Instead we may, for example, apply machine learning to approximate the running costs.

A thorough discussion of machine learning algorithms for approximating minimal load disruption cost is out of the scope of this article. The interested reader are referred to the early textbook~\cite{WehenkelML} on application of machine learning to various problems in power system security.

\subsubsection{Conditional expectation approximation}
Following the seminal paper by Longstaff and Schwartz~\cite{Longstaff}, that introduced regression Monte Carlo in optimal stopping to a wider audience, a number of different Monte Carlo techniques have been proposed for the valuation of American options. Beginning with \cite{CarmLud} some of these Monte Carlo techniques have also been applied to solve optimal switching problems.

For high-dimensional optimal switching problems the predominant algorithms are based on Monte Carlo regression combined with either, policy approximation\cite{CarmLud} or value function approximation\cite{RAid}.

In policy approximation the value function in the expected value on the right hand side of \eqref{ekv:diskBellman} is approximated by regression on the cost-to-go of each trajectory when following the approximated optimal policy. Applied to our problem, policy approximation would give
\begin{align}
\hat V_{\Pi}^{\text{P}}(T,x,a,r)&=0,\\
\hat V_{\Pi}^{\text{P}}(t_m,x,a,r)&=\min \limits_{\beta\in \overline\Lambda(a)}\Bigg\{ K_{\beta}+\hatE\Bigg[\sum\limits_{j=m}^{N_\Pi}f(X_{t_j}^{t_m,x},\alpha^{t_m,a,r}(t_j,\hat v^*),\xi^{t_m,a,r}(t_j,\hat v^*))+\sum_{t_m<\hat \tau_j < T} K_{\beta_j}\Bigg]\Bigg\}, \label{ekv:policyAPPR}
\end{align}
where $\hatE$ is an approximation of the conditional expectation operator and $\hat v^*$ is the approximation of the optimal control obtained through the scheme.

Value function approximation is based on using an estimate of the value function on the right hand side of \eqref{ekv:diskBellman}, \ie the value function at time $t_{m+1}$, when computing the value function at time $t_m$. When applying value function approximation to our problem we get
\begin{align}
\hat V_{\Pi}^{\text{VF}}(T,x,a,r)&=0,\\
\hat V_{\Pi}^{\text{VF}}(t_m,x,a,r)&=f(x,a,r)\Delta t+\min \limits_{\beta\in \overline\Lambda(a)}\Big\{K_\beta+\hatE\left[\hat V_{\Pi}^{\text{VF}}(t_{m+1},X_{t_{m+1}}^{t_m,x},A_{t_{m+1}}^{t_m,a+\beta,r},R_{t_{m+1}}^{t_m,a+\beta,r})\right]\Big\}. \label{ekv:vfAPPR}
\end{align}
Policy approximation is biased towards overestimation of the value function, since policy approximation evaluates the approximated (and hence sub-optimal) policy. Value function approximation, on the other hand, gives a bias towards underestimation of the optimal cost. To realize this, assume that the conditional approximation expectation is $\mcF_{t_m}$--adapted and unbiased in the sense that $\E[\hatE [\cdot | \mcF_{t_m}]|\mcF_{t_m}]=\E[\cdot|\mcF_{t_m}]$. Then, by applying Jensen's inequality,
\begin{align*}
\E[\hat V_{\Pi}^{\text{VF}}(t_m,x,a,r)]&=f(x,a,r)\Delta t+\E\left[\min \limits_{\beta\in \overline\Lambda(a)}\Big\{K_\beta+\hatE\left[\hat V_{\Pi}^{\text{VF}}(t_{m+1},X_{t_{m+1}}^{t_m,x},A_{t_{m+1}}^{t_m,a+\beta,r},R_{t_{m+1}}^{t_m,a+\beta,r})\right]\Big\}\right]\\
&\leq f(x,a,r)\Delta t+\min \limits_{\beta\in \overline\Lambda(a)}\Big\{K_\beta+\E\left[\hatE\left[\hat V_{\Pi}^{\text{VF}}(t_{m+1},X_{t_{m+1}}^{t_m,x},A_{t_{m+1}}^{t_m,a+\beta,r},R_{t_{m+1}}^{t_m,a+\beta,r})\right]\right]\Big\}.
\end{align*}
Backward induction combined with the fact that $\hat V_{\Pi}^{\text{VF}}(T,\cdot,\cdot,\cdot)=V_{\Pi}(T,\cdot,\cdot,\cdot)$ now implies that
\begin{equation}
\E[\hat V_{\Pi}^{\text{VF}}(t_m,x,a,r)]\leq V_{\Pi}(t_m,x,a,r),\quad \forall t_m\in\Pi,\, x\in\R^{n_G},\,r\in \bar H(a),\,\forall a\in \Gamma.
\end{equation}

In~\cite{Bouchard2012} policy approximation and value function approximation are numerically compared on several multi-dimensional American options problems and policy approximation seems to give faster convergence in all of the investigated cases.

The drawback of applying policy approximation to solve problems with delays, such as Problem 1, is that we are forced to compute partial state trajectories corresponding to the estimated optimal control. If, as in~\cite{CarmLud}, the state trajectory is independent of the applied controls this does not pose any problems. For problems with delays, after augmentation, the delay time is part of the state. Hence, for each time, $t_m\in\Pi$, we may have to re-calculate the optimal cost-to-go during entire delay periods as the optimal trajectory may change.

Instead, a mixture of policy approximation and value function approximation may be adequate. This would give rise to the dynamic programming algorithm,
\begin{align}
\hat V_{\Pi}^{\text{M}}(T,x,a,r)&=0,\nonumber\\
\hat V_{\Pi}^{\text{M}}(t_m,x,a,r)&=\min \limits_{\bar\tau\in \overline{\mcT}_{t_{m+1}}}\hatE\Bigg[\sum_{t_{m}\leq t_j<\bar\tau} f(X_{t_j}^{t_m,x},A_{t_j}^{t_m,a,r},R_{t_j}^{t_m,a,r})\Delta t \\
&\hspace{100pt}  +\min\limits_{\beta\in\Lambda(A_{\bar\tau}^{t_m,a,r})}\left\{K_\beta+\hat V_{\Pi}^{\text{M}}(\bar\tau,X_{\bar\tau}^{t_m,x},A_{\bar\tau}^{t_m,a,r}+\beta,R_{\bar\tau}^{t_m,a,r})\right\}\Bigg]\label{ekv:mixAPPR}.
\end{align}
Note that, unlike the preceding methods, the mixture does not guarantee an under- or an over-bias of the estimate of the optimal cost-to-go. During periods where no switches are made on a certain trajectory the optimal running cost for that trajectory will be overestimated since a sub-optimal control is applied. On the other hand, when a switch is made the cost-to-go may be underestimated due to the concavity of the min operator, as mentioned above.\\

Now to the approximation of the conditional expectation operator. For each $a\in \Gamma$ we define a set of basis functions $b_j^a:\R^{n_L}\times H_{\Pi}(a)\to\R$, for $j=1,\ldots,M^a$ with $M^a\in \mathbb{N}$ and let, for each $\varphi:\Pi\times \R^{n_L}\times H_{\Pi}(a)\to\R$
\begin{align*}
\gamma^a_{t_m}(\varphi)=\mathop{\arg\max}\limits_{\gamma\in\R^{M^a}}\E\left[\left(\varphi(t_{m+1},X_{t_{m+1}},R_{t_m})-\sum_{j=1}^{M^a}\gamma_j b_j^a(X_{t_m},R_{t_m})\right)^2\right].
\end{align*}

Using the basis functions and the simulated trajectories $(X_{t_m}^l)_{1\leq m\leq N_\Pi}^{1\leq l \leq M_T}$ and $(R_{t_m}^l)_{1\leq m\leq N_\Pi}^{1\leq l \leq M_T}$ we can estimate $\gamma^a_{t_m}(\varphi)$ as
\begin{align}
\hat\gamma^a_{t_m}(\varphi)=\mathop{\arg\max}\limits_{\gamma\in\R^{M^a}}\frac{1}{M_T}\sum_{l=1}^{M_T}\left(\varphi(t_{m+1},X^l_{t_{m+1}},R^l_{t_m})-\sum_{j=1}^{M^a}\gamma_j b_j^a(X^l_{t_m},R^l_{t_m})\right)^2.\label{ekv:leastSQ}
\end{align}

As pointed out in~\cite{RAid} it is often necessary to truncate the conditional expectations. This can be done using known lower and upper bounds $\varphi^{\rm{LB}}$ and $\varphi^{\rm{UB}}$ of $\varphi$ for which $\Upsilon^{t_m,x,r}(\varphi^{\rm{LB}}):=\E\left[\varphi^{\rm{LB}}(t_{m+1},X^{t_m,x}_{t_{m+1}},r)\right]$ and $\Upsilon^{t_m,x,r}(\varphi^{\rm{UB}}):=\E\left[\varphi^{\rm{UB}}(t_{m+1},X^{t_m,x}_{t_{m+1}},r)\right]$ can be computed.

We then get the regression Monte Carlo estimate
\begin{equation*}
\hatE\left[\varphi(t_{m+1},X_{t_{m+1}}^{t_m,x},r)\right]:=\Upsilon^{t_m,x,r}(\varphi^{\rm{LB}}) \vee \sum_{j=1}^{M^a}(\hat\gamma^a_{t_m})_j b_j^a(x,r) \wedge \Upsilon^{t_m,x,r}(\varphi^{\rm{UB}})
\end{equation*}

In~\cite{Bouchard2012} the use of local polynomials was suggested as these can naturally be extended to high dimensional problems, and the corresponding least-squares problem in \eqref{ekv:leastSQ} becomes sparse. We will later apply local polynomial bases when solving actual optimal operation problems but let us for now assume that the basis is made up of a number of indicator functions $b_j=1_{B_j}$ where $(B_j)_{1\leq j\leq M}$ is a partition of $\R^{n_L}\times H(a)$. Then $\hat\gamma^a_{t_m}(\varphi)$ is given by
\begin{equation*}
(\hat\gamma^a_{t_m})_j(\varphi)=\frac{\sum_{l=1}^{M_T}\varphi(t_{m+1},X^l_{t_{m+1}},R^l_{t_m})1_{B_j}(X^l_{t_{m}},R^l_{t_m})}{\sum_{l=1}^{M_T}1_{B_j}(X^l_{t_{m}},R^l_{t_m})}.
\end{equation*}
This gives the variance
\begin{align*}
\text{Var}\left[(\hat\gamma^a_{t_m})_j(\varphi)\right]=\frac{1}{M_T^{j,t_m}}\text{Var}\left[\left.\varphi(t_{m+1},X_{t_{m+1}},R_{t_m})\right|(X_{t_{m}},R_{t_m})\in B_j\right],
\end{align*}
where $M_T^{j,t_m}$ is the expected number of sample paths for which $(X^l_{t_{m}},R^l_{t_m})\in{B_j}$. Concerning the $L^2$--error of the estimate we have, before truncation,
\begin{align*}
\|(\hat\gamma^a_{t_m})_j-\E&\left[\varphi(t_{m+1},X^{t_m,x}_{t_{m+1}},r)\right]\|^2\\
&=\E\left[\left((\hat\gamma^a_{t_m})_j-\E\left[\varphi(t_{m+1},X^{t_m,x}_{t_{m+1}},r)\right]\right)^2\right]\\
&=\E\left[\left((\hat\gamma^a_{t_m})_j-(\bar\gamma^a_{t_m})_j\right)^2\right] + \left((\bar\gamma^a_{t_m})_j-\E\left[\varphi(t_{m+1},X^{t_m,x}_{t_{m+1}},r)\right]\right)^2.
\end{align*}
Hence, the efficiency of the estimate depends on the incompleteness of the set of basis functions, but also on the variance of the vector $\hat\gamma^a_{t_m}(\varphi)$.

To improve the efficiency of the algorithm when valuing American options, it was suggested in~\cite{Bouchard2012} that the corresponding European option value could be used as a control variable. For more general optimal switching problems it may be difficult to find a suitable counterpart, however, for Problem~1 the lower bound may serve as an adequate substitute.

Hence, instead of estimating $\E\left[\varphi(t_{m+1},X^{t_m,x}_{t_{m+1}},r)\right]$ we may want to estimate\\ $\E\left[\varphi(t_{m+1},X^{t_m,x}_{t_{m+1}},r)-\varphi^{\rm{LB}}(t_{m+1},X_{t_{m+1}},r)\right]$, to get the new approximation
\begin{align*}
\hatE_{\rm{CV}}^{t_m,x}\left[\varphi(t_{m+1},X_{t_{m+1}},r)\right]:= &\Upsilon^{t_m,x,r}(\varphi^{\rm{LB}}) + \bigg(0 \wedge \sum_{j=1}^{M^a}(\hat\gamma^a_{t_m})_j(\varphi-\varphi^{\rm{LB}}) b_j^a(x,r)\\
&\vee (\Upsilon^{t_m,x,r}(\varphi^{\rm{UB}})-\Upsilon^{t_m,x,r}(\varphi^{\rm{LB}}))\bigg).
\end{align*}
As long as $\text{Var}\left[(\hat\gamma^a_{t_m})_j(\varphi-\varphi^{\rm{LB}})\right]<\text{Var}\left[(\hat\gamma^a_{t_m})_j(\varphi)\right]$, this will lead to an improvement in the efficiency of the estimate.

To apply this in our algorithm we let $\varphi(t_{m+1},x,r)=V_\Pi(t_{m+1},x,A_{t_m+1}^{t_m,a,r},R_{t_m+1}^{t_m,a,r})$ and $\varphi^{\rm{LB}}(t_{m+1},x,r)=V_\Pi^\eta(t_{m+1},x,A_{t_m+1}^{t_m,a,r},R_{t_m+1}^{t_m,a,r})$, for $\eta\in U$. With $\Delta V_\Pi^{\eta} := V_\Pi- V^{\eta}_\Pi$, we get, for all $\eta\in U$,
\begin{align*}
\Delta V_\Pi^{\eta}(T,x,a,r)=\,&0,\\
\Delta V_\Pi^{\eta}(t_m,x,a,r)=\,&\Delta f^{\eta}(x,a,r)\Delta t + \min\limits_{\beta\in \overline\Lambda(a)}\big\{K_{\beta}+\E\big[V^{\eta}_\Pi(t_{m+1},\eta_L^\top X^{t_m,x}_{t_{m+1}},A_{t_{m+1}}^{t_m,a+\beta,r},S(R_{t_{m+1}}^{t_m,a+\beta,r}))
\\
&+\Delta V_\Pi^{\eta}(t_{m+1},X^{t_m,x}_{t_{m+1}},A_{t_{m+1}}^{t_m,a+\beta,r},R_{t_{m+1}}^{t_m,a+\beta,r})\big]\big\}
\\
- \min\limits_{k\in\overline{\mathcal{K}}^T_{a,S(r),0}}&\min\limits_{\beta\in \overline\Lambda(a^{t,k})}\big\{K_{\beta}+\E\big[V^{\eta}_\Pi(t_{m+1},\eta_L^\top X^{t_m,x}_{t_{m+1}}, a^{t,k}+\beta,S(r)+N_T(a^{t,k}+\beta)\Delta t-\delta_{k,a_k})\big]\big\},
\end{align*}
where $\overline{\mathcal{K}}^T_{a,S(r),0}:={\mathcal{K}}^T_{a,S(r),0}\cup\emptyset$ and
\begin{equation*}
\Delta f^{\eta}(x,a,r):=f(x,a,r)-f^{\eta}((\eta_L)^\top x,a,S(r)).
\end{equation*}

In the numerical example we will compare the control variable approximation $V^{\eta}_\Pi + \Delta\hat{ V}_\Pi^{\text{M}}$ to the regular approximation $\hat{V}_\Pi^{\text{M}}$ for a small scale problem against a solution obtained by Markov-Chain approximation.

The choice of basis functions can be categorized into two classes. In spectral methods the basis functions have global support, while in finite element methods the support of each basis function is restricted to one of the sets in a partition of the state space. Many suggestions of global basis polynomials appear in the literature and good results can be obtained for specific problems, but as pointed out in~\cite{Bouchard2012} a choice of appropriate global basis polynomials is made very difficult as the dimension of the problem increases. In~\cite{Bouchard2012} and~\cite{RAid} basis polynomials with local support are used, motivated both by the fact that they scale naturally to higher dimensions and that the least squares problem in \eqref{ekv:leastSQ} becomes sparse, making computation efficient. A more detailed discussion on the benefits of local support basis functions can be found in~\cite{Bouchard2012}.

When partitioning the state space it is reasonable the make sure that each subset in the partition gets approximately the same number of particles. The partition of the state space may be either deterministic, \ie based on the density of\footnote{Or rather an estimate of the density as the density of $R_{t_m}$ is unknown prior to knowing the optimal control on $\{0,t_1,\ldots,t_{m-1}\}$.} $(X_{t_m},R_{t_m})$, or adapted to the trajectories of the particles.

\subsection{Algorithm}
The algorithm implemented to solve Problem 1 is outlined in Algorithm~\ref{TheAlg}. The remainder of this section is aimed at detailing and analyzing the different steps of Algorithm~\ref{TheAlg}.

\begin{algorithm}[h!]
 \KwData{lower and upper bounds on the value function, market data, expected costs from load disruption}
 \KwResult{optimal operation strategy}
 randomize $R_0$\;
 \For {$m=1:N_\Pi-1$}{
 randomize $X_{t_{m+1}}$\;
 save seed($m$)\;
 clear $X_{t_{m}}$\;
 }
 \For {$m=N_\Pi:-1:1$}{
  \If {$m<N_\Pi$}{
   using seed($m$) and $X_{t_{m+1}}$ recover $X_{t_m}$\;
  }
  using $R_0$ compute $R_{t_m}$\;
  \For {$a\in \Gamma$}{
    \For {$l:=1$ \emph{to} $M_T$}{
      running cost := find running cost in $(X^l_{t_m},a,R^{a,l}_{t_m})$\;
      \eIf {$m=N_\Pi$}{
        cost-to-go($l$) = running cost\;
      }{
        $a_{\text{next}}:=A_{t_{m+1}}^{t_m,a,R^{a,l}_{t_m}}$\;
        cost optimal next := get approximate cost-to-go in $(t_{m+1},X^l_{t_{m+1}},a_{\text{next}},R^{a,l}_{t_{m+1}})$\;
        \For {$\beta \in \Lambda(a_{\rm{next}})$}{
          cost if switch := $K_\beta$ + get approximate cost-to-go in $(t_{m+1},X^l_{t_{m+1}},a_{\text{next}}+\beta,R^{a,l}_{t_{m+1}})$\;
          \If {\emph{cost if switch} $<$ \emph{cost optimal next}}{
            cost optimal next := cost if switch\;
          }
        }
        cost-to-go($l$) := running cost + cost optimal next\;
      }
    }
    set approximate cost-to-go in mode $a$\;
  }
}
 \caption{Solving the optimal switching problem\label{TheAlg}}
\end{algorithm}

\subsubsection{Computing the lower and upper bounds}
When computing the lower bound we first pick a set $U^\eta=\{\eta^1,\ldots,\eta^{N_\eta}\}\subset S^{n_p-1}$, with $[\eta^j_{n_G+1}\:\cdots\: \eta^j_{n_G+n_L}] > 0$ and find $\hat V^{\eta^1},\ldots,\hat V^{\eta^{N_\eta}}$. The lower bound is then the pointwise minimum of these functions.

The time complexity of computing $\hat V_\Pi^{\text{LB}}$ is $\bigO (N_\eta N_\Pi^2 N_\Gamma n_G N_x^\eta)$ and for $\hat V_\Pi^{\text{UB}}$ it is $\bigO (N_\Pi^2 N_{\Theta} 2^{n_G} N_x^B)$, where $N_\Gamma =|\Gamma|$ and $N_x^\eta$ and $N_x^B$ are the sizes of the grids used to discretize the one-dimensional state spaces of $X^\eta$ and $X^B$ respectively.

\subsubsection{Drawing sample trajectories}
Let $H:=\mathop{\prod}\limits_{k=1}^{n_G}\left(\mathop{\prod}\limits_{(l,m)\in \Lambda_k}[0,\delta_{k,(l,m)}]\right)$ and $d_H:=\sum\limits_{k=1}^{n_G}|\Lambda_k|$, so that $H\subset \R_+^{d_H}$. We define $H_{\Pi}:=\{r\in H:\: r_j\in \Pi,\,j=1,\ldots,d_H \}$. Furthermore we assume that the sets $\Lambda_k$ are ordered so that for each $(l,m)\in\Lambda_k$ and $(l',m')\in\Lambda_k$, with $(l,m)\neq(l',m')$, either $(l,m)<(l',m')$ or $(l,m)>(l',m')$. This allows us to define a correspondence between the elements of vectors in $H_{\Pi}$ and the transition progression vector through the bijective map $\iota:\Lambda\to \{1,\ldots,d_H\}$, given by
\begin{equation}
\iota(k,(l,m)):=\sum\limits_{k'=1}^{k-1}|\Lambda_{k'}|+|\{(l',m')\in\Lambda_k:\: (l',m')\leq (l,m)\}|.
\end{equation}

To approximate the conditional expectations in the Bellman equation, a number of sampled state-trajectories $(X_{t_m}^l,R_{t_m}^l)_{1\leq m \leq N_\Pi}^{1\leq l \leq M_T}$ are used. The demand sample trajectories, $X_{t_m}^l$ are created from the demand SDE but the corresponding samples of $\xi$ cannot be obtained at this stage as this would require knowing the optimal control beforehand\footnote{This is exactly why the policy approximation algorithm becomes intractable.}. Instead we create samples\footnote{The notation $R_{t_m}^l$ is used to emphasize that samples do not correspond to any specific control.} $R_{t_m}^l\in H_\Pi$ for each $t_m\in \Pi$ and $l=1,\ldots,M_T$ from a different distribution. This is done by picking $R_{0}^l$ uniformly in $H_\Pi$ and letting $(R_{t_m}^l)_j:=\left[(R_{0}^l)_j+t_m\right] \mod \delta_{\iota^{-1}(j)}$.

Now, we define the map $\bar \iota:\Gamma \to 2^{\{1,\ldots,d_H\}}$ as
\begin{equation}
\bar \iota(a):=(\iota(\kappa^T_1,a_{\kappa^T_1}),\ldots,\iota(\kappa^T_{N_T(a)},a_{\kappa^T_{N_T(a)}})),
\end{equation}
where $\kappa^T\in 2^{\Node_G}$ is an ordered enumeration the set $\mathcal{K}^T(a)$, \ie the set of all generators that are in transition in mode $a$. Furthermore, we let $R_{t_m}^{a,l}$ be the vector with
\begin{equation*}
(R_{t_m}^{a,l})_k:=\left\{\begin{array}{l l} (R_{t_m}^{l})_{\iota(k,a_k)},\quad &\text{if }k\in K^T(a), \\
0,& \text{otherwise},\end{array}\right.
\end{equation*}
for $k=1,\ldots,n_G$, so that $R_{t_m}^{a,l}$ is uniformly distributed on $H^\Pi(a)$.

Creating the sample trajectories is done in $\bigO(M_T n_L N_\Pi + M_T d_H)$. Concerning the memory requirement, as described in~\cite{RAid}, we do not need to store the entire set of trajectories, only the seeds used at each time instance in the Euler scheme. The trajectories $(X_{t_m}^l)_{1\leq m\leq N_\Pi}^{1\leq l\leq M_T}$ can then be recalled in the backward induction scheme when computing the value function approximations.

\subsubsection{State space partition}
To define the set of local basis functions we need to make a partition of the state space. A tensor partition, based on an interval partition of each state-space dimension, was suggested in~\cite{Bouchard2012}. This approach was shown to be successful in examples of up to six dimensions. In higher dimensions the curse of dimensionality becomes noticeable as the number of hypercubes in the partition grows exponentially. To prevent this explosion in computational burden, in the form of an exponentially increasing number of samples needed for the algorithm to converge, we use an area partition of the system to define the state space partition. The rationality behind this is that, if the area partition of the system is wisely chosen, the net import/export from each area is generally an adequate indicator of the stress that the system is under.

We, thus, assume that the system is divided into $n_A$ areas based on \eg net structure, geography or electrical distance. We let $I:=[I^D_1\:\cdots\:I^D_{n_A}\:\: I^R_1\:\cdots\:I^R_{d_H}]^\top \in \Zed^{n_A+d_H}$ be the number of intervals in the partition of each dimension. For $l=1,\ldots,n_A$ we define the intervals $[-\infty, y^{D}_{l,1}]\times [y^{D}_{l,1}, y^{D}_{l,2}]\times\cdots \times [y^{D}_{l,I^D_l-1},\infty]$ and for $k=1,\ldots,d_H$ we define the intervals $[-\infty, y^{R}_{k,1}]\times [y^{R}_{k,1}, y^{R}_{k,2}]\times\cdots \times [y^{R}_{k,I^R_k-1},\infty]$. For each $a\in\Gamma$ we divide the state space of $(X_t,R^a_t)$ into $\prod_{l=1}^{n_A} I^D_l\prod_{k\in\mathcal{K}^T(a)}I^{R}_{\iota^{-1}(k,a_k)}$ rectangular sets
\begin{align*}
D_{i^D_1\cdots i^D_{n_A}i^{R,a}_1\cdots i^{R,a}_{N_T(a)}}:=\bigg\{(x,r)\in \R^{n_L}\times H^\Pi(a) :\: y^{D}_{l,i^D_l-1}&\leq\sum_{j\in A_l}x_j < y^{D}_{l,i^D_l},\,l=1,\ldots,n_A\\
y^{R}_{(\bar\iota(a))_k,i^{R,a}_k-1}&\leq r_{(\bar\iota(a))_k} < y^{{R}}_{(\bar\iota(a))_k,i^{R,a}_k},\,k=1,\ldots,N_T(a),
\bigg\},
\end{align*}
where $A_l\subset \{1,\ldots,n_N\}$ is the subset of nodes in Area $l$. The local basis is then given by functions of the form
\begin{equation*}
1_{D_{i^D_1\cdots i^D_{n_A}i^{R,a}_1\cdots i^{R,a}_{N_T(a)}}}(x,r)\tilde b^a(x,r),
\end{equation*}
where $\tilde b^a:\R^{n_L}\times H(a)\to \R$ is a polynomial.

When defining the intervals in the partition we would like to get approximately the same number of samples in each hypercube. We can either use the adaptive method proposed in~\cite{Bouchard2012} where, at each time step, partial sorts of the samples are performed and intervals are chosen such that exactly the same number of samples falls in each interval. Alternatively the density of $(X_t,R_t)$ can be used to define intervals with the same expected number of samples in each interval.

If the correlation is week between the aggregate load in different areas tensor partitions will give approximately the same number of samples in each hypercube.

The computational complexity of computing the hypercubes using a partial sorting algorithm is $\bigO (\sum_{j=1}^{n_A+d_H}(I_j-1)M_T)$. With a density based (non-adaptive) partition the complexity is $\bigO (\sum_{j=1}^{n_A+d_H}(I_j-1))$

If the non-adaptive method is applied we are able to use different grids for different $a\in\Gamma$ and still retain computational tractability. This would allow us to keep a fixed expected number of particles in each hypercube irrespective of the dimensionality of the state space corresponding to the operating mode.

\subsubsection{Dynamic programming}
In the dynamic programming algorithm we build the value function by, for each time $t_m\in\Pi\setminus \{T\}$ and each operating mode $a\in \Gamma$, finding the optimal $\mcF_{t_m}$--measurable decision for each $(X_{t_m}^l,R_{t_m}^{a,l})^{1\leq l\leq M_T}$ with respect to the approximated costs-to-go in the corresponding states.

To do this we need to evaluate the $|\Lambda(a)|+1=\bigO(n_G)$ different feasible actions that we can take when the operating mode is $a\in\Gamma$. Assuming that the particles have already been sorted into their respective hypercubes this has computational complexity of order $\bigO(M_T n_G M^a_{\textit{pol}})$, where $M_{\textit{pol}^a}$ is the number of polynomials, $\tilde b^a$, used in the approximation of the conditional expectation.

The overall computational complexity of this step is thus $\bigO(M_T N_\Pi N_\Gamma n_G M^a_{\textit{pol}})$.

\subsubsection{Conditional expectation approximation}
At the core of the algorithm is solving the least squares problem giving the coefficients in the approximation of the conditional expectations of the value function. The generic least squares problem can be written
\begin{equation*}
  \min\limits_{\gamma\in \R^{M}} \|A\gamma-b\|,
\end{equation*}
where in our problem $A$ is a sparse $M_T\times M^a$--matrix when approximating the conditional expectation of the value function in operating mode $a\in\Gamma$. The optimal solution to this problem is found by solving the normal equation $A^\top A \gamma =A^\top b$~\cite{BV03}. The main computational effort here goes into computing the matrix $A^\top A$. Since each row of $A$ only has $M_{\textit{pol}}^a$ non-zero elements, the computation time for $A^\top A$ is of order $\bigO(M_T M_{\textit{pol}}^a)$. Furthermore, the matrix $A$ never needs to be stored. After computing $A^\top A$, the complexity of solving the normal equation is negligible.

The total computational complexity of conditional expectation approximations is $\bigO(M_T N_\Pi\sum_{a\in\Gamma}M_{\textit{pol}}^a)$.

\subsubsection{Overall computational complexity}

The overall computational complexity is thus $\bigO (N_\eta N_\Pi^2 N_\Gamma n_G N_x^\eta+N_{\Theta} 2^{n_G} N_x^B + M_T N_\Pi\Big(\sum _{j=1}^{n_A+d_H}(I_j-1)$\\
$+n_L+N_\Gamma n_G M^a_{\textit{pol}}\Big))$.

\section{Numerical example\label{Sec:NumEX}}
In this section we consider two numerical examples where the proposed method is applied to solve production planning problems.

The aim of the first example is to show the different concepts of the problem. The size of the underlying system is small enough that a tractable solution can be found by a Markov-chain approximation method. This gives us a benchmark to which the outlined algorithm can be compared.

In the second example we investigate a more realistic situation where the number of loads and the number of operating modes are substantially higher.


\subsection{A 7-node system with two controllable generators}
In the first example we want to optimize the operation of a system with 7 nodes during an operating period of $T=60$ minutes.

\subsubsection{Network model}
\begin{figure}[h!]
  \centering
  \includegraphics[width=0.4\textwidth]{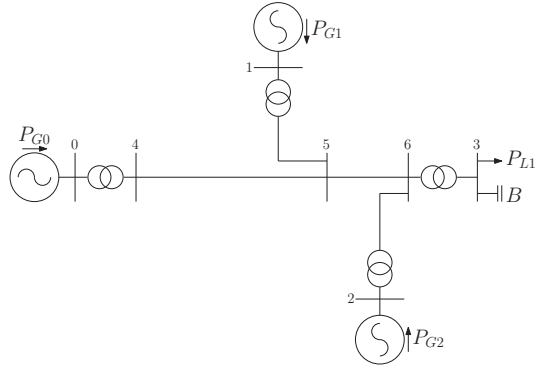}
  \caption{The 7-node system has three generators and one load.}\label{KsysFIG}
\end{figure}

The one-line diagram of the system is plotted in Fig.~\ref{KsysFIG}.
This system, taken from~\cite{Karystianos}, has three generators and one load. Node 0 is a slack node where production is balanced and the reactive power is limited in all three generators.

\subsubsection{Feasible sets and contingencies}
\begin{figure}
        \centering
        \begin{subfigure}[b]{0.48\textwidth}
                \includegraphics[width=\textwidth]{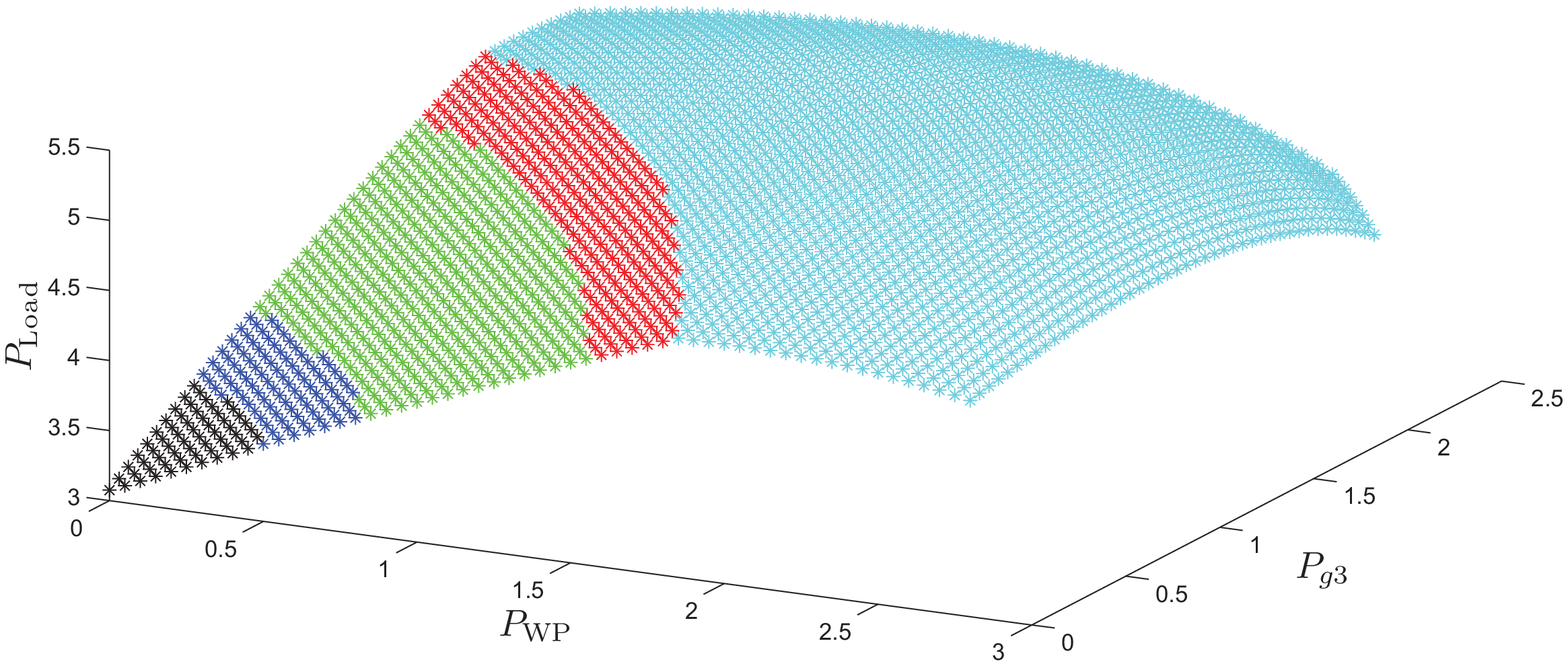}
                \caption{Pre-contingency}
                \label{fig:cont0}
        \end{subfigure}%

        \begin{subfigure}[b]{0.48\textwidth}
                \includegraphics[width=\textwidth]{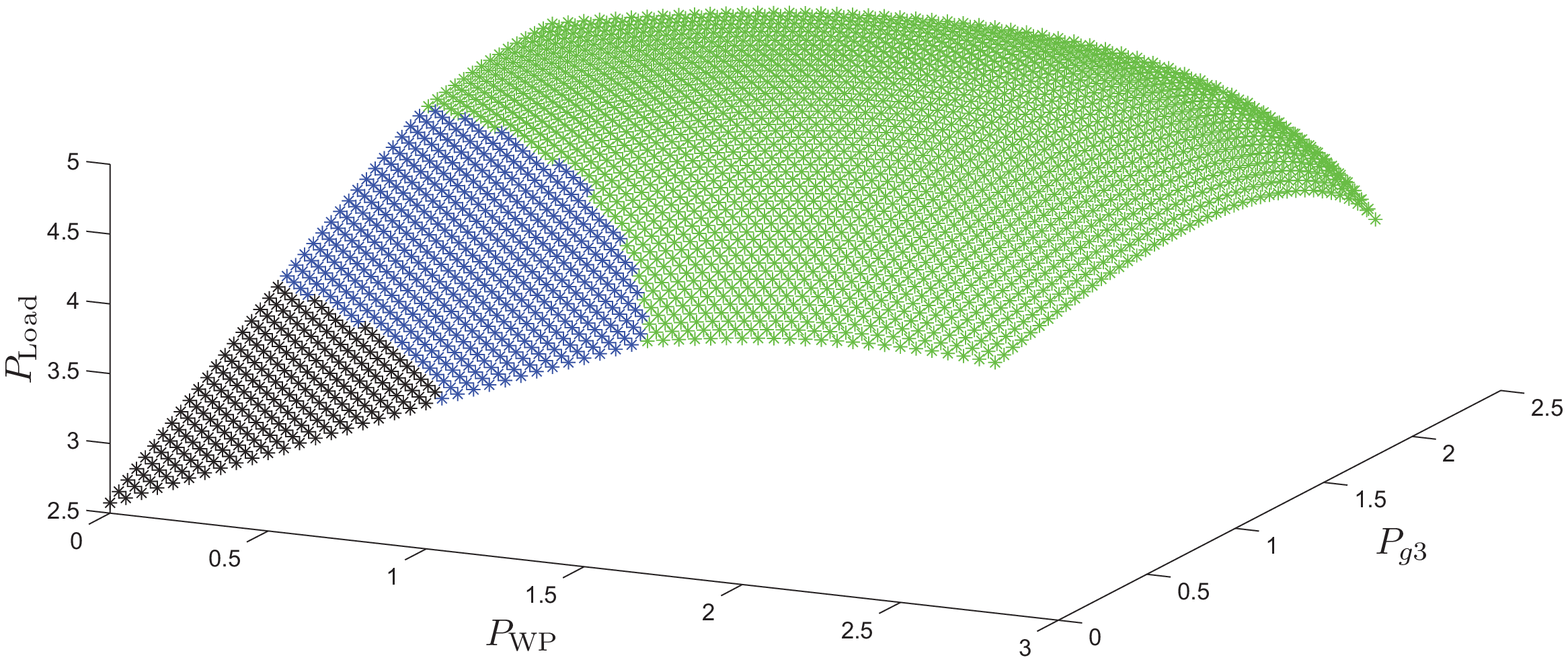}
                \caption{Post-contingency 1}
                \label{fig:cont1}
        \end{subfigure}
        ~ 
        \begin{subfigure}[b]{0.48\textwidth}
                \includegraphics[width=\textwidth]{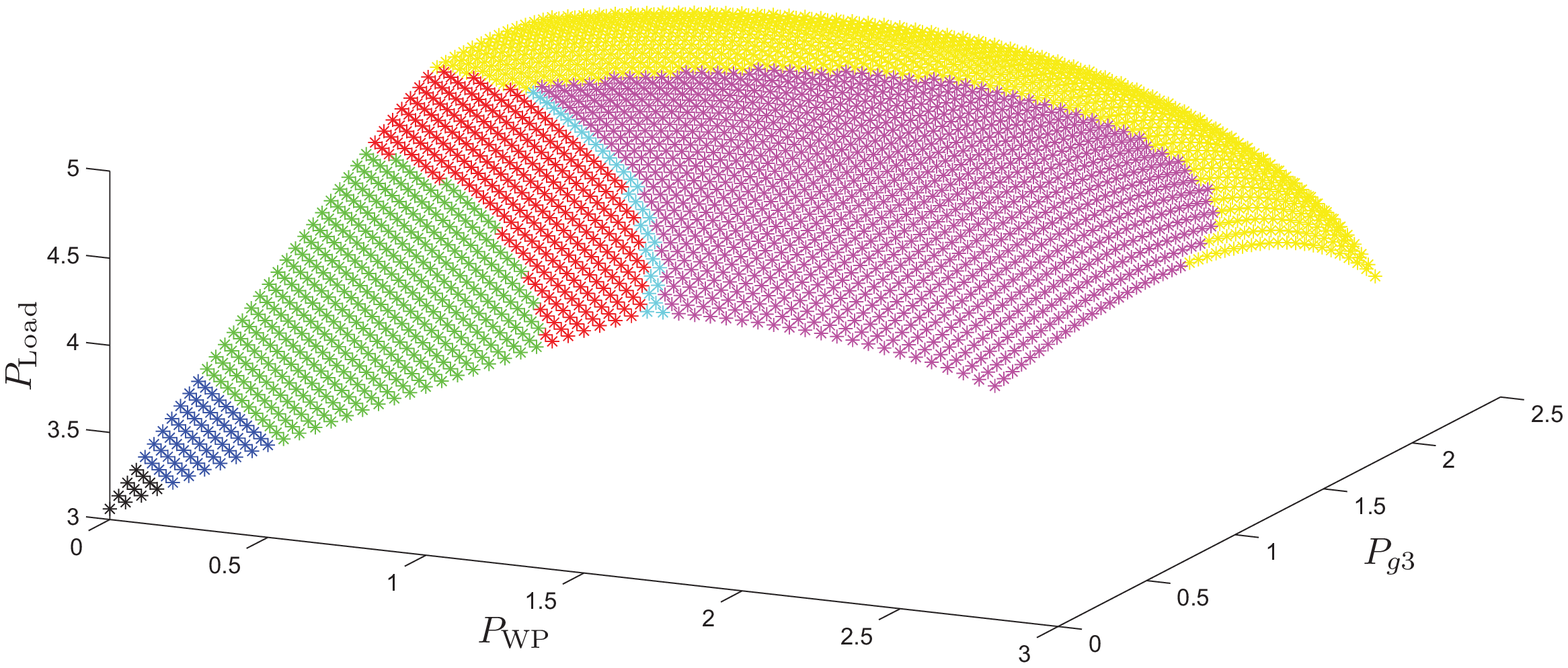}
                \caption{Post-contingency 2}
                \label{fig:cont2}
        \end{subfigure}
        \caption{Boundaries of the feasibility set in the different contingency states, \ie  $\partial G_i$, for $i=0,1,2$. Colors are used to represent different binding constraints in the power flow equations.}\label{fig:feasBOUNDksys}
\end{figure}

We consider two different contingencies. Due to the lack of redundancy in this small network we let the first contingency, case $i = 1$, be a doubling of the impedance of the transmission line connecting nodes 5 and 6. The second contingency case, $i = 2$ is a failure leading to the loss of the capacitor at the load bus. The boundaries of the feasibility domain for these different system configurations are plotted in Fig~\ref{fig:feasBOUNDksys}.

From the figures it appears that the feasible domains all meet the requirements of Assumption~\ref{AssG}. We set the failure weights to\footnote{The rates are deliberately exaggerated to make sure that contingencies have a noticeable impact on the optimal policies.} $q_0=1$ min$^{-1}$ and $q_i=0.01$ min$^{-1}$, for $i=1,2$.

\subsubsection{Stochastic model}
We assume that the load in Node 3 is given by the sum of the forecasted trajectory $m(t)=350+100\tfrac{t}{T}$ and the stochastic process $(Z_t:\,0\leq t\leq T)$ which solves
\begin{align*}
dZ_t&=-0.01 Z_tdt + 5 dW_t,\\
Z_0&=0,
\end{align*}
where $W_t$ is a standard Wiener process.

\subsubsection{Production modes and ramp rates}
We assume that the operator can switch between two production levels in Generator 1 and three production levels in Generator 2. The stationary production levels are given in Table~\ref{tab:thetaEX1}.

\begin{table}[h!]\centering \setlength{\tabcolsep}{10pt}
\begin{tabular}{| c | c | c | c | c |}\hline
  $\theta_{1,1}$ & $\theta_{1,2}$ & $\theta_{2,1}$ & $\theta_{2,2}$ & $\theta_{2,3}$ \\
  \hline
  100 & 150 & 0 & 100 & 200 \\ \hline
\end{tabular}
\caption{Stationary production levels [MW] in Example 1.}\label{tab:thetaEX1}
\end{table}
For the ramp function we assume that
\begin{equation}\label{ekv:RFex1}
\theta_{k,(l,m)}(s)=\theta_{k,l}+1_{[\delta'_{k,(l,m)},\delta_{k,(l,m)})}(s)\frac{s-\delta'_{k,(l,m)}}{\delta_{k,(l,m)}-\delta'_{k,(l,m)}} (\theta_{k,m}-\theta_{k,l}),
\end{equation}
where the delay times $\delta'_{k,(l,m)}$ and $\delta_{k,(l,m)}$ are given in Table~\ref{tab:deltasEX1}.
\begin{table}[h!]\centering \setlength{\tabcolsep}{10pt}
\begin{tabular}{| c | c | c | c | c | c | c |}\hline
  $k,(l,m)$ & $1,(1,2)$ & $1,(2,1)$ & $2,(1,2)$ & $2,(2,1)$ & $2,(2,3)$ & $2,(3,2)$ \\
  \hline
  $\delta'_{k,(l,m)}$ & 5 & 5 & 10 & 5 & 1 & 1 \\
  $\delta_{k,(l,m)}$ & 10 & 7.5 & 14 & 7 & 5 & 3 \\ \hline
\end{tabular}
\caption{Delay times in Example 1.}\label{tab:deltasEX1}
\end{table}

At the start of the operating period we assume that we are in the stationary mode $\alpha(0)=[1\:\: 1]^\top$.

\subsubsection{Operating costs}
The operating cost is the sum of the production cost and the expected cost of disconnecting load. We assume that the cost of disconnecting load is $c=10\,$k\currency/MW, that the production cost is $20\,$\currency/MWh~in Generator 1 and $25\,$\currency/MWh~in Generator 2 and that the switching costs are as in Table~\ref{tab:CsEX1}.

\begin{table}[h!]\centering \setlength{\tabcolsep}{10pt}
\begin{tabular}{| c | c | c | c | c | c | c |}\hline
  $k,(l,m)$ & $1,(1,2)$ & $1,(2,1)$ & $2,(1,2)$ & $2,(2,1)$ & $2,(2,3)$ & $2,(3,2)$ \\
  \hline
  $K_{k,(l,m)}\,$[\currency] & 500 & 250 & 2000 & 500 & 500 & 200 \\ \hline
\end{tabular}
\caption{Switching costs in Example 1.}\label{tab:CsEX1}
\end{table}

\subsubsection{Dimensions}
This problem has $2\cdot 3=6$ stationary modes. There are $2\cdot 4 + 2\cdot 3=14$ modes with one generator in transition and $2\cdot 4=8$ modes with both generators in transition, giving a total of 28 modes.

In the modes where two generators are in transition simultaneously there are $\mathcal{O} (\Delta t^{-2})$ possible values for the transition progression vectors (\ie $|H^\Pi(a)|=\mathcal{O} (\Delta t^{-2})$), giving a total computational time that is $\mathcal{O} ((N_\Pi)^3)$. Hence, already this small problem becomes very computationally demanding when classical approximation-methods are applied.

\subsubsection{Markov chain approximation}
\begin{figure}
        \centering
        \begin{subfigure}[b]{0.3\textwidth}
                \includegraphics[width=\textwidth]{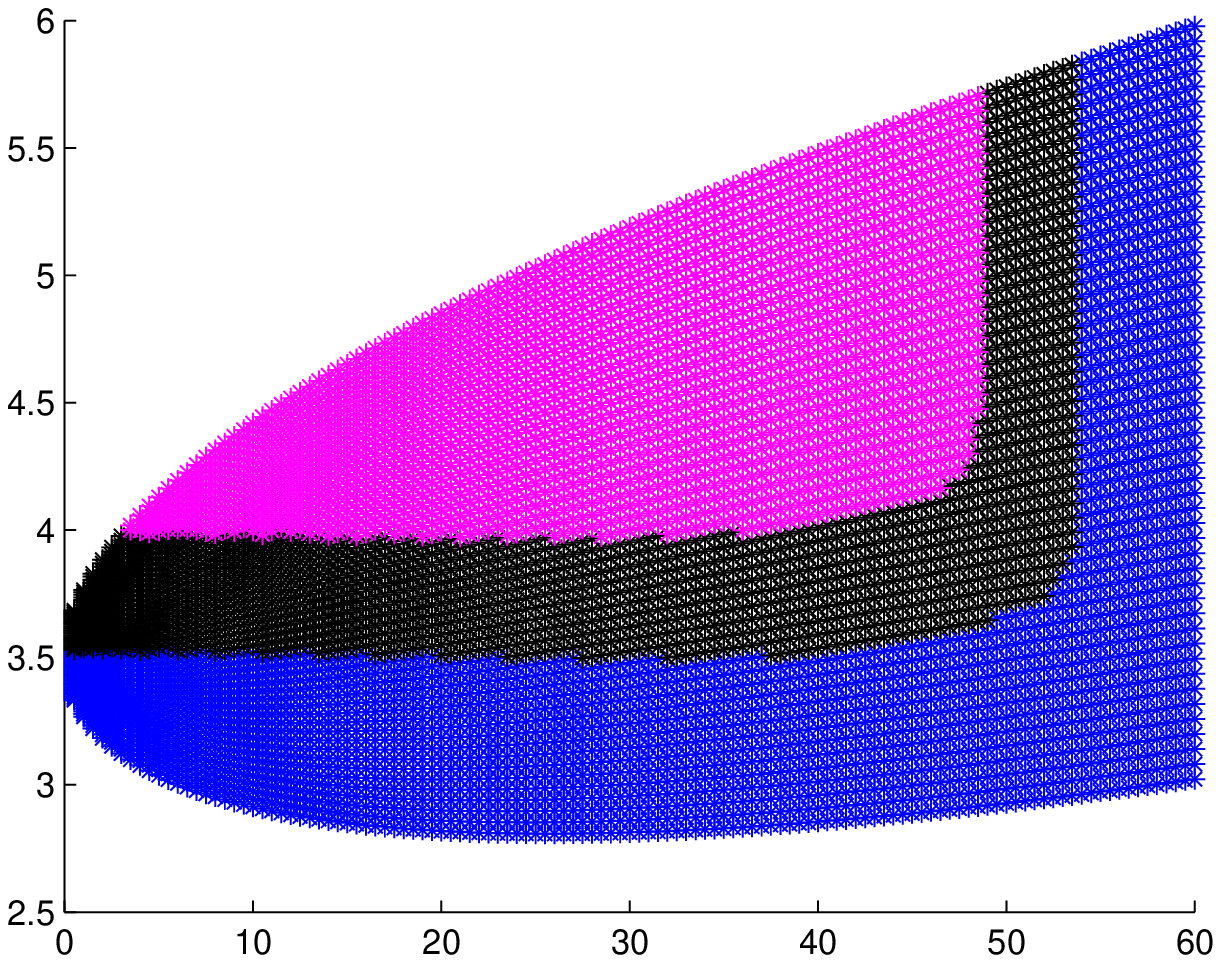}
                \caption{Switching regions from mode $[1\:\:  1]^\top$.}
                \label{fig:swRGS1}
        \end{subfigure}%
        \qquad
        \begin{subfigure}[b]{0.3\textwidth}
                \includegraphics[width=\textwidth]{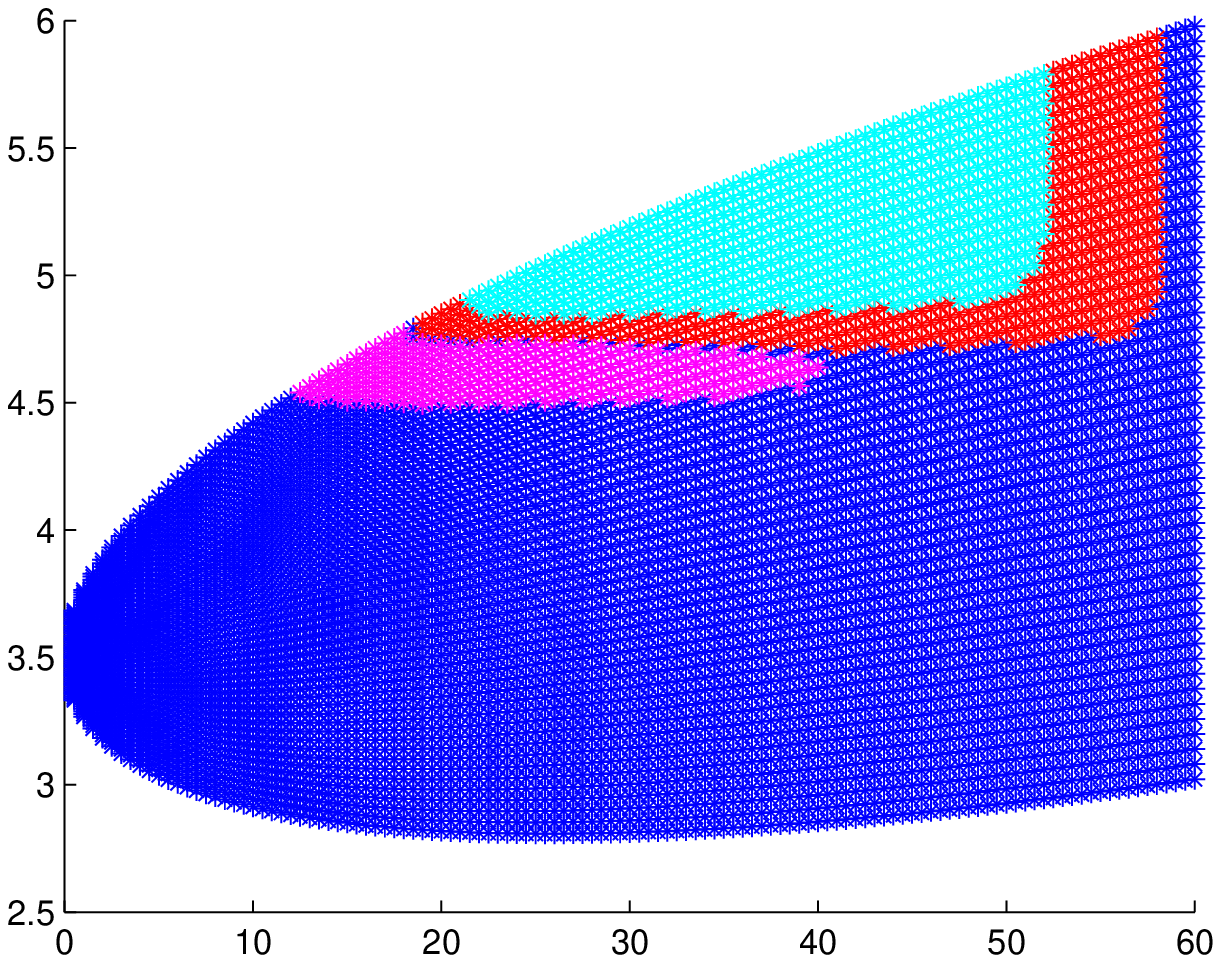}
                \caption{Switching regions from mode $[1\:\:  2]^\top$.}
                \label{fig:swRGS2}
        \end{subfigure}
        \qquad
        \begin{subfigure}[b]{0.3\textwidth}
                \includegraphics[width=\textwidth]{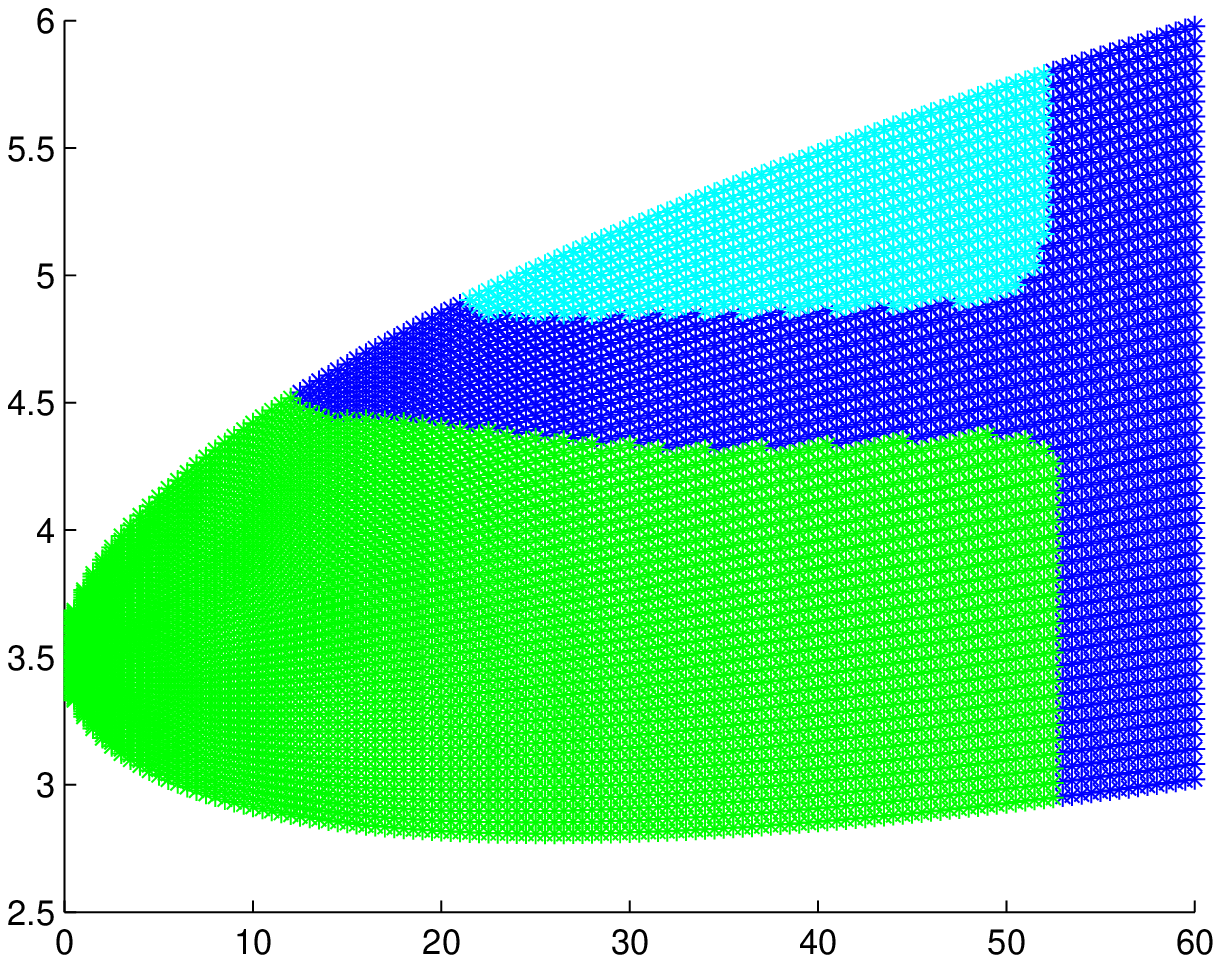}
                \caption{Switching regions from mode $[1\:\:  3]^\top$.}
                \label{fig:swRGS3}
        \end{subfigure}
        \vspace{0mm}
        
        \begin{subfigure}[b]{0.3\textwidth}
                \includegraphics[width=\textwidth]{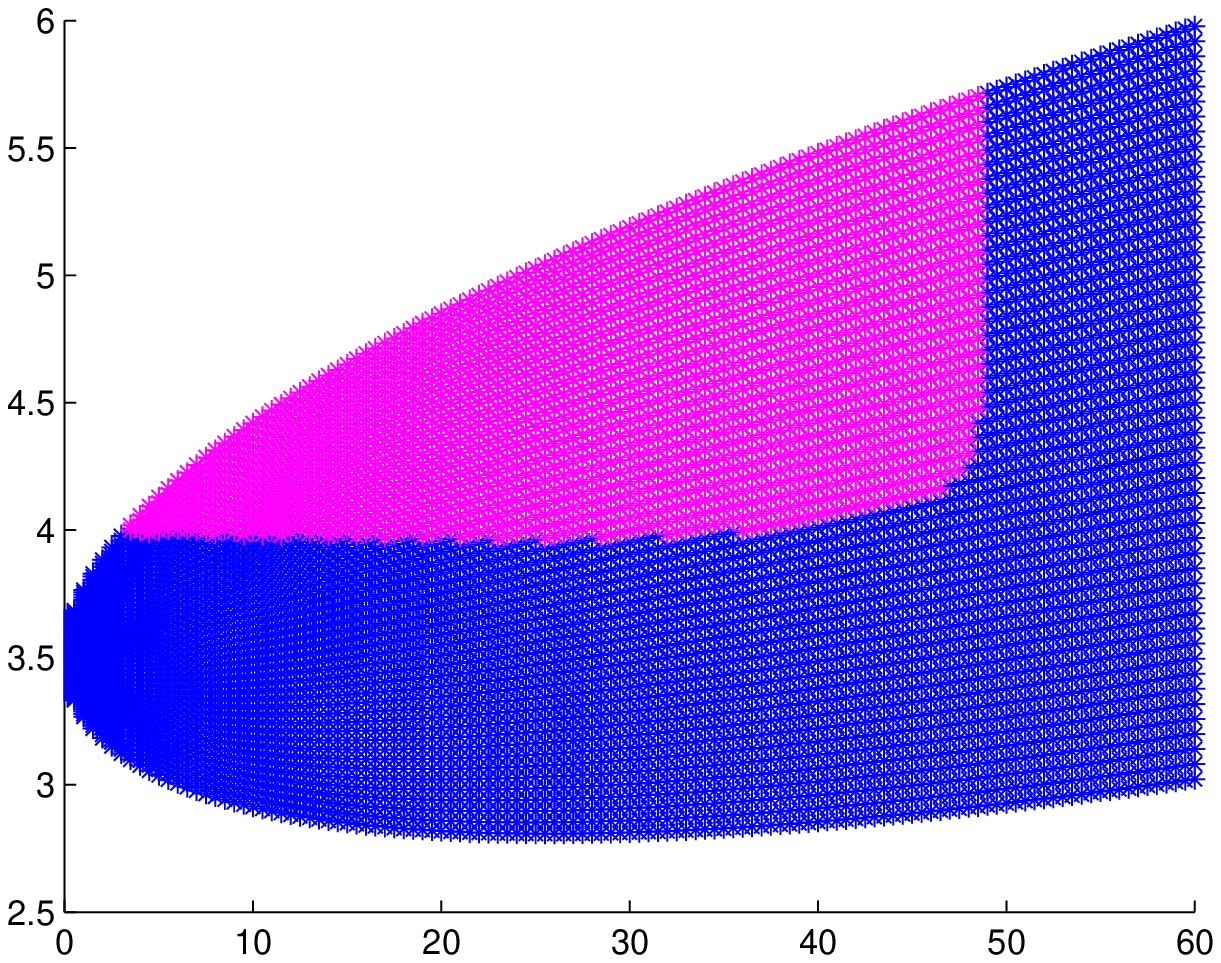}
                \caption{Switching regions from mode $[2\:\:  1]^\top$.}
                \label{fig:swRGS4}
        \end{subfigure}%
        \qquad
        \begin{subfigure}[b]{0.3\textwidth}
                \includegraphics[width=\textwidth]{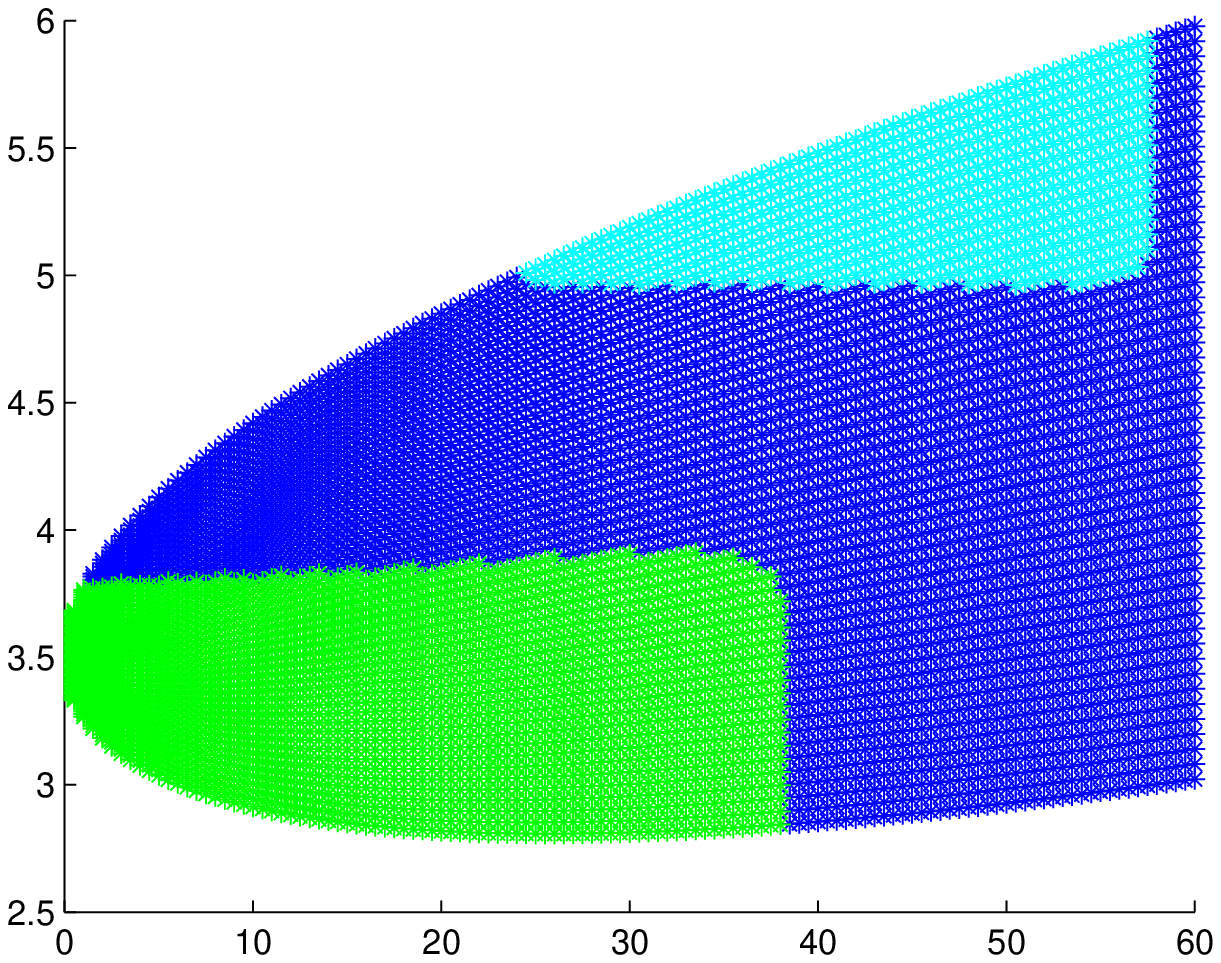}
                \caption{Switching regions from mode $[2\:\:  2]^\top$.}
                \label{fig:swRGS5}
        \end{subfigure}
        \qquad
        \begin{subfigure}[b]{0.3\textwidth}
                \includegraphics[width=\textwidth]{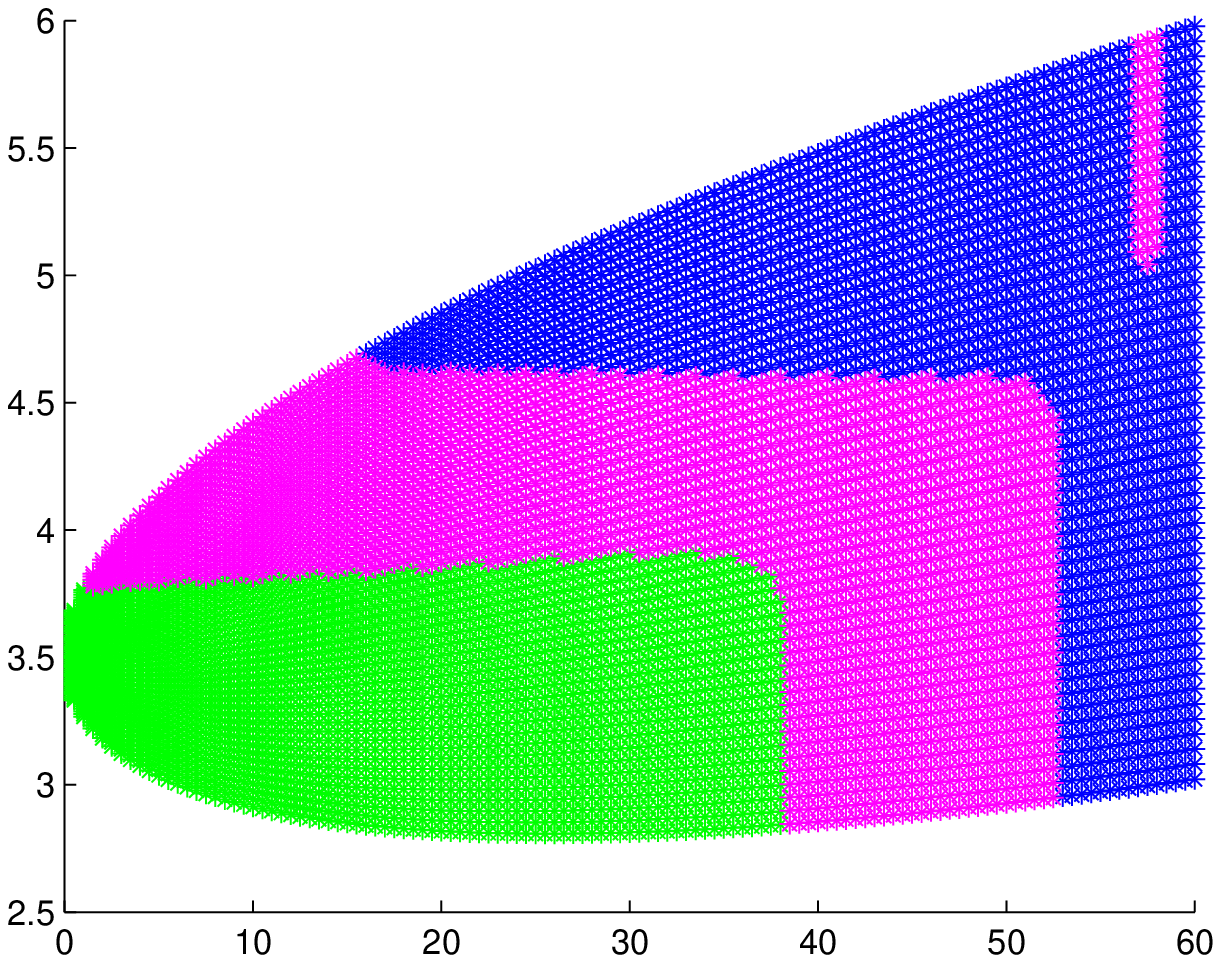}
                \caption{Switching regions from mode $[2\:\: 3]^\top$.}
                \label{fig:swRGS6}
        \end{subfigure}
        \caption{Switching regions for the stationary modes.}\label{fig:swRGS}
\end{figure}

We let $\Delta t = 0.5$ and solve the optimal switching problem using a Markov chain approximation of the demand process. The expected cost of operation is 6696~\currency~under the optimal control. The optimal actions in the stationary operation modes are plotted in Fig.~\ref{fig:swRGS}.

An explanation of the color scheme used in the figures is given in Table~\ref{tab:colorSCHEME}.
\begin{table}[h!]\centering \setlength{\tabcolsep}{10pt}
\begin{tabular}{| c | c | c | c | c | c | c |}
\hline
  {\color{blue}blue}  & {\color{yellow}yellow} & {\color{green}green} & {\color{red}red} & {\color{black}black} & {\color{magenta}magenta} & {\color{cyan}cyan} \\ \hline
  remain              &    $[1\:\: 1]^\top$    &   $[1\:\: 2]^\top$   & $[1\:\: 3]^\top$ &   $[2\:\: 1]^\top$   & $[2\:\: 2]^\top$ & $[2\:\: 3]^\top$ \\ \hline
\end{tabular}
\caption{Color and corresponding action (switching towards the stationary mode) in Fig~\ref{fig:swRGS}.}\label{tab:colorSCHEME}
\end{table}

\subsubsection{Upper and lower bounds on the value function\label{sect:LB_UB}}
When computing an approximation, $\hat V^{\text{LB}}$, of the lower bound we realize that with one load there is only one possible choice of $\eta_L$ ($\eta_L=1$). We choose the vectors $\eta_G$ such that each $(\eta_G^j,\eta_L)$ is normal to one of the feasibility boundaries (for $i=0,1,2$) at one of the production set-points. Hence, we get $N_{\eta_G}=6\cdot 3$. To approximate $X^\eta_t$ we use a Markov chain with a state-space of $N_x^\eta=1001$ points. The optimal cost-to-go at zero is estimated to 6471\currency (which is 225\currency~below the actual optimal cost-to-go at zero).

To compute an approximation, $\hat V^{\text{UB}}$, of the upper bound we approximate $X^B_t$ by a Markov chain with a state-space of $N_x^B=1001$ points. The upper bound is estimated to be 7672\currency~(\ie 976\currency~above the actual optimal cost-to-go).
\begin{figure}[h!]
  \centering
  \includegraphics[width=0.5\textwidth]{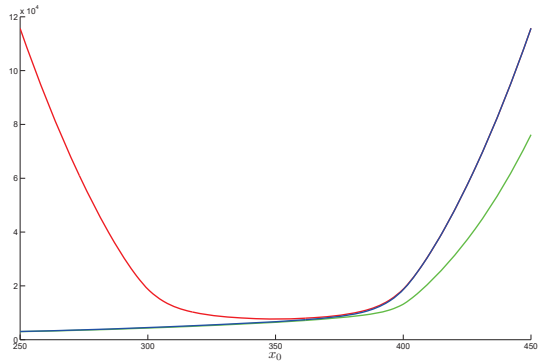}
  \caption{Value function at $t=0$ in $\alpha_0=[1\:\: 1]^\top$ and the lower and upper bounds.}\label{fig:LBexUB}
\end{figure}

\subsubsection{Regression Monte Carlo solution}
We apply the regression Monte Carlo method for two different settings:
\begin{enumerate}[A)]
  \item With the lower and upper bounds computed in Section~\ref{sect:LB_UB}
  \item With lower and upper bounds set to 0 and $\infty$ respectively
\end{enumerate}
We use a set of affine local basis functions with an adaptive state space partition. For the basis functions different partitions where used for different numbers of samples $M_T$. The basis used was decided based on an ad-hoc testing procedure where the number of hypercubes was increased until the numerical efficiency deteriorated due to over-fitting. Table~\ref{tab:locbaseEX1} gives a summary of the basis functions used in the example, with $M^a_{\max}=\max_{a\in\Gamma}M^a$.
\begin{table}[h!]\centering \setlength{\tabcolsep}{10pt}
\begin{tabular}{| c | c | c | c | c | c | c |}\hline
  & \multicolumn{3}{ c |}{A} & \multicolumn{3}{ c |}{B}\\
\hline
  $M_T$      &                  $I$                  &  $\tilde b^a$   & $M^a_{\max}$ &             $I$               &  $\tilde b^a$  & $M^a_{\max}$ \\  \hline
  1000       &  $[2\:\:2\:\:2\:\:1\:\:1\:\:1\:\:3]$  &  $\{1,r^a,x\}$  &   48   &   $[1\:\:1\:\:1\:\:1\:\:1\:\:1\:\:3]$  &  $\{1,r^a,x\}$  &   12    \\ \hline
  5000       &  $[2\:\:2\:\:3\:\:2\:\:1\:\:1\:\:8]$  &  $\{1,r^a,x\}$  &  192   &   $[2\:\:2\:\:2\:\:1\:\:1\:\:1\:\:5]$  &  $\{1,r^a,x\}$  &   80   \\ \hline
  10000      &  $[2\:\:2\:\:3\:\:2\:\:1\:\:1\:\:10]$ &  $\{1,r^a,x\}$  &  240   &   $[2\:\:2\:\:2\:\:1\:\:1\:\:1\:\:10]$ &  $\{1,r^a,x\}$  &  160   \\ \hline
  25000      &  $[2\:\:2\:\:3\:\:2\:\:1\:\:1\:\:15]$ &  $\{1,r^a,x\}$  &  360   &   $[2\:\:2\:\:3\:\:2\:\:1\:\:1\:\:10]$ &  $\{1,r^a,x\}$  &  240   \\ \hline
  50000      &  $[3\:\:2\:\:4\:\:2\:\:2\:\:1\:\:22]$ &  $\{1,r^a,x\}$  & 1056   &   $[2\:\:2\:\:3\:\:2\:\:1\:\:1\:\:15]$ &  $\{1,r^a,x\}$  &  360   \\ \hline
\end{tabular}
\caption{Estimated optimal costs in Example 1.}\label{tab:locbaseEX1}
\end{table}
The results, presented in Table~\ref{tab:CsEX1}, are based on 1000 different solutions by Monte Carlo regression using the mixture of policy approximation and value function approximation given in \eqref{ekv:mixAPPR}. The optimal controls $\hat v^*$ were approximated with $M_T$ ranging from 1000 to 50000. For each solution $J(\hat v^*)$ (the cost-to-go at time $t=0$ when applying the control $\hat v^*$) is computed by implementing $\hat v^*$ in a backward-scheme where $(X_t:\,0\leq t\leq T)$ is approximated by a Markov chain with 1001 states. Table~\ref{tab:CsEX1} also gives the estimated standard deviation of $J(\hat v^*)$.
\begin{table}[h!]\centering \setlength{\tabcolsep}{10pt}
\begin{tabular}{| c | c | c | c | c |}\hline
  & \multicolumn{2}{ c| }{A} & \multicolumn{2}{ c| }{B}\\
\hline
  $M_T$      & $ J(0,\hat v^*)$ & $\hat\sigma$ & $J(0,\hat v^*)$ & $\hat\sigma$ \\
  \hline
  1000       &      6818        &     15.04    &       7468      &     73.15    \\ \hline
  5000       &      6710        &     3.99     &       7135      &     21.86    \\ \hline
  10000      &      6708        &     3.28     &       7047      &      8.71    \\ \hline
  25000      &      6707        &     2.37     &       7010      &      7.47    \\ \hline
  50000      &      6706        &     1.40     &       6997      &      6.74    \\ \hline
\end{tabular}
\caption{Estimated optimal costs in Example 1.}\label{tab:CsEX1}
\end{table}


\subsection{Example on the IEEE 39--node system}
We now turn to an example in a system of more realistic size, namely the IEEE 39--bus New England system. All power system computations performed in this part are done using the MATLAB\textsuperscript{\textregistered}-package MATPOWER~\cite{MATPOWERref}.

\subsubsection{Network model}
\begin{figure}[h!]
  \centering
  \includegraphics[width=0.5\textwidth]{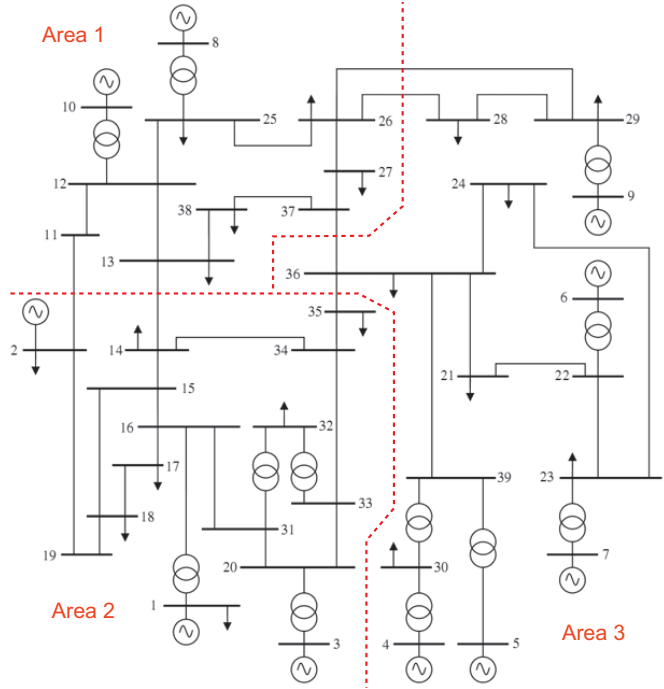}
  \caption{One-line diagram of the IEEE 39-bus system, with a three-area partition.}\label{fig:ieee39}
\end{figure}

The IEEE 39--bus system, depicted in Fig.~\ref{fig:ieee39}, is a model of the New England power system where the generator and load at Node 2, represent the connection to the rest of The Eastern Interconnection.
The system is divided into three areas based on geographical location.

\subsubsection{Feasible sets and contingencies}
We consider 41 different contingencies, 33 line contingencies with one contingency for each line, except the lines connecting the generators at nodes 1-10 to the rest of the systems and the lines connecting buses 36 to 39 and 30 to 39, and 8 generator contingencies with one for each generator except Generator 1 where the primary control is located, and Generator 2.

When a line contingency occurs we assume that the affected line is disconnected and when a contingency on a generator occurs we assume that the generator is taken out of service.

We set the failure weight $q_i=1e-4$, for $i=1,\ldots,33$, \ie for the lines, and $q_i=5e-5$, for $i=34,\ldots,41$, \ie for the generators.

\begin{figure}[h!]
  \centering
  \includegraphics[width=0.5\textwidth]{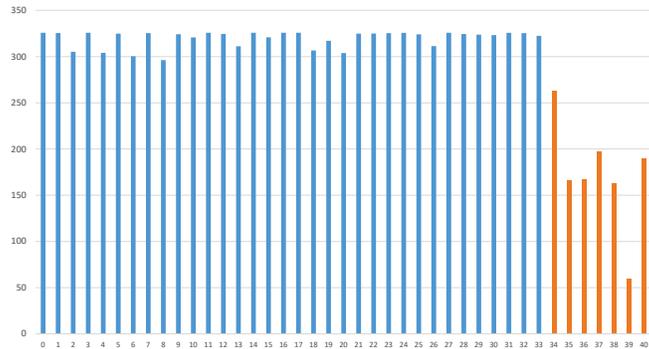}
  \caption{Loadability limit from the base case operating point for the different contingencies.}\label{fig:ieee39CCP}
\end{figure}

In Fig.~\ref{fig:ieee39CCP} the load-space distance from the base-case\cite{Pai} operating point to the feasibility boundary $\partial G_i$ is plotted for $i=0,\ldots,n_c$.
As we can tell from the figure the generator contingencies are more severe than the line contingencies.

\subsubsection{Stochastic model}
We assume that the stochastically varying loads are all major loads\footnote{We define a major load as a load with a default active power demand of at least 50 MW.} except the load at node 2 which is assumed to stay at the default value\footnote{This would correspond to a constant transmission of power to the other subsystems of the Eastern interconnection.}, which gives us 16 variable loads. We assume that these are also the curtailable loads and, thus, get $n_L=16$.

We assume that the active power demand is the sum of the forecasted trajectory $m(t)=p_D^0(1+0.3t/T)$, where $p_D^0$ is the base-case demand as given in \cite{Pai}. We set $\gamma=0.01$ and let $\sigma$ be such that the standard deviation at time $T$ is $0.1p_D^0$ and loads in the same area have correlation coefficient $0.1$, while loads in different areas are uncorrelated.

\subsubsection{Production modes and ramp rates}
We assume that the operator can switch between two production levels in Generator 2, representing the possibility of trading over the  transmission corridor that connects New England to the rest of the system. The generator at Node 10 has a large reserve capacity (base-case production 250MW compared to an installed capacity of 1040 MW) and we assume that this interval contains three levels of locally optimal efficiency which the operator can order switches between. Furthermore, we assume that there are two identical fossil-fuel production units at Node 4 only one of which is operating in base case.  The stationary production levels are given in Table~\ref{tab:thetaEX2}.

\begin{table}[h!]\centering \setlength{\tabcolsep}{10pt}
\begin{tabular}{| c | c | c | c | c | c | c |}\hline
  $\theta_{2,1}$ & $\theta_{2,2}$ & $\theta_{4,1}$ & $\theta_{4,2}$ & $\theta_{10,1}$ & $\theta_{10,2}$ & $\theta_{10,3}$ \\
  \hline
  1000 & 1500 & 508 & 1016 & 250 & 750 & 1040 \\ \hline
\end{tabular}
\caption{Stationary production levels [MW] in Example 2.}\label{tab:thetaEX2}
\end{table}
For the ramp function we assume, as in the previous example, that \eqref{ekv:RFex1} holds where the delay times for this example are given in Table~\ref{tab:deltasEX2}.
\begin{table}[h!]\centering \setlength{\tabcolsep}{5pt}
\begin{tabular}{| c | c | c | c | c | c | c | c | c |}\hline
  $k,(l,m)$           & $2,(1,2)$ & $2,(2,1)$ & $4,(1,2)$ & $4,(2,1)$ & $10,(1,2)$ & $10,(2,1)$ & $10,(2,3)$ & $10,(3,2)$ \\ \hline
  $\delta'_{k,(l,m)}$ & 5         & 5         & 15        & 3         & 2          & 2          & 2          & 2 \\
  $\delta_{k,(l,m)}$  & 10        & 7.5       & 20        & 5         & 5          & 4          & 4          & 3 \\ \hline
\end{tabular}
\caption{Delay times in minutes for Example 2.}\label{tab:deltasEX2}
\end{table}

At the start of the operating period we assume that we are in the stationary mode $\alpha_0=[1\:\: 1\:\: 1]^\top$.

\subsubsection{Operating costs\label{Sec:OC39}}
The operating cost is the sum of the production cost and the expected cost of disconnecting load. We assume that the production cost is $40\,$\currency/MWh~in Generator 2, $80\,$\currency/MWh~in Generator 4 and $30\,$\currency/MWh~in Generator 10. The switching costs are given in Table~\ref{tab:CsEX2}.

\begin{table}[h!]\centering \setlength{\tabcolsep}{5pt}
\begin{tabular}{| c | c | c | c | c | c | c | c | c |}\hline
  $k,(l,m)$                  & $2,(1,2)$ & $2,(2,1)$ & $4,(1,2)$ & $4,(2,1)$ & $10,(1,2)$ & $10,(2,1)$ & $10,(2,3)$ & $10,(3,2)$ \\ \hline
  $K_{k,(l,m)}\,$[\currency] & 2000      & 1000      & 10000      & 2000      & 200        & 100        & 300        & 100        \\ \hline
\end{tabular}
\caption{Switching costs in Example 2.}\label{tab:CsEX2}
\end{table}

The cost of disconnecting load is assumed to be the same in each node. We will compare the result of applying three different load disruption costs $c_j=1\,$k\currency/MW, $c_j=10 \,$k\currency/MW and $c_j=100\,$k\currency/MW, for $j=1,\ldots,n_L$.

\subsubsection{Dimensions}
There are 4 operating modes each for generators 2 and 4, and 7 modes for Generator 10, giving a total of $4^2\cdot 7=112$ operating modes. Of these only $2\cdot 2 \cdot 3=12$ are stationary modes.

In the modes where three generators are in transition simultaneously there are $\mathcal{O} (\Delta t^{-3})$ possible values for the transition progression vectors, giving a total computational time that is $\mathcal{O} ((N_\Pi)^{4})$. Combined with the 16-dimensional state space of the demand-process we realize that this problem is far out of reach of classical numerical methods.

\subsubsection{Upper and lower bounds on the value function\label{sect:LB_UB39}}
To compute an approximation of the lower bound we compute the normals to $\partial G_i$ for the non-contingent system ($i=0$) in the demand vectors closest to the base case for the production levels corresponding to the stationary modes. We thus let $p_L^{l}$ solve
\begin{align*}
\min\limits_{p_L\in \R^{n_L}_+}  \quad\quad & \|p_L-p_D^0\|\\
\text{subj. to}\quad\:\: & (\zeta(a^l,0),p_L) \in \partial G_0,
\end{align*}
where $a^l$ is an enumeration of the elements of $\Theta$. We then let $U=\{\eta_L^1,\ldots,\eta_L^{12}\}$,  where
\begin{equation*}
\eta_L^l=\frac{p_L^l-p_D^0}{\|p_L^l-p_D^0\|}.
\end{equation*}
We choose the vectors $\eta_G^{l,j}$ such that each $(\eta_G^{l,j},\eta_L^l)$ is normal to $\partial G_0$ at $(\zeta(a^l,0),p_D)$ for some $p_D\in\R^{n_L}$. This gives a total of $N_{\eta_G}\cdot |U|=12\cdot 12=144$ different vectors $\eta$ that have to be computed for the lower bound. To approximate $X^\eta_t$ we use a Markov chain with a state-space of $N_x^\eta=1001$ points.

To compute an approximation, $\hat V^{\text{UB}}$, of the upper bound we compute $d_i(t,p_G)$, for all $(t,p_G)\in \{0,6,12,\ldots,T\}\times \cup_{a\in\Theta}\zeta(a,0)$ and use linear interpolation to extend to $\Pi$ and then \eqref{ekv:dBext} to extend to $D_G$. We approximate $X^B_t$ by a Markov chain with a state-space of $N_x^B=1001$ points.

The results for the different disruption costs are presented in Table~\ref{tab:bndsEX2}.
\begin{table}[h!]\centering \setlength{\tabcolsep}{5pt}
\begin{tabular}{| c | c | c | c |}\hline
                & $c=1\,$k\currency/MW  & $c=10\,$k\currency/MW & $c=100\,$k\currency/MW \\ \hline
  Lower bound   &         6963          &          18434        &          22786         \\ \hline
  Upper bound   &        22539          &          44079        &          94865         \\ \hline
\end{tabular}
\caption{Lower and upper bounds on the value function.}\label{tab:bndsEX2}
\end{table}

\subsubsection{Approximating the expected cost of energy not served}
To represent the expected cost from energy not served, \ie $\sum_{i=0}^{n_c} q_i c^\top z^*_i(p_G,p_D)$, we build an ANN based on a training set of 10 000 randomly chosen points in the parameter space where the expected energy not served is given by the solution to the $n_c+1$ optimal power flow problems corresponding to the different system configurations.

\subsubsection{Regression Monte Carlo solution}
As in the previous example we will try two different settings:
\begin{enumerate}[A)]
  \item With the lower and upper bounds computed in Section~\ref{sect:LB_UB39}
  \item With lower and upper bounds set to 0 and $\infty$ respectively
\end{enumerate}
The operators problem will be solved for the three different disruption costs listed in Section~\ref{Sec:OC39}.

When making a state space partition for the 39-bus system, it would not be tractable to base the partition on an interval partition of each dimension. With $n_L=16$ an interval partition with 3 intervals for each load-space dimension would give $3^{16}\approx 4.3 e7$ hypercubes, presumably requiring billions of trajectories for convergence.

Instead we apply a state-space partition based on the areas outlined in Fig.~\ref{fig:ieee39}. Table~\ref{tab:locbaseEX1} gives a summary of the basis functions used in the example, with $M^a_{\max}=\max_{a\in\Gamma}M^a$.
\begin{table}[h!]\centering \setlength{\tabcolsep}{10pt}
\begin{tabular}{| c | c | c | c | c |}
\hline
  $M_T$      &       $I^D$       &                 $I^R$                 &   $\tilde b^a$     & $M^a_{\max}$ \\  \hline
  5000       &  $[2\:\:2\:\:2]$  &  $[2\:\:2\:\:2\:\:2\:\:2\:\:2\:\:2]$  &   $\{1,x^A\}$      &       256    \\ \hline
  10000      &  $[2\:\:2\:\:2]$  &  $[2\:\:2\:\:2\:\:2\:\:2\:\:2\:\:2]$  &   $\{1,x^A\}$      &       256    \\ \hline
  50000      &  $[2\:\:2\:\:2]$  &  $[2\:\:2\:\:2\:\:2\:\:2\:\:2\:\:2]$  &   $\{1,x\}$        &      1088    \\ \hline
  100000     &  $[2\:\:2\:\:2]$  &  $[2\:\:2\:\:2\:\:2\:\:2\:\:2\:\:2]$  &   $\{1,x\}$        &      1088    \\ \hline
\end{tabular}
\caption{Estimated optimal costs in Example 2.}\label{tab:locbaseEX1}
\end{table}
In the table $x^A\in\R^{n_A}$ refers to the vector giving the total demand in each of the $n_A$ areas, \ie $x_k^A=\sum_{j\in A_k}x_j$.\\

To evaluate the approximated strategies a Monte Carlo simulation, using the lower bound as a control variable, is performed with one million cases. The results are given in Table~\ref{tab:resEX2}.
\begin{table}[h!]\centering \setlength{\tabcolsep}{10pt}
\begin{tabular}{| c | c | c | c | c | c | c |}\hline
  & \multicolumn{2}{ c| }{$c=1\,$k\currency/MW} & \multicolumn{2}{ c| }{$c=10\,$k\currency/MW} & \multicolumn{2}{ c| }{$c=100\,$k\currency/MW} \\
\hline
  $M_T$      &    A    &    B    &    A    &    B    &    A    &    B    \\
  \hline
  5000       &      8023        &     15979     &     19193        &      34286       &     36279       &     55564       \\ \hline
  10000      &      7893        &     14731     &     18967        &      29784       &     34997       &     54018       \\ \hline
  50000      &      7555        &     12519     &     18970        &      28592       &     32525       &     50088       \\ \hline
 100000      &      7472        &     11836     &     18958        &      27253       &     32312       &     48939       \\ \hline
\end{tabular}
\caption{Estimated optimal costs in Example 2.}\label{tab:resEX2}
\end{table}

\subsection{Discussion}
In both of the above examples we see a significant improvement when using the proposed method compared to crude regression Monte Carlo. The origin of the improvement is probably to a large extent reduction of variance in the estimate. Another possible factor is that it is easier to
fit the local polynomials to $\Delta V^\eta$ than to $V$, as $V$ may be highly non-linear.

The higher dimensionality of the second problem is reflected in a slower convergence, requiring a larger number of Monte Carlo samples to obtain an acceptable accuracy.

The proposed method may lead to tractable computations for more realistic systems, however, an alarming fact is that the numerical complexity is linear the cardinality of $\Gamma$, which grows exponentially in the number of controllable generators. This seriously threatens the usefulness of the method as many producers will in general be interested in participating in the real-time market (where regulating power is sold). To manage situations with many controllable generators we may have to do some sort of screening of operating modes before attempting to solve the problem with regression Monte Carlo. A possible approach would be to use the upper and lower bounds on the value function to sort out operating modes that cannot be optimal.

In the above examples a tighter upper bound for the value function at $t=0$ could be estimated by simulating trajectories and applying the control induced by the lower bound. An interesting line of future research would be to determine if the estimated optimal control can be further improved by using the upper and lower bounds when generating the samples of $(R_{t_m};\: t_m\in \Pi)$.

\bibliographystyle{plain}
\bibliography{pFORM_ref}
\end{document}